\documentclass[12pt,reqno]{amsart}
\usepackage{verbatim}
\usepackage{stmaryrd}
\usepackage{amsmath}
\usepackage[mathscr]{eucal}
\newenvironment{pf}{\noindent{\it Proof.} }{\fbox{}\\}
\newtheorem{thm}{Theorem}[subsection]
\newtheorem{cor}{Corollary}[subsection]
\newtheorem{prop}{Proposition}[subsection]
\newtheorem{lem}{Lemma}[subsection]
\newtheorem{defi}{Definition}[subsection]
\newtheorem{rem}{Remark}[subsection]
\newcommand{\1}{\mathbf{1}}
\title{ White Noise Calculus and Hamiltonian of a Quantum Stochastic Process}
\author{Wilhelm von Waldenfels}
			
\address{Institut f\"ur Angewandte Mathematik\\Universit\"at
Heidelberg\\ Im Neuenheimer Feld 294 \\ 69120 Heidelberg \\ Germany \\ 
 e-mail Wilhelm.Waldenfels@T-Online.de}
\begin{document}
\begin{abstract}A white noise quantum stochastic calculus is developped using classical measure theory as mathematical tool.
Wick's and Ito's theorems have been established. The simplest quantum stochastic differential equation has been solved,
unicity and the conditions for unitarity have been proven. The Hamiltonian of the associated one parameter strongly
continuous group has been calculated explicitely.
\end{abstract} 
\maketitle
\tableofcontents

\section*{Introduction}
The object of this paper is at first to establish a white noise quantum stochastic calculus 
using  classical measure theory as mathematical tool. Secondly we solve the simplest quantum stochastic 
equation written in normal ordered form \cite{acc2002}

\begin{align}
\frac{d}{dt}U^t_s &= A_1 a^+_t U^t_s + A_0 a^+_t U^t_s a_t + A_{-1} U^t_s a_t + B U^t_s \\
U^s_s &= 1 \nonumber
\end{align}
and show the unitarity of the solution. Thirdly we derive the Hamiltonian of the associated one-parameter 
strongly continuous group in the explicite form
\begin{equation} H = i\hat{\partial} + M_1 \hat{\mathfrak{a}}^+ +M_0 \hat{\mathfrak{a}}^+\hat{\mathfrak{a}}+
 M_{-1} \hat{\mathfrak{a}} + G \end{equation}
 where
\begin{align}A_1 &= \frac{1}{i-M_0/2}M_1\nonumber\\
A_0&=\frac{M_0}{i-M_0/2}\nonumber\\
A_{-1}&= M_{-1}\frac{1}{i-M_0/2}\nonumber\\
B&= -iG - \tfrac{i}{2}M_{-1}\frac{1}{i-M_0/2}M_1
\end{align}
Quantum stochastic calculus can be considered as a kind of one dimensional field theory. The basic quantities
are the creation operator $a_t^+$ and the annihilation operator $a_t$ of a particle at $t \in \mathbb{R}$.In the
context of quantum probability these operators are called  quantum white noise operators.  Treating them
in a mathematical rigorous way is not straight forward. Accardi et al.\cite{acc99} use Schwartz's distribution theory,
Obata et al.\cite{obata08}use the theory of Gelfand triples. We apply just classical measure theory exploiting the positivity
of the white noise operators. Equation (1) has been stated and solved in an other form at first by
Hudson and Parthasarathy \cite{partha}. The kernel formalism developped by Maassen and Meyer \cite{meyer}\cite{vw96}
\cite{vw2000}.
 is dual to white noise calculus.
There the equation can be solved too.

Equation (3) occurs in what Accardi calls the Hamiltonian form of the quantum stochasic differential equation
\cite{acc96}\cite{gough}\cite {gough97}\cite{Ayed2005}. The Hamiltonian form of a quantum stochastic differential 
equation is a related but different 
object to the Hamiltonian considered here. Neither the obvious connection nor Accardi's time consecutive
principle have  been regarded in this paper.

Denote by $\Theta(t)$ the  time shift by $t$, then
\begin{align}W(t)=&\Theta(t)U_0^t;& W(-t)&=W(t)^+;&t& \in \mathbb{R}\end{align}
forms a unitary, strongly continuous one parameter semigroup on the Fock space. By
Stone's theorem there exists a selfadjoint operator $H$, called the Hamiltonian, such
that
\begin{equation} W(t)=e^{-iHt}.\end{equation}
So the existence of $H$ is clear, what was not clear,was an explicit
analytic expression representing $H$. The problem was posed by Accardi \cite{acc89}in 1989,
solved for the case of commuting $L$ and $L^+$ by Chebotarev \cite{chebotarev} in 1997 and
finally solved by Gregoratti in 2000 in his thesis \cite{gregoratti}. I presented in a previous paper for the case
$A_0$ a 
different approach using Maassen-Meyer calculus and obtained Gregoratti's result with a different proof and in a different,
more intuitive representation. Here I treat the general case, use white noise calculus, calculate
the resolvent and arrive at formula (2). Our result is again equivalent to that of Gregoratti.

The basic equations for the white noise operators are
\begin{align} [a_s,a^+_t] &= \delta(s-t)\\a_s\Phi&=0\nonumber ,
\end{align}
where $\Phi$ is the vacuum. Let us as first approximate the creation and annihilation operators. Assume a continuous
function $\varphi$ of compact support in $\mathbb{R}$ and a continuous \emph{symmetric}  function $f$ on $\mathbb{R}^n$
then
\begin{align*} 
(a^+(\varphi)f)(t_0,t_1,\cdots ,t_n) =& \varphi(t_0)f(t_1,\cdots ,t_n)+\varphi(t_1)f(t_0,t_2,\cdots ,t_n)+\\&\cdots +
\varphi(t_n)f(t_0,t_1,\cdots ,t_{n-1})\\
 (a(\varphi)f)(t_2,\cdots ,t_n)=& \int_{\mathbb{R}} dt_1\varphi(t_1)f(t_1,\cdots ,t_n).\end{align*}
Shrink $\varphi = \varphi(t)$ to $\delta(t-s)$ and obtain
\begin{align*}(a^+(s)f)(t_0,t_1,\cdots ,t_n)
= &\delta(t_0-s)f(t_1,\cdots ,t_n)+ \\&\cdots +\delta(t_n-s)f(t_0,\cdots ,t_{n-1})\\
 a(s)f(t_2,\cdots ,t_n) =& \int dt_1 \delta(t_1-s) f(t_1,\cdots ,t_n)= f(s,t_2,\cdots ,t_n)\end{align*} 
The $\delta$-function can interpreted in three different ways. Denote by $\varepsilon_s(dt)$ the point measure at the
point $s \in \mathbb{R}$ 
and by $\varepsilon(dt)$ the function $s \mapsto \varepsilon_s(dt)$. Then 
\begin{enumerate}
\item $\delta(t-s)ds = \varepsilon_t(ds)$
\item $ \delta(t-s)dt = \varepsilon_s(dt)$
\item $\delta(t-s)dtds = \varepsilon_s(dt)ds = \varepsilon_t(ds) dt = \lambda(dt,ds)$\\
with $\int \lambda(dt,ds)g(t,s)= \int dt g(t,t)$
\end{enumerate}

We use mostly the first possibility and arrive at
\begin{multline}
(a^+(ds)f)(t_0,\cdots ,t_n) = (a^+(\varepsilon(ds)f)(t_0,\cdots ,t_n)\\ = \varepsilon_{t_0}(ds) f(t_1,\cdots ,t_n) +\cdots +
\varepsilon_{t_n}(ds)f(t_0,\cdots ,t_{n-1})
\end{multline}
So the function $\varphi$ is replaced by the function $t \mapsto \varepsilon_t(ds)$. For the second possibility we use
a gothic $\mathfrak{a}$ in order to distinguish.
$$ (\mathfrak{a}^+(s)f)(t_0,\cdots ,t_n)= \varepsilon_s(dt_0)f(t_1,\cdots ,t_n)+\cdots 
+\varepsilon_s(dt_n)f(t_0,\cdots ,t_{n-1})$$
This expression has to be interpreted in the notation of L.Schwartz, who identifies measures having a density with
respect to Lebesgue measures with their densities. So in usual notation
\begin{multline} (\mathfrak{a}^+(s)f)(dt_0,\cdots ,dt_n)= \varepsilon_s(dt_0)f(t_1,\cdots ,t_n)dt_1\cdots dt_n
+\cdots \\
+\varepsilon_s(dt_n)f(t_0,\cdots ,t_{n-1})dt_0\cdots dt_{n-1}.\end{multline}
We shall use this notation only in the last section. The third possibility occurs in hidden form in calulations further on.

 For the annihilation operator there is no problem
 \begin{equation}(a(s)f)(t_2,\cdots ,t_n) = \int \varepsilon_s(dt_1) f(t_1,\cdots ,t_n)= f(s,t_2,\cdots ,t_n)
\end{equation}
We obtain the commutation relation
 \begin{equation} [a(s),a^+(dt)] = \varepsilon_s(dt). \end{equation}
Define the infinitesimal number operator
\begin{equation} \mathfrak{n}(ds) = a^+(ds)a(s), \end{equation}
then
$$(\mathfrak{n }(ds)f)(t_1,\cdots ,t_n) = \sum _{i=1}^n \varepsilon_{t_i}(ds)f(t_1,\cdots ,t_n).$$

We introduce the space
\begin{equation} \mathfrak{R} = \{\emptyset\} + \mathbb{R}+\mathbb{R}^2+\cdots  \end{equation}
where the $+$ -sign indicates disjoint union .We introduce in $\mathfrak{R}$ the Lebesgue measure 
\begin{equation} \int f(w)dw = f(\emptyset) +\sum _{n=1}^\infty
\frac{1}{n!} \int _{\mathbb{R}^n}dt_1\cdots dt_n f(t_1,\cdots ,t_n) \end{equation}
The \emph{Fockspace}
$\Gamma$ is the space of all square integrable \emph{symmetric} functions on $\mathfrak{R}$.

A normal ordered monomial of creation and annilation operators is of the form
\begin{multline} a^+(dt_{a_1})\cdots a^+(dt_{a_l})a^+(dt_{b_1})\cdots a^+(dt_{b_m})\\\times
a(dt_{b_1})\cdots a(dt_{b_m})a(t_{c_1})\cdots a(t_{c_n}) = a^+(t_\alpha)a^+(t_\beta)a(t_ \beta)a(t_\gamma)\end{multline}
with
$$ \alpha =\{a_1,\cdots ,a_l\},\beta=\{b_1,\cdots ,b_m\},\gamma = \{c_1,\cdots ,c_n\}$$

We make for (1) the Ansatz 
\begin{equation} U^t_s = \sum_{l,m,n =0}^\infty \frac{1}{l!m!n!}\iiint _{\alpha,\beta,\gamma}u^t_s(t_\alpha,t_\beta,t_\gamma)
a^+(dt_\alpha)a^+(dt_\beta)a(t_ \beta)a(t_\gamma)dt_\gamma\end{equation}
and assume that
$$ (t,t_\alpha,t_\beta,t_\gamma) \mapsto u^t_s(t_\alpha,t_\beta,t_\gamma) $$
is symmetric in $t_\alpha$ and in $t_\beta$ and in $t_\gamma$ and is locally integrable. Then there exists a unique 
solution which can be explicitely given as follows. Assume that all points $t,t_\alpha,t_\beta,t_\gamma$ are different
,denote $p=l+m+n$,
and order the set 
$$t_\alpha\cup t_\beta \cup t_\gamma  = \{s_1 < \cdots <s_p \}$$
in increasing order and define $i_j = 1,0,-1$ if $s_j \in t_\alpha,t_\beta,t_\gamma$ resp..Then
\begin{multline} u^t_s(t_\alpha,t_b,t_\gamma)= \1\{t_\alpha\cup t_\beta \cup t_\gamma \subset ]s,t[ \}\\\times
e^{(t-s_p)B}A_{i_p}e^{(s_p-s_{p-1})B}A_{i_{p-1}}\cdots A_{i_2}e^{(s_2-s_1)B }A_{i_1}e^{(s_1-s)B}\end{multline}
The function $u^t_s$ has remarkable analytical properties, e.g. it is smooth except a finite number of jumps.
We define a class of functions $\mathcal{C}^1$  which includes $u^t_s$ and has the same essential properties. 

In order to prove unitarity one has to deal with products of normal  ordered monomials. These are not normal
ordered any more, but they are \emph{admissible}: if $a^+(dt)$ and $a(t)$ occur in the same monomial then $a^+(dt)$
must be left of $a(t)$. The theory of admissible monomials can be done in any locally compact space. 
We improve some previous
studies \cite{vW2007} , prove Wick's theorem and establish duality.

Using that theory and the properties of $\mathcal{C}^1$ functions we are able to show a generalization of Ito's
 theorem. With its help we establish the  unitarity of $U^t_s$. Denote
by $\Gamma_k \subset \Gamma,\,k=0,1,2,\cdots$ the subspace of those functions for which
$$ \langle f |(N+1)^k|f \rangle = \sum_{n=0}^\infty \frac{(n+1)^k}{n!} \int \|f(t_1,\cdots,t_n)\|^2
 dt_1 \cdots dt_n
< \infty .$$
 Then $\Gamma_k$ is invariant with respect to $U^t_s$.
 
 In order to describe the Hamiltonian we split the point $0 \in \mathbb{R}$ into two points and introduce 
 $\mathbb{R}_0 = ]-\infty,-0]+[+0,\infty[$ with the usual topology and Lebesgue measure and define 
 $\mathfrak{R}_0$ accordingly. We define the measures
 $\varepsilon_{\pm 0}$ and 
 \begin{align}(\mathfrak{a}_{\pm 0}f)(t_2,\cdots ,t_n)& = f(\pm 0,t_1,\cdots ,t_n)\\
  (\mathfrak{a}_{\pm 0}^+f)(t_0,\cdots ,t_n) &= \varepsilon_{\pm 0}(dt_0)f(t_1,\cdots ,t_n)+\cdots  
 + \varepsilon_{\pm 0}(dt_n)f(t_0,\cdots ,t_{n-1})\nonumber\end{align}
 We define
 \begin{align} \hat \varepsilon_0 =& \tfrac {1}{2}(\varepsilon_{+0}+\varepsilon_{-0})&
\hat{\mathfrak{a}} = &\tfrac {1}{2}(\mathfrak{a}_+ +\mathfrak{a}_-)&
\hat{\mathfrak{a}}^+ = &\tfrac {1}{2}(\mathfrak{a}^+_+ +\mathfrak{a}^+_-).
\end{align}
 
The symmetric derivative $\hat{\partial}$ is defined for functions on $\mathfrak{R}_0$ \cite{vW2005a}.
 It is the  second 
quantization of its restriction to functions on $\mathbb{R}_0$, which are continuously differentiable outside
0. Take as example the Heaviside function
$$Y(t)= \left\{\begin{array}{rcl}1&\mbox{ for }&t>0\\
0&\mbox{ for }&t<0 \end{array}\right.$$
and not defined for $t=0$.Then
$$ \hat{\partial}Y(t)= (Y(t+0)-Y(t-0))\hat{\varepsilon}_0(dt)=\hat{\varepsilon}_0(dt)$$
This notion is very much related to the notion of a $\delta$-
function on simplices, introduced by 
J.Gough\cite{gough} and by Accardi ,Lu and Volovich\cite {acc99} .

We introduce the preliminary Hamiltonian
\begin{equation}\hat {H} = i\hat{\partial} + M_1 \hat{\mathfrak{a}}^+ +M_0 \hat{\mathfrak{a}}^+\hat{\mathfrak{a}}+
 M_{-1} \hat{\mathfrak{a}} + G .
 \end{equation}
 with
 \begin{align}M_1&=M_{-1}^+&M_0&=M_0^+&G&=G^+.\end{align}
We define
$$Z= \int _0^\infty e^{-t} \Theta(t) dt ,$$
where $\Theta(t)$ is the time shift.The domain $D$ of $\hat{H}$ is the space of all functions of the form
$$f=Z(f_0+\mathfrak{a}^+f_1)$$
with
$$(\mathfrak{a}^+f)(t_0,\cdots ,t_n)=
 \varepsilon_0(dt_0)f(t_1,\cdots ,t_n)+\cdots 
+\varepsilon_0(dt_n)f(t_0,\cdots ,t_{n-1})$$
using again Schwartz's convention
and with $f_0 \in \Gamma_1,\,f_1\in \Gamma_2$. So $\hat{H}$ maps $D$ into a space of measures on $\mathfrak{R}_0$,
which might be singular at the origine $\pm 0$.

 Denote by $D_0 \subset D $ the subspace where $\hat{H}$ is non-singular,
then $D_0 \subset D_H$, where $D_H$ is the
domain of $H$ and $D_0$ is dense in $\Gamma$.  On $D_0$ the operator $H$ equals the restriction of $\hat{H}$ to 
$D_0$ and $H$
 is the closure of the restriction.
 
 Gregoratti uses as well the operators $\mathfrak{a}_\pm$ but he does not employ the singular operators
 $\mathfrak{a}_\pm^+$. Therefor his expression for $H$ is not symmetric. To prove symmetry, he has to impose
 boundary  conditions. These turn out to be equivalent to our restriction to $H_0$. We base  the
 characterization of the Hamiltonian on the resolvent, where as Gregoratti uses a direct approach.

\section{Admissible monomials}
\subsection{Some notations}

Let $X$ be a set . We denote by $\mathfrak{W}(X)$ the set of all finite
sequences of elements of $X$ or words formed by elements of $X$ and write
for short $\mathfrak{X} = \mathfrak{W}(X)$ and we have
$$ \mathfrak{X}=\emptyset+ X + X^2+ \cdots $$

We use the plus sign for denoting
disjoint union of sets. We denote by $\mathfrak{S}(X)$ the set of finite 
subsets of $X$ and by $\mathfrak{M}(X)$ the set of finite multisets of $X$.
A multiset is on $X$ is a mapping $m:X \rightarrow \mathbb{N}= \{ 0,1,3,\cdots \}$
The cardinality of $m$ is $\sharp m = \sum _{x \in X} m(x)$.. The multiset is
 finite if its cardinality is finite. We may write multisets in the form of an enumaration
 $ m = \{x_1,\cdots ,x_n\}^\bullet$, where $x \in X $ occurs $m(x)$ times. So not all $x_i$ must be different.
If $w = (x_1,\cdots,x_n)$ is a word and $\sigma$ is a 
permutation,then
$$ \sigma w = ( x_{\sigma^{-1}(1)}, \cdots ,  x_{\sigma^{-1}(n)})$$
The word $w$ defines a multiset 
$$m_w = \{x_1, \cdots, x_n\}^\bullet$$
If $w'$ is another word, then $m_{w'}= m_w$ iff there exists a permutation changing $w$ 
into $w'$.

An $n$-chain is a totally ordered set of $n$ elements 
$$\underline{a}=( a_1, \cdots , a_n ) $$ 

We denote by $ x_{\underline{a}}$ the word 
$$x_{\underline{a}}=( x_{a_1}, \cdots , x_{a_n} ) $$ 
If $\alpha$ is the underlying set of $ \underline{a}$, we denote by
$$x_\alpha=\{ x_{a_1}, \cdots , x_{a_n} \}^\bullet $$ 
 the corresponding multiset, as it does not depend on the order of $ \underline{a}$. That is the usual way we are writing multisets.

A function $f$ on $\mathfrak{X}$ is called symmetric, if $f(w) = f( \sigma w )$ for all
 permutations of $w$. If $\alpha$ is an $n$-set, i.e. a set with n elements , then $f(x_{\alpha})$ is well defined.

A function $f$ on $\mathfrak{X}$ is called symmetric, if $f(w) = f( \sigma w )$ for all
 permutations of $w$. If $\alpha$ is an $n$-set, i.e. a set with n elements , then $f(x_{\alpha})$ is well defined.

In order to avoid unnecessary complications, we shall only consider locally compact spaces countable at infinity.
Assume now, that $X$ is a locally compact space, provide $X^n$ with the
product topology and $\mathfrak{X} = \mathfrak{W}(X)$ with that topology where 
the $X^n$ are as well open and closed and where the restriction to $X^n$
coincides with the topology of $X^n$. Then $\mathfrak{X}$ is locally compact as well,
 its compact sets are 
contained in  a finite union of the $X^n$ and their intersection with the $X^n$ are compact.

If $S$ is a locally compact space, denote by $\mathcal{K}(S)$ the space of complex valued 
continuous functions on $S$ with compact support and by $\mathcal{M}(S)$ the space of complex
 measures on $S$. If $\mu$ is a complex measure on $\mathfrak{X}$ , we write
$$\mu=\mu_0 +  \mu_1 + \mu_2 + \cdots $$
\noindent
where $\mu_n$ is the restriction of $\mu$ to $X^n$. We denote by $\Psi$ the measure given
by 
\begin{equation*}\label{eq1}\Psi(f) = f(\emptyset).\end{equation*}
 Then $\mu_0$ is a multiple of $\Psi$. If $\underline{a}$ is an
 $n$-chain, we denote
$$\mu (dx_{\underline{a}})=\mu_n(dx_{a_1}, \cdots dx_{a_n})$$

If $\mu$ is a symmetric measure, then $\mu(d x_{\alpha})$, where $\alpha$ is an $n$-set
is well defined. 

We define
$$\int_{\Delta(\mathfrak{X})} \mu(dw) f(w) = \mu_0f(\emptyset)
 + \sum_{n=1}^\infty \frac{1}{n!}\int \mu_n(dx_1,\cdots,dx_n)f(x_1,\cdots,x_n) $$
 Instead of $\int_{\Delta(\mathfrak{X})}$ we write often $\int_\Delta$ ore simply $\int$.
 Similar 
 $$\int_{\Delta(X^n)} \mu (dx_1,\cdots ,dx_n)f(x_1,\cdots ,x_n)= \frac {1}{n!} \int _{X^n}
 \mu (dx_1,\cdots ,dx_n)f(x_1,\cdots ,x_n)$$ 
 If $X=\mathbb{R}$ and
 \begin{equation*} \mathfrak{R} = \{\emptyset\} + \mathbb{R} + \mathbb{R}^2 + \cdots  \end{equation*}
 and if $\mu$ is the Lebesge measure
 $$\int_{\Delta(\mathbb{R}^n)} \mu dx_1\cdots dx_n f(x_1,\cdots ,x_n)= \int _{x_1<\cdots <x_n}
 dx_1\cdots dx_nf(x_1,\cdots ,x_n)$$
 
A \emph{hierarchy} is a family of finite index sets
$$ A= \{\alpha_0,\alpha_1,\alpha_2,\cdots\},$$
such that $\#\alpha_n = n$. We write 
\begin{eqnarray*} \int_{\Delta(\mathfrak{X})} \mu(dw) f(w)&=& \int _{\Delta,\alpha \in A}\mu(dx_\alpha) f(x_\alpha)\\
&=& \mu_0f(\emptyset)
 + \sum_{n=1}^\infty \frac{1}{n!}\int \mu_n (dx_{\alpha_n})f(x_{\alpha_n})\end{eqnarray*}
or simply
$$\int _{\Delta,\alpha }\mu(dx_\alpha) f(x_\alpha) $$
or
$$\int_\alpha \mu(\alpha)f(\alpha)$$
if the context and the variable $x$ is clear. 
\subsection{Sum Integral Lemma for Measures}
We recall the sum integral lemma .
We shall need a stronger result, a sum integral lemma for measures, as it was written in our papers \cite{vw2000}
\cite{vW2005a}
For the sake of completeness we shall reformulate and prove it.
\begin{lem}\label{lem2.1}[Sum-integral lemma for measures].Let be given a  measure
$$\mu(dw_1, \cdots , dw_k)$$
on 

$$\mathfrak{X}^k= \sum _{n_1,\cdots,n_k}X^{n_1}\times\cdots X^{n_k},$$ 
symmetric in each of the variables $w_i$ and assume that
$$\sum\frac{1}{n_1! \cdots n_k!}\, \mu_{n_1, \cdots, n_k} $$
is a bounded measure on $\mathfrak{X}^k$
where $  \mu_{n_1, \cdots, n_k} $
is the restriction of $\mu$ to $X^{n_1} \times\cdots \times X^{n_k}$. Then 
$$\int_{\Delta(\mathfrak{X})} \cdots \int_{\Delta(\mathfrak{X})} \mu(dw_1,\cdots,dw_k)=
\int_{\Delta(\mathfrak{X})}\nu(dw)$$
where $\nu$  is a  measure on  $\mathfrak{X}$ and 
$\sum (1/n!)\nu_n$ is a bounded measure,where
$\nu_n$ is the restriction of $\nu$ to $X^n$ and 
$$\nu_n(dx_1,\cdots, dx_n)=\sum_{\beta_1+\cdots+\beta_k=[1,n]}
\mu_{\#\beta_1, \cdots, \#\beta_k}(dx_{\beta_1}, \cdots, dx_{\beta_k}).$$ 
\noindent
where $\beta_1, \cdots, \beta_k$ are disjoint sets.
 
Using hierarchies $A_1,\cdots ,A_k,B$ we may write
$$\int_{\Delta,\alpha_1\in A_1}\cdots \int_{\Delta,\alpha_k\in A_k}\mu(dx_{\alpha_1},\cdots ,dx_{\alpha_k})
= \int_{\Delta,\beta \in B} \nu(dx_\beta)$$
With
\begin{eqnarray*}\nu(dx_\beta)&=& \sum_{\beta_1+\cdots +\beta_k=\beta}\mu(dx_{\beta_1},\cdots ,dx_{\beta_k})\\
\mu(dx_{\beta_1},\cdots ,dx_{\beta_k})&=& \mu_{\#\beta_1,\cdots ,\#\beta_k}(dx_{\beta_1},\cdots ,dx_{\beta_k}).
\end{eqnarray*}
\end{lem}
\begin{pf}
\begin{multline*}\int_{\Delta(\mathfrak{X})} \cdots \int_{\Delta(\mathfrak{X})} \mu(dw_1,\cdots,dw_k)\\
=\sum_{n_1,\cdots  ,n_k}\int_{X^{n_1}}\cdots \int_{X^{n_k}}
\frac{1}{n_1!\cdots n_k!}\mu_{n_1,\cdots ,n_k}(dx_{\alpha_1},\cdots ,dx_{\alpha_n})
\end{multline*}
where the $\alpha_i$ are the intervals
$$\alpha_1= [1,n_1],\alpha_2=[n_1+1,n_1+n_2],\cdots ,\alpha_k = [n_1+\cdots n_{k-1}+1,n_1+\cdots +n_k]$$
Fix $n_1,\cdots ,n_k$ and put $n=n_1+\cdots n_k$.Then
 Then
\begin{multline*}\int_{X^{n_1}} \cdots \int_{X^{n_k}}
\mu_{n_1,\cdots,n_k}(dx_{\alpha_1},\cdots,dx_{\alpha_k})\\ = \frac{1}{n!} \sum_\sigma
\int_{X^{n_1}} \cdots \int_{X^{n_k}}
\mu_{n_1,\cdots,n_k}(dx_{\sigma(\alpha_1)},\cdots,dx_{\sigma(\alpha_k)})\end{multline*}
where the sum runs through all permutations of $n$ elements. The subsets $\sigma(\alpha_i)
=\beta_i$ have the property
\[\beta_1+ \cdots \beta _k =[1,n], \#\beta_i=n_i\tag{$*$}.\]
Fix $\beta_1,\cdots \beta_k $ with property $(*)$.	
There are exactly $n_1!\cdots n_k!$ permutations $\sigma$, such that
$$\sigma(\alpha_i)=\beta_i \mbox{ for } i=1,\cdots,k.$$
Hence the last expression equals
$$\frac{n_1!\cdots n_k!}{n!} \sum_{\beta_1, \cdots ,\beta _k}
\int \cdots \int \mu_{n_1,\cdots,n_k}(dx_{\beta_1},\cdots,dx_{
\beta_k}),$$
where the $\beta_i$ have the property $(*)$.
\end{pf}

\subsection{Creation and Annihilation Operators on Locally Compact Spaces}

We define creation and annihilation operators for symmetric functions and measures on 
$\mathfrak{X}$. Assume a function $\varphi\in \mathcal{K}(X)$, a  function
$f\in \mathcal{K}$, a measure $\nu \in \mathcal{M}(X)$, a symmetric measure
$\mu\in \mathcal{M}(\mathfrak{X})$. We define
$$(a(\nu)f)(x_1,\cdots,x_n)=\int\nu(dx_0)f(x_0,x_1, \cdots,x_n)$$
or in another notation,
where $\alpha + c = \alpha + \{c\}$ means that the point $c$ is added to the set $\alpha$ and similar using
$\alpha \setminus c = \alpha \setminus \{c\}$ 
\begin{eqnarray*}(a(\nu)f)(x_\alpha)&=&\int\nu(dx_c) f (x_{\alpha+c})\\
 (a^+(\varphi)f)(x_\alpha)&=&\sum_{c\in\alpha}{\varphi(x_c)f(x_{\alpha\setminus c})}\\
 (a^+(\nu)\mu) (dx_\alpha)&=&\sum_{c\in\alpha}{\nu(dx_c) \mu(dx_{\alpha\setminus c})}\\
(a(\varphi)\mu)(dx_\alpha)&=&\int \varphi(x_c)\mu(dx_{\alpha+c})\end{eqnarray*}

If $\Phi$ is the function 
\begin{equation*}\Phi(\emptyset)=1,\, \Phi(x_\alpha) = 0 \mbox{ for } \alpha \neq \emptyset \end{equation*}
then 
$$a(\nu)\Phi=0$$
Similar if $\Psi$ is the measure defined by
$$ \Psi(f)=f(\emptyset),$$
then
$$ a(\varphi)\Psi=0.$$
One finds the commutation relations
\begin{eqnarray*}\left[a(\nu),a^+(\varphi)\right]&=&\int \nu (dx) \varphi(x)\\
\left[a(\varphi),a^+(\nu)\right]&=&\int \nu (dx) \varphi(x)\end{eqnarray*}
and obtains 
\begin{eqnarray*}\int (a^+(\nu)\mu)(dw)f(w)&=&\int\mu(dw)(a(\nu)f)(w)\\
\int (a(\varphi)\mu)(dw)f(w)&=&\int\mu(dw)(a^+(\varphi)f)(w)\end{eqnarray*}
We define the exponential measures and functions
\begin{eqnarray*}e(\varphi)&=&\Phi + \varphi + \varphi^{\otimes2} +\cdots =e^{a^+(\varphi)}\Phi\\
e(\nu)&=&\Psi + \nu + 
\nu^{\otimes2} +\cdots =e^{a^+(\nu)}\Psi.\end{eqnarray*}
So for $\alpha =\{a_1,\cdots ,a_n\}$
\begin{eqnarray*}
 e(\varphi)(x_\alpha) & = &\varphi(x_{a_1})\cdots \varphi(x_{a_n})\\
e(\nu)(dx_\alpha)& = &\nu(dx_{a_1})\cdots \nu(dx_{a_n}).\end{eqnarray*}

\subsection{Multiplication of Diffusions}
\begin{defi}
Let $X$ and $Y$ be two locally compact spaces. A continuous diffusion is a  vaguely continuous 
mapping
\begin{align*}
 \kappa : X &\to \mathcal{M}_+(X) \\
x& \mapsto \kappa_x \end{align*}
Using the oldfashioned way of writing we note
$$\kappa = \kappa_x(dy) = \kappa(x,dy)$$
\end{defi}
We consider three types of multiplication of  diffusions
\begin{enumerate}
\item Let $X_1,X_2,Y_1,Y_2$ be four locally compact spaces,and let
$$ \kappa_1:X_1 \to \mathcal{M}_+(Y_1)$$
$$ \kappa_2:X_2 \to \mathcal{M}_+(Y_2)$$
be continuous diffusions.Then 
\begin{align*}\kappa:X_1 \times X_2 &\to \mathcal{M}_+(Y_1\times Y_2)\\
 \kappa(x_1,x_2; dy_1,dy_2)& = \kappa_1(x_1,dy_1)\kappa_2(x_2,dy_2)\end{align*}
or 
$$ \kappa_{(x_1,x_2)} = \kappa_{1,x_1}\otimes\kappa_{2,x_2}. $$
\item Let $X,Y,Z$ be three locally compact spaces and
$$ \kappa_1: X \to \mathcal{M}_+(Y)$$
$$\kappa_2: Y \to \mathcal{M}_+(Z)$$
be continuous  diffusions, then
\begin{align*}\kappa: X &\to \mathcal{M}_+(Y\times Z)\\
\kappa(x;dy,dz)&= \kappa_1(x,dy)\kappa_2(y,dz)\end{align*}
So
$$ \iint \kappa(x;dy,dz)f(x,y) = \int \kappa_1(x,dy)\int \kappa_2(y,dz)f(y,z) $$
\item Let $X,Y,Z$ be three locally compact spaces and
$$ \kappa_1: X \to \mathcal{M}_+(Y)$$
$$\kappa_2: X \to \mathcal{M}_+(Z)$$
be continuous  diffusions, then
\begin{align*}\kappa: X &\to \mathcal{M}_+(Y\times Z)\\
\kappa(x;dy,dz)&= \kappa_1(x,dy)\kappa_2(x,dz)\end{align*}
So
$$\kappa_x = \kappa_{1,x}\otimes \kappa_{2,x}$$
\end{enumerate}

Using the positivity of the diffusions it is easy to see, that all three types of multiplications
yielt again positive continuous diffutions. We shall not introduce symbols for the multiplications, but use
instead the indication by differentials.

We are mostly interested into the diffusion
\begin{eqnarray*} \varepsilon:x \in X &\mapsto& \varepsilon_x \in\mathcal{M}_+(X)\\
\int \varepsilon_x(dy) \varphi(y)& = &\varphi(x) ,\end{eqnarray*}
so $\varepsilon_x$ is the point measure in the point $x$.
We have the three ways of defining the product of two point measures.
If the 4 variables $x_1,x_2,x_3,x_4$ are different, then we may define  at first the tensor product
\begin{eqnarray*}\varepsilon_{x_1}(dx_2)\varepsilon_{x_3}(dx_4)&=& \varepsilon_{x_1}\otimes \varepsilon_{x_3}(dx_2,dx_4)\\
\iint \varepsilon_{x_1}\otimes \varepsilon_{x_3}(dx_2,dx_4)f(x_2,x_4)&=& f(x_1,x_3).\end{eqnarray*}
The second way is
\begin{eqnarray*}\varepsilon_{x_1}(dx_2)\varepsilon_{x_2}(dx_3)= E_{x_1}(dx_2,dx_3)&=&
\varepsilon_{x_1}\otimes \varepsilon_{x_1}(dx_2,dx_3)\\ 
 \iint E_{x_1}(dx_2,dx_3)f(x_2,x_3) &=& f(x_1,x_1).\end{eqnarray*}

 The third possibility is
$$\varepsilon_{x_1}(dx_2)\varepsilon_{x_1}(dx_3)= \varepsilon_{x_1}\otimes \varepsilon_{x_1}(dx_2,dx_3)
.$$
We omit the variable x and write only the indices and denote $$\varepsilon_{x_b}(dx_c)= \varepsilon(b,c) \text{ and }
 E_{x_1}(dx_2,dx_3)= E(1;2,3).$$
 
 We want to define the product of the set $$\{\varepsilon(b_i,c_i):i=1,\cdots ,n\}.$$ Consider the set
 $$S=\{(b_1,c_1),\cdots ,(b_n,c_n)\},$$
 where all the $b_i$ and all the $c_i$ are different and $b_i\neq c_i$. We introduce in $S$ the structure of an oriented
 graph by defining the relation of right neighbor
 $$ (b,c)\triangleright (b',c')\Longleftrightarrow c=b'.$$
 An element $(b,c)$ has atmost one right neighbor, as $(b_i,c_i)\triangleright (b_j,c_j)$ and
 $(b_i,c_i)\triangleright (b_k,c_k)$ implies $b_j = b_k$ and $j=k$.
 So the
 components of the graph $S$ are either chains or circuits.
 
 We have to avoid the
expression $\varepsilon_x(dx)$. This notion makes no sense and if one wants to give it a sense, one runs into problems.
If $X$ is discrete, then $\varepsilon_x(dy) = \delta_{x,y}$ and $\varepsilon_x(dx)= \delta_{x,x} =1$. If $X=\mathbb{R}$,then
$\varepsilon_x(dy)= \delta(x-y)dy$ and if one wants to approximate Dirac's delta funtion one obtains 
$\varepsilon_x(dx)=\infty$.

Consider a circuit $$(1,2),(2,3),\cdots ,(k-2,k-1),(k-1,1).$$ It corresponds to a product 
 $$\varepsilon(1,2)\varepsilon(2,3)\cdots \varepsilon(k-2,k-1)\varepsilon(k-1,1).$$ Integrating over $x_2,\cdots ,x_{k-1}$
 one obtains $\varepsilon(1,1)$, which cannot be defined. So in order that $\prod_{i=1}^n \varepsilon(b_i,c_i)$ 
 can be defined,
 it is necessary, that the graph $S$ contains no circuits.

    In the contrary,if  $(1,2),(2,3),\cdots ,(k-2,k-1),(k-1,k)$ is a chain, then by extending the second way
    of multiplication we have
 \begin{multline*}\varepsilon(1,2)\varepsilon(2,3)\cdots \varepsilon(k-2,k-1)\varepsilon(k-1,k)\\= E(1;2,3,\cdots ,k) = 
 \varepsilon_{x_1}^{\otimes(k-1)}(dx_2,\cdots ,dx_k)\end{multline*}
 $$ \idotsint_{2,3,\cdots ,k} E(1;2,3,\cdots ,k)f(2,3,\cdots ,k) = f(x_1,\cdots ,x_1)$$
  above Denote
  $S_-=\{b_1,\cdots ,b_n\}$ and $S_+=\{c_1,\cdots ,c_n\}$. If $S$ contains no circuits, then  any $p \in S_-\setminus S_+$
  is the starting point of a (maximal) chain $(p,c_{p,1}),(c_{p,1},c_{p,2}),\cdots ,(c_{p,k-1},c_{p,k})$. Denote
  $\pi_p =\{c_{p,1},\cdots ,c_{p,k}\}$.
  We have 
  $$\varepsilon(p,c_{p,1})\varepsilon(c_{p,1},c_{p,2})\cdots \varepsilon(c_{p,k-1},c_{p,k})= E(p,\pi_p)$$
  where explicitely
  $$E(p,\pi_p)= \varepsilon_{x_p}^{\otimes \#\pi_p}(dx_{\pi_p}).$$
  Finally we define
  \begin{defi} If S contains no circuits, then
  $$E_S =\prod_{i=0}^n \varepsilon(b_i,c_i)= \prod_{p\in S_-\setminus S_+} E(p,\pi_p).$$
  \end{defi}
\subsection{Measure Valued Finite Particle Vectors}

We generalize the creation operator $a^+$ to the diffusion $\varepsilon:x\mapsto \varepsilon_x$ and define
$$ (a^+(\varepsilon(dy))f)(x_\alpha) = \sum_{c\in \alpha}\varepsilon_{x_c}(dy)f(x_{\alpha\setminus c}).$$
We write for short,if $b$ is an index 
$$a^+(\varepsilon(dx_b)) = a^+(dx_b)=a^+_b.$$
Similar
$$ a(\varepsilon_{x_b}) = a_b.$$
If $\alpha=\{b_1,\cdots ,b_n\}$ is a set, then
 \begin{align*}a_\alpha^+ &= a^+_{b_1}\cdots a^+_{b_n}\\
 a_\alpha&= a_{b_1}\cdots a_{b_n}\end{align*}
 We define the continuous diffusion
 \begin{align*}\Phi_\alpha :\mathfrak{X}&\to X^\alpha\\
 \Phi_\alpha &= a_\alpha^+\Phi \end{align*}
 where
$$\Phi_\alpha(x_\upsilon) = \varepsilon(\upsilon,\alpha)$$
for
 $\upsilon\cap \alpha =\emptyset$ and
\begin{equation*}\varepsilon(\upsilon,\alpha)=\begin{cases}
\sum_{h:\alpha\rightarrowtriangle \upsilon}\prod_{i=1}^n \varepsilon(h(b_i),b_i) 
&\text{ if } \#\alpha = \#\upsilon \\ 
 \,\,\,0 &\text{  otherwise  } \end{cases} 
 \end{equation*}
Here the sign $ \rightarrowtriangle $ signifies a bijective mapping. So the sum runs over all bijections from $\alpha$
to $\upsilon$. We call $\Phi_\alpha$ is a \emph{measure valued finite particle vector}.

Assume a set $\sigma=\{s_1,\cdots ,s_m\}$ and a set $S=\{(b_i,c_i):i=1,\cdots ,n\}$, where all the elements $b_i$ and
$c_i$ are different. Denote as above $S_-=\{b_1,\cdots ,b_n\}$ and $S_+=\{c_1,\cdots ,c_n\}$ and assume that
$\sigma \cap S_+ = \emptyset$. We extend the relation $\triangleright$ of right neighbor from $S$ to the pair
$(\sigma,S)$ by defining
$$ s \triangleright (b,c) \Longleftrightarrow s =b.$$
 If the graph $(\sigma,S)$ is without circuits and $(\sigma\cup S_+ \cup S_-)\cap \upsilon =\emptyset,$  
 then for any
 $f:\upsilon\rightarrowtriangle \sigma$, the graph   $$S\cup \{(c,f(c)),c \in\upsilon\}$$ 
 is without circuits, so
 we can define without problems $E_S \Phi_\sigma = \Phi_\sigma E_S$.
 
 The the set of components of the graph $(\sigma,S)$ is
 \begin{align*}\Gamma& = \Gamma(\sigma,S)= \Gamma_1 +\Gamma_2\\
 \Gamma_1&=
 \{\{ s,(s,c_{s,1}),(c_{s,1},c_{s,2}),\cdots,(c_{s,k_s-1},c_{s,k_s})\}; s \in \sigma \}
\\
\Gamma_2& = 
\{\{(t,c_{t,1}),(c_{t,1},c_{t,2}),\cdots,(c_{t,k_t-1},c_{t,k_t})\}; t \in S_-\setminus (S_+ + \sigma)\}.
\end{align*}
Put
\begin{align*} \xi_s& = \{s,c_{s,1},\cdots ,c_{s,k_s}\}\text{ for } s \in \sigma\\
\pi_t &= \{c_{t,1},\cdots ,c(t,k_t)\}  \text{ for } t\in S_-\setminus (S_+ + \sigma)\\
\pi &= S_+ +\sigma\\
\varrho& =  S_-\setminus (S_+ + \sigma)
\end{align*}
Then 
$$ \pi = \sum _{s \in \sigma} \xi_s + \sum _{t \in \varrho} \pi_t$$
So $\Phi_\sigma E_S$ is a continuous diffusion
\begin{align*}\Phi_\sigma E_S : \mathfrak{X }\times X^\varrho & \to \mathcal{M}_+(X^\pi)\\
\Phi_\sigma E_S(\upsilon+\varrho)&=\sum_{f:\upsilon\rightarrowtriangle\sigma}\prod_{c\in\upsilon}E(c,\xi_{f(c)})
\prod_{t\in\varrho}E(t,\pi_t)
\end{align*}
 \begin{defi} We denote by $\mathcal{G}_{n,\pi,\varrho}$ the set sums of elements of the form $\Phi_\sigma E_S$, such that
 $\sigma \cap S_+ = \emptyset $ and
the graph $(\sigma,S)$ is without circuits  
and that $$\varrho = S_-\setminus(S_+ +\sigma),\,\,\, \pi = S_+ +\sigma,\,\,\,
n = \# \sigma.$$\end{defi}
We extend the definitions of subsection1.3
 to $\Phi_\sigma$ and obtain for $b \neq c;\,\,\, b,c \notin \sigma$
\begin{eqnarray*} a_b^+a_c^+\Phi_\sigma &=&a_c^+a_b^+\Phi_\sigma \\ 
a_ba_c\Phi_\sigma &=&a_ca_b\Phi_\sigma \\
a_ba_c^+\Phi_\sigma &=& \varepsilon(b,c)\Phi_\sigma + a_c^+a_b\Phi_\sigma 
\end{eqnarray*}
We obtain for $c \notin \sigma$ using $a_c\Phi=0$
\begin{eqnarray*}a_c\Phi_\sigma& = &\sum _{b \in \sigma} \varepsilon(c,b)\Phi_{\sigma\setminus b}\\
a_c^+\Phi_\sigma&=& \Phi_{\sigma+c}.
\end{eqnarray*}
\begin{prop} Assume 
$$f = \Phi_\sigma E_S \in \mathcal{G}_{n,\pi,\varrho}.$$
Then for $c \notin \pi$ we have
$$(a^+_c\Phi_\sigma) E_S \in \mathcal{G}_{n+1,\pi+c,\varrho\setminus c}$$
and we define
$$ a^+_cf =(a^+_c\Phi_\sigma) E_S $$ 
If $c \notin \pi + \varrho $, then
$$(a_c\Phi_\sigma) E_S \in \mathcal{G}_{n-1,\pi,\varrho + c }$$
and we define
$$ a_cf =(a_c\Phi_\sigma) E_S $$
\end{prop}
\begin{pf}
The graph of $\Phi_\sigma E_S$ is $(\sigma,S)$.Its set of components is $\Gamma=\Gamma_1+\Gamma_2$ as above.
Assume $c \notin \pi$ and consider $(a^+_c\Phi_\sigma )E_S $. The corresponding graph is $(S',\sigma ')=(\sigma+c,S)$.
Denote by $\Gamma '=\Gamma_1'+\Gamma_2'$ the corresponding set of components of $(S',\sigma ')$. There are two cases: 
\begin{itemize}\item[a)] $c\notin S_-$, then $\Gamma_1' = \Gamma_1 +\{c\}, \Gamma_2' = \Gamma_2$ and $\pi '= \pi + c,
 \varrho ' = \varrho $.
\item[b)] $ c=t \in S_-$, then \begin{align*}
\Gamma_1' &= \Gamma_1 + \{t,(t,c_{t,1}),(c_{t,1},c_{t,2}),\cdots,(c_{t,k_t-1},c_{t,k_t}) \}\\
\Gamma_2'&= \Gamma_2\setminus \{(t,c_{t,1}),(c_{t,1},c_{t,2}),\cdots,(c_{t,k_t-1},c_{t,k_t})\}\end{align*}
 and  
 $\pi ' = \pi +c$ and $\varrho '= S_-'\setminus  (\sigma '+ S'_+)
= \varrho\setminus \{c\}$.
\end{itemize}

 In both cases the graph 
$(\sigma+c, S)$ contains no circuits and $(a^+_c\Phi_\sigma) E_S $ is defined and we set
$$ a_c^+(\Phi_\sigma E_S) = (a^+_c\Phi_\sigma) E_S .$$

Assume $c \notin \pi+\varrho$ and consider $(a_c\Phi_\sigma)E_S$. It consists of a sum of terms with a graph of the form 
 $(S'',\sigma '')=(\sigma\setminus b,S + (c,b))$. Denote the corresponding sets of components by $\Gamma_1'',\Gamma_2''
 $. Then we have
 \begin{align*} \Gamma_1''& = \Gamma_1\setminus \{b,(b_{c_1},b_{c_2}),\cdots \cdots (b_{c_{k-1}},b_{c_k})\}\\
 \Gamma_2''& = \Gamma_2+\{(c,b),(b_{c_1},b_{c_2}),\cdots \cdots (b_{c_{k-1}},b_{c_k})\}\end{align*}
 and $\pi '' = \pi, \varrho '' = \varrho+\{c\}$. The graph $(\sigma '',S'')$ has no circuits.
 \end{pf}
 \subsection{Admissible Sequences}
\begin{defi}A sequence
$$ W =(a^{\vartheta_n}_{c_n},\cdots ,a^{\vartheta_1}_{c_1})$$
with $c_1,\cdots ,c_n$ indices and $\vartheta_i = \pm 1$ and
$$ a^\vartheta_c= \left\{ \begin{array}{rcl} a^+_c &\textrm { for }& \vartheta=+1\\
a_c &\textrm { for }& \vartheta=-1 \end{array}\right.$$
is called \emph{admissible} if
$$ i>j \;\;\Longrightarrow \:\;\left\{ c_i \neq c_j \textrm{ or } \{c_i = c_j \textrm{ and } \vartheta_i=1,\vartheta_j =-1 \}
\right\} .$$
\end{defi}

So $W$ is admissible, if it contains only pairs ( not necessarily  being neighbors) of the form $(a_c^{\vartheta},
a_{c'}^{\vartheta '})$ with $c \ne c'$ or $(a_c^+,a_c)$ and no pairs of the form $(a_c,a_c),(a_c^+,a_c^+)$ or
$(a_c,a^+_c)$.

 \begin{defi}If $W$
is an admissible sequence,define
\begin{eqnarray*}\omega(W)&=& \{c_1,\cdots ,c_n\}\\
\omega_+(W)& = &\{c_i,1 \leq i \leq n: \vartheta_i = +1 \}\\
\omega_-(W)& = &\{c_i,1 \leq i \leq n: \vartheta_i = -1 \}.\end{eqnarray*}
\end{defi}
\begin{lem} If $W$ is admissible and $W_1 \subset W$, then $W_1$ is admissible. If $W$ is admissible,then
\begin{eqnarray*}
(a^+_c,W) \textrm{ is admissible }& \Longleftrightarrow & c \notin \omega_+(W)\\ 
(a_c,W) \textrm{ is admissible }& \Longleftrightarrow & c \notin \omega(W)\\ 
(W,a^+_c) \textrm{ is admissible }& \Longleftrightarrow & c \notin \omega(W)\\ 
(W,a_c) \textrm{ is admissible }& \Longleftrightarrow & c \notin \omega_-(W). 
\end{eqnarray*}
If $W=(W_2,W_1)$ is admissible,then
\begin{eqnarray*}
(W_2,a^+_c,W_1) \textrm{ is admissible }& \Longleftrightarrow & c \notin \omega(W_2)\cup \omega_+(W_1)\\ 
(W_2,a_c,W_1) \textrm{ is admissible }& \Longleftrightarrow & c \notin \omega_-(W_2)\cup \omega(W_1).
\end{eqnarray*}
\end{lem}
\begin{prop}
 Assume
$$ W =(a^{\vartheta_n}_{c_n},\cdots ,a^{\vartheta_1}_{c_1 }) $$
to be an addmissible sequence. Assume disjoint index sets $\pi$ and $\varrho$ and 
\begin{eqnarray*}\omega_+(W) \cap \pi& = &\emptyset \\ \omega_-(W)\cap (\pi +\varrho)&= &\emptyset .\end{eqnarray*}
Define for $k=1,\cdots ,n$
$$ W_k =(a^{\vartheta_k}_{c_k},\cdots ,a^{\vartheta_1}_{c_1 }) $$
Set 
 $\pi_0 = \pi,\varrho_0 = \varrho$ and
\begin{eqnarray*} \pi_k &= &\pi + \omega_+(W_k)\\
\varrho_k &=& \varrho\setminus(\omega_+(W_k)\setminus \omega_-(W_k)) +\omega_- (W_k)\setminus \omega_+(W_k),\end{eqnarray*}
where for sets $\alpha,\beta$
$$\alpha \setminus \beta = \alpha \setminus (\alpha \cap \beta).$$

Then 
$$ a^{\vartheta_k}_{c_k}: \mathcal{G}_{\pi_{k-1},\varrho_{k-1}} \to  \mathcal{G}_{\pi_{k},\varrho_{k}}$$
and the iterated application
$$ M =a^{\vartheta_n}_{c_n}\cdots a^{\vartheta_1}_{c_1 } : \mathcal{G}_{\pi,\varrho} \to \mathcal{G}_{\pi ',\varrho '} $$
with
\begin{eqnarray*} \pi ' &= &\pi + \omega_+(W)\\
\varrho ' &=& \varrho\setminus(\omega_+(W)\setminus \omega_-(W)) +\omega_- (W)\setminus \omega_+(W).\end{eqnarray*}
\end{prop}
\begin{defi} If
$$ W =(a^{\vartheta_n}_{c_n},\cdots ,a^{\vartheta_1}_{c_1 }) $$
is an admissible sequence we call 
$$ M =a^{\vartheta_n}_{c_n}\cdots a^{\vartheta_1}_{c_1 }$$
an admissible monomial.\end{defi}
\begin{pf}
We prove by induction. The case $k=1$ is trivial. Put $\omega_k = \omega(W_k)$ etc..Assume
\begin{eqnarray*} \pi _{k-1} &= &\pi + \omega_{+,k-1}\\
\varrho _{k-1} &=& \varrho\setminus(\omega_{+,k-1}\setminus \omega_{-,k-1}) +\omega_{-,k-1}\setminus \omega_{+,k-1}
\end{eqnarray*}
Assume $\vartheta_k = +1$. In order that $a^+_{c_k}$ is defined, $c_k \notin \pi_{k-1}$. But $ c_k \notin \pi$ by assumption
and $c_k \notin \omega_{+,k-1}$, as $W_k$ is admissible. So
$$ a^+_{c_k} : \mathcal{K}_{\pi_{k-1},\varrho_{k-1}} \to \mathcal{K}_{\pi_{k-1} + c_k, \varrho_{k-1}\setminus c_k}$$
Now 
$ \pi_{k-1} + c_k = \pi + \omega_{+,k} = \pi_k $
and it can be seen easily, that \newline
$\varrho_{k-1}\setminus c_k = \varrho_k$.

Assume now,that $\vartheta_k = -1$.In order that $a_{c_k}$ is defined, 
$$ c_k \notin \pi_{k-1}+ \varrho_{k-1} \subset \pi +\varrho + \omega_{k-1}.$$
But $c_k \notin \pi + \varrho $ by assumption and $c_k \notin \omega_{
k-1}$, as $W_k$ is admissible.
$$ a_{c_k} : \mathcal{K}_{\pi_{k-1},\varrho_{k-1}} \to \mathcal{K}_{\pi_{k-1} , \varrho_{k-1} + c_k}$$
But $\omega_{+,k} = \omega_{+,k-1}$ and 
\begin{eqnarray*} \pi_k &= &\pi+ \omega_{k-1} =\pi_{k-1}\\
\varrho_k & = & \varrho \setminus (\omega_{+,k}\setminus\omega_{-,k}) + (\omega_{-,k}\setminus\omega_{+,k})\\
& = & \varrho \setminus (\omega_{+,k-1}\setminus\omega_{-,k-1}) + (\omega_{-,k}\setminus\omega_{+,k}) + c_k\\
& = & \varrho_{k-1} + c_k .\end{eqnarray*}
\end{pf}
\begin{lem}
Assume $W=(W_2,W_1)$ to be admissible ,denote by $M_2$ and $M_1$ the correponding monomials and let $c \ne c'$  be indices.
If $(W_2, a^{\vartheta}_c, a^{\vartheta '} _{c'},W_1)$ is admissible, so is $(W_2, a^{\vartheta '}_{c'}, a^{\vartheta } _{c},W_1)$  and
\begin{eqnarray*} M_2 a^+_{c }a^+_{c'}M_1& = &M_2 a^+_{c '}a^+_{c}M_1 \\
M_2a_{c }a_{c'}M_1 & = & M_2 a^+_{c' }a^+_{c}M_1 \\
M_2 a_{c }a^+_{c'}M_1 & = & M_2 a^+_{c' }a_{c}M_1 + \varepsilon_{x_c}(dx_{c'}) M_2M_1 .\end{eqnarray*}
\end{lem}

\subsection{Wick's theorem}

Assume two finite index sets $\sigma,\tau$ and a finite set of pairs $S=\{(b_i,c_i): i \in I\}$, such that all $b_i$ and 
all $c_i$ are different and $b_i \neq c_i$.We extend the relation of right neighborhood to the triple $(\sigma,S,\tau)$
by putting for $(b,c) \in S, t \in \tau$
$$ (b,c)\triangleright t \Longleftrightarrow c=t $$
Consider a triple $(\sigma,S,\tau), \sigma \cap \tau = \emptyset$ and two finite sets $\upsilon,\beta$ such that the 
three sets $\sigma \cup S_+ \cup S_- \cup \tau$ and $\upsilon$ and $\beta$ are pairwise disjoint.
As 
$$(a^+_\sigma a_\tau \Phi_\upsilon)(\beta)= \sum_{\upsilon_1 +\upsilon_2 = \upsilon}
\varepsilon(\tau,\upsilon_1)\varepsilon(\beta,\sigma+\upsilon_2)$$
we find that the product $(a^+_\sigma a_\tau \Phi_\upsilon)(\beta)E_S$ is defined if the graph $(\sigma,S,\tau)$
is without circuits and we define the operator
$$a^+_\sigma a_\tau E_S = a^+_\sigma E_S a_\tau= E_S a^+_\sigma a_\tau $$ 
that way.

Consider an admissible sequence $W= (a^{\vartheta_n}_{c_n},\cdots ,a^{\vartheta_1}_{c_1})$ and the associated sets
$\omega_+,\omega_-$. We define the set $\mathfrak{P}(W)$ of all decompositions of $[1,n]$ .i.e.all sets
of subsets, of the form
\begin{align*}\mathfrak{p} &= \{\mathfrak{p}_+,\mathfrak{p}_-,\{q_i,r_i\}_{i \in I}\}\\
[1,n] & =  \mathfrak{p}_+ +\mathfrak{p}_- + \sum_{ i \in I}\{q_i,r_i\}\\
\mathfrak{p}_+ &\subset \;\omega_+,\; \mathfrak{p}_- \subset \omega_-,\;q_ i\in \omega_-,\;r_i \in \omega_+,\;q_i > r_i
\end{align*}
\begin{lem}
 Assume $W$ to be admissible and $ \mathfrak{p} \in \mathfrak{P}(W)$. Then the graph 
$(\sigma,S,\tau)$ with 
\begin{align*} \sigma &= \{c_s: s\in \mathfrak{p}_+\},&
S =&\{(c_{q_i},c_{r_i}):i \in I \},& \tau & = \{c_t: t \in \mathfrak{p}_-\} \end{align*}
is without circuits. \end{lem}
\begin{pf}
Let $(c_q,c_r) \triangleright (c_{q'},c_{r'})$, then $c_r =  c_{q'}$ and $r > q'$ as $W$ is admissible. By the definition
of $\mathfrak{p}$ we have $q > r$ and $q' > r'$, so $q>r>q'>r'$, and if we have a sequence		 
$(c_{q_1},c_{r_1})\triangleright \cdots \triangleright (c_{q_k},c_{r_k})$, then $q_1 > r_1 >\cdots q_k > r_k$
and as $\vartheta_1 = -1,\vartheta_k = +1$ we have $c_{q_1} \ne c_{r_k}$ as $W$ is admissible. This proves that 
 $S$ is without circuits. For the other components of  the graph one uses similar arguments.\end{pf}
\begin{defi} For $\mathfrak{p} \in \mathfrak{P}(W)$ we define
$$ \lfloor W \rfloor_\mathfrak{p} = \prod_{s \in \mathfrak{p}_+} a^+_{c_s}\prod_{i \in I} \varepsilon(c_{q_i},c_{r_i})
\prod_{t \in \mathfrak{p}_-} a_{c_t}$$
\end{defi}
\begin{thm}[Wick's theorem] If $W$ is admissible and if $M$ is the corresponding monomial, then
$$ M  = \sum_{\mathfrak{p} \in \mathfrak{P}(W)}\lfloor W \rfloor_\mathfrak{p}.$$
\end{thm}
\begin{pf}
We proceed by induction. The case $n=1$ is clear. We write for short $\mathfrak{p}_i = (q_i,r_i),
 \varepsilon(c(\mathfrak{p}_i)) = \varepsilon(c_{q_i},c_{r_i})$.Assume $$V= 
 (a^{\vartheta_n}_{c_n},\cdots ,a^{\vartheta_1}_{c_1})$$
to be admissible and $$N= a^{\vartheta_n}_{c_n}\cdots a^{\vartheta_1}_{c_1}.$$
Consider $W=(a_{c_{n+1}},V)$ and define an application
$ \varphi_-: \mathfrak{P}(W) \to \mathfrak{P}(V)$ consisting in erasing $n+1$. Then $n+1$ may occur in one of the
$\mathfrak{p}_i$, say in $\mathfrak{p}_{i_0}$, or in $\mathfrak{p}_-$. In the first case
$$ \varphi_-\mathfrak{p} = \{\mathfrak{p}_+ + \{ r_{i_0}\},\; (\mathfrak{p}_i)_{i\in I\setminus i_0},\; \mathfrak{p}_- \}$$
in the second case
$$ \varphi_-\mathfrak{p} = \{\mathfrak{p}_+ ,\; (\mathfrak{p}_i)_{i\in I},\; \mathfrak{p}_- \setminus
\{n+1\}\}.$$
Assume
$$\mathfrak{q}= \{ \mathfrak{q}_+,\;(\mathfrak{q}_j)_{j \in J},\; \mathfrak{q}_-\} \in \mathfrak{P}(V)$$
then
\begin{align*} \varphi^{-1}_- \mathfrak{q}& = \{ \mathfrak{p}: \varphi_- \mathfrak{p}= \mathfrak{q}\} =
\{ \mathfrak{p}^{(0)}, \mathfrak{p}^{(l)},l \in \mathfrak{q}_+ \}\\
 \mathfrak{p}^{(0)}& =\{\mathfrak{q}_+,(\mathfrak{q}_j)_{j \in J},\mathfrak{q}_- + \{n+1\}\} \\
 \mathfrak{p}^{(l)}& =\{\mathfrak{q}_+ \setminus l,(\mathfrak{q}_j)_{j \in J},(n+1,l),\mathfrak{q}_- \}.
 \end{align*}
 Consider
 \begin{align*} a_{c_{n+1}} \lfloor V \rfloor _\mathfrak{q}&= a_{c_{n+1}} \prod_{s \in \mathfrak{q}_+}a^+_{c_s}\prod_{j\in J}
  \varepsilon (c(\mathfrak{q}_j))
 \prod_{t \in \mathfrak{q}_-} a_{c_t}\\
 & =\prod_{s \in \mathfrak{q}_+}a^+_{c_s}\prod_{j\in J}
  \varepsilon (c(\mathfrak{q}_j))
 \prod_{t \in \mathfrak{q}_-+ \{n+1\}} a_{c_t}\\
  & \,\,\,\,\,+ \sum_{l \in \mathfrak{q}_+} \prod_{s \in \mathfrak{q}_+\setminus l}a^+_{c_s}\prod_{j\in J}
  \varepsilon (c(\mathfrak{q}_j)\varepsilon(c_{n+1},c_l)
 \prod_{t \in \mathfrak{q}_-} a_{c_t}\\
 &= \sum_{\mathfrak{p}\in \varphi^{-1}_-(\mathfrak{q})}\lfloor W \rfloor _\mathfrak{p}\end{align*}
 Finally
 $$a_{c_{n+1}} N  = \sum_{\mathfrak{q}\in \mathfrak{P}(V)}a_{c_{n+1}} \lfloor V \rfloor = 
 \sum _{\mathfrak{q} \in\mathfrak{P}(V)}\sum_{\mathfrak{p}\in \varphi^{-1}_-(\mathfrak{q})}\lfloor W \rfloor _\mathfrak{p}
 =\sum_{\mathfrak{p}\in \mathfrak{P}(W)}\lfloor W \rfloor _\mathfrak{p}$$

 Cosider now
 $W= (a^+_{c_{n+1}} ,V)$ and define an application
$ \varphi_+: \mathfrak{P}(W) \to \mathfrak{P}(V)$ consisting in erasing $n+1$ then
 \begin{align*} \varphi_+\mathfrak{p}&
  = \{\mathfrak{p}_+ \setminus\{n+1\},\; (\mathfrak{p}_i)_{i\in I},\; \mathfrak{p}_- \}\\
 \varphi_+^{-1} \mathfrak{q}& =
\{\mathfrak{q}_+ + \{n+1\},(\mathfrak{q}_j)_{j \in J},\mathfrak{q}_- + \}\\
 a_{c_{n+1}}^+ \lfloor V \rfloor _\mathfrak{q}& = \lfloor W \rfloor _{ \varphi_+^{-1}\mathfrak{q}}.\end{align*}
By the same reasoning as above one fishes the proof. \end{pf}

\subsection{Representation of Unity}
We extend the functional $\Psi$ to $\mathcal{G}_{n,\pi,\varrho}$ by putting
$$ \Psi\Phi_\sigma=\begin{cases}1& \text{ for } \sigma = \emptyset\\
0 &\text{ otherwise }\end{cases}$$
and
$$\Psi\Phi_\sigma E_S = (\Psi\Phi_\sigma)E_S.$$
\begin{defi}
Assume $$W= 
 (a^{\vartheta_n}_{c_n},\cdots ,a^{\vartheta_1}_{c_1})$$
to be admissible and $$M
= a^{\vartheta_n}_{c_n}\cdots a^{\vartheta_1}_{c_1}.$$
Then we define
$$ \langle M \rangle = \sum_{\mathfrak{p }\in \mathfrak{P}_0(W)} \lfloor W \rfloor _{\mathfrak{p}}$$
Here $\mathfrak{P}_0(W)$ is the set of partitions of $[1,n]$ into pairs $\{q_i,r_i\}_{i = 1,\cdots ,n/2}$ such that 
$ q_1 > r_i, \vartheta_{q_i} = -1, \vartheta_{r_i} = +1$. 
$$ \lfloor W \rfloor _{\mathfrak{p}} = \prod_i \varepsilon(b_i,c_i)$$
If $ n$ is odd or $\mathfrak{P}_0(W)$ is empty, then
$\langle M \rangle  =0$.
\end{defi}
As a consequence of Wick'theorem 1.7.1 we obtain 
\begin{prop}
We obtain
\begin{equation*} (M\Phi)(\emptyset) = \Psi M \Phi = \begin{cases}0 & \text { if $\vartheta_1 +\cdots+  \vartheta_n \ne 0$}\\
\langle M \rangle & \text{ if }\vartheta_1+ \cdots + \vartheta_n = 0\end{cases}\end{equation*}

\end{prop}
If $M$ is admissible, then
$$ (M\Phi_\beta)(\alpha) = \Psi a_\alpha M a^+_\beta\Phi = \langle a_\alpha M a^+_\beta \rangle. $$
We shall use this notation very often.
\begin{thm}If $M=M_2M_1$ is admissible,then
$$ \langle M \rangle = \int_\alpha \langle M_2 a^+_\alpha\rangle\langle a_\alpha M_1\rangle. $$
\end{thm}
\begin{pf}
Assume
\begin{align*} M&=a^{\vartheta_n}_{c_n}\cdots a^{\vartheta_1}_{c_1}\\
M_2&=a^{\vartheta_n}_{c_n}\cdots a^{\vartheta_k}_{c_k}\\
M_1&=a^{\vartheta_{k-1}}_{c_{k-1}}\cdots a^{\vartheta_1}_{c_1}
\end{align*}
We prove by induction with respect to $k$.For $k=n$ we have
$$\Psi a^+_\alpha a_\alpha M \Phi = \begin{cases} \langle M\rangle &\text{ for } \alpha=\emptyset\\
0& \text{otherwise}\end{cases}$$
Integration yields the result. Put $M_2'=a^{\vartheta_n}_{c_n}\cdots a^{\vartheta_{k+1}}_{c_{k+1}}$.
Assume $\vartheta_k = -1$. Then
\begin{multline*} \int_\alpha \langle M_2 a^+_\alpha\rangle\langle a_\alpha M_1\rangle =
\int_\alpha \langle M_2'a_{c_k} a^+_\alpha\rangle\langle a_\alpha M_1\rangle\\=
\int_\alpha \sum_{b\in \alpha}\langle M_2' a^+_{\alpha\setminus b}\rangle\langle a_\alpha M_1\rangle
\varepsilon(c_k,b)
=\int_\alpha \int_b\langle M_2' a^+_\alpha\rangle\langle a_{\alpha+b} M_1\rangle \varepsilon(c_k,b)\\
=\int_\alpha \langle M_2' a^+_\alpha\rangle\langle a_\alpha a_{c_k} M_1\rangle\end{multline*}
In a similar way one proves
$$\int_\alpha \langle M_2' a^+_\alpha\rangle\langle a_\alpha a^+_{c_k}M_1\rangle=
\int_\alpha \langle M_2' a^+_{c_k} a^+_\alpha\rangle\langle a_\alpha M_1\rangle$$
\end{pf}

\subsection{Duality}We repete the definitions of subsection 1.5
Assume $\Phi_\sigma E_S \in \mathcal{G}_{n,\pi,\varrho}$, then $n = \# \sigma, \pi = \sigma + S_+,
\varrho =S_-\setminus(\sigma + S_+)$
and the set
of components of the graph $(\sigma,S)$ is
 \begin{align*}\Gamma& = \Gamma_1 +\Gamma_2\\
 \Gamma_1&=
 \{\{ s,(s,c_{s,1}),(c_{s,1},c_{s,2}),\cdots,(c_{s,k_s-1},c_{s,k_s})\}; s \in \sigma \}
\\
\Gamma_2& = 
\{\{(t,c_{t,1}),(c_{t,1},c_{t,2}),\cdots,(c_{t,k_t-1},c_{t,k_t})\}; t \in S_-\setminus (S_+ + \sigma)\}.
\end{align*}
Put
\begin{align*} \xi_s& = \{s,c_{s,1},\cdots ,c_{s,k_s}\}\text{ for } s \in \sigma\\
\pi_t &= \{c_{t,1},\cdots ,c_{t,k_t}\}  \text{ for } t\in S_-\setminus (S_+ + \sigma)\\
\end{align*}

Then 
$$ \pi = \sum _{s \in \sigma} \xi_s + \sum _{t \in \varrho} \pi_t$$
So $\Phi_\sigma E_S$ is a continuous diffusion
\begin{align*}\Phi_\sigma E_S : \mathfrak{X }\times X^\varrho & \to \mathcal{M}_+(X^\pi)\\
\Phi_\sigma E_S(\alpha+\varrho)&=\sum_{f:\alpha\rightarrowtriangle\sigma}\prod_{c\in\upsilon}E(c,\xi_{f(c)})
\prod_{t\in\varrho}E(t,\pi_t)
\end{align*}
with
$\varrho = S_-\setminus (\sigma + S_+)$ and 
$$\Phi_\sigma E_S = \Phi_\sigma \prod _{s \in \sigma} E(s,\pi_s) \prod_{r \in \pi_r} E(r,\pi_r)$$
$$(\Phi_\sigma E_S)(x_\alpha)  = \sum_{f:\alpha \rightarrowtriangle \sigma}\prod_{c \in \alpha} \varepsilon (c,f(c)) 
 E(f(c),\pi_{f(c)} )\prod_{r \in \varrho} E(r,\pi_r)$$ $$ =
  \sum_{f:\alpha \rightarrowtriangle \sigma}\prod_{c \in \alpha} 
 E(c,\xi_{f(c)} )\prod_{r \in \varrho} E(r,\pi_r)$$
 with $\xi_b = \{b\} + \pi_b $.
 
 We have
 $$ \pi = \sum_{s \in \sigma}\xi_s + \sum_{r \in \varrho} \pi_r.$$
 
 We may define the product
\begin{align*}(F_1,F_2) &\in \mathcal{G}_{n,\pi_1,\varrho_1} \times \mathcal{G}_{n,\varrho_2,\pi_2} \mapsto
F_1 F_2\in \mathcal{G}_{n,\pi_1+ \pi_2,\varrho_1 \cup \varrho_2}\\
 F_1F_2(x_\alpha)& = \sum_{\substack {f_1:\alpha \rightarrowtriangle \sigma_1\\
 f_2:\alpha \rightarrowtriangle \sigma_2}}
\prod_{c \in \alpha} E (c;\xi_{f_1(c)} + \xi_{f_2(c)})
\prod_{r \in \varrho_1 \cup \varrho_2} E(r,\pi_{1,r} + \pi_{2,r})
\end{align*}
 assuming $\pi_1 \cap \pi_2 = \emptyset$
  and using $E(r,\pi_1)E(r,\pi_2)= E(r;\pi_1 + \pi_2)$ for 
 $\pi_1 \cap \pi_2 = \emptyset$ setting$E(r;\emptyset)=1$.
 
 Assume a positive measure $\lambda$ on $X$ and denote $e(\lambda)= \lambda$ and $$e(\lambda)(dx_\alpha) =
 \lambda^{\otimes\alpha}(dx_\alpha)=
 \lambda(dx_\alpha) = \lambda(\alpha) = \lambda_\alpha.$$ Define
 $$ \mathcal{G}_{\pi,\varrho} = \bigoplus_{n=0}^\infty{G}_{n,\pi,\varrho} $$
For $\pi_1 \cap \pi_2 = \emptyset$ we define the application
$$(F_1,F_2) \in \mathcal{G}_{\pi_1,\varrho_1} \times \mathcal{G}_{\pi_2,\varrho_2} \mapsto\langle
 F_1 ,F_2\rangle_\lambda,$$
where$\langle F_1 ,F_2\rangle_\lambda $ is a
continuous diffusion 
\begin{align*}X^{\varrho_1 \cup \varrho_2} &\to\mathcal{M}_+(X^{\pi_1+\pi_2})\\
\langle F_1 ,F_2\rangle_\lambda &= \int (F_1F_2)(x_\alpha)\lambda(dx_\alpha)\end{align*}
Using only the essential information we write
$$ \langle F_1,F_2 \rangle_\lambda (\pi_1 + \pi_2, \varrho_1 \cup \varrho_2) =
\int \lambda_\alpha F_1(\alpha, \pi_1,\varrho_1)F_2(\alpha,\pi_2,\varrho_2).$$
\begin{defi}We define the measure
$\Lambda$ on $X^k$ given by
\begin{multline*}\int\Lambda(1,\cdots ,k)f(1,\cdots ,k)\\=  
\int\Lambda(dx_1,\cdots ,dx_k) f(x_1,\cdots ,x_k)= \int\lambda(dx)f(x,\cdots ,x)\end{multline*}
\end{defi}
So we have
$$\lambda(1)\varepsilon(1,2)\cdots \varepsilon(k-1,k) = \Lambda (1,2,\cdots ,k)$$
If $F=\Phi_\sigma E_S \in \mathcal{G}_{\pi,\varrho}$, then $F\lambda_\varrho$ is a continuous diffusion
$\mathfrak{X}\to \mathcal{M}^{\pi+\varrho}$ given by
$$(\Phi_\sigma E_S)(x_\alpha)  = \sum_{f:\alpha \rightarrowtriangle \sigma}\prod_{c \in \alpha}  
 E(c,\xi_{f(c)} )\prod_{r \in \varrho} \Lambda(\{r\}+\pi_r)$$
 Assume $$ W =(a^{\vartheta_n}_{c_n},\cdots ,a^{\vartheta_1}_{c_1 }), $$ to be an admissible sequence assume
 $\vartheta_1+\cdots \vartheta_n = 0 $. Define a usual $\omega_\pm(W)= \{c_i:\vartheta_i=\pm 1\}$.
Recall prop 1.8.1
 $$ \langle M \rangle = \sum_{\mathfrak{p }\in \mathfrak{P}_0(W)} \lfloor W \rfloor _{\mathfrak{p}}$$

Here $\mathfrak{P}_0(W)$ is the set of partitions of $[1,n]$ into pairs $\{q_i,r_i\}_{i = 1,\cdots ,n/2}$ such that 
$ q_1 > r_i, \vartheta_{q_i} = -1, \vartheta_{r_i} = +1$. 
$$ \lfloor W \rfloor _{\mathfrak{p}} = \prod_i \varepsilon(b_i,c_i)$$
Call $S(\mathfrak{p})$ the graph related to $\mathfrak{p}$ and $\Gamma(S(\mathfrak{p}))$ the set of components of the graph.
To any $s\in S_-(\mathfrak{p})\setminus S_+(\mathfrak{p})$ there is a component. As 
$$ S_-(\mathfrak{p})\setminus S_+(\mathfrak{p}) = \omega_-(W)\setminus\omega_+(W) = \varrho$$
for all $\mathfrak{p}$ we obtain
$$ \langle M\rangle \lambda(\varrho) = \sum_{\mathfrak{p }\in \mathfrak{P}_0(W)} \lfloor W \rfloor _{\mathfrak{p}}
\lambda(\varrho)= \sum_{\mathfrak{p }\in \mathfrak{P}_0(W)}\prod_{\gamma \in \Gamma(S(\mathfrak{p}))}\Lambda(\gamma).$$
\begin{defi}
Assume
$$ W =(a^{\vartheta_n}_{c_n},\cdots ,a^{\vartheta_1}_{c_1 }), $$
to be an admissible sequence, then define the formal adjoint sequence by
$$ W^+  =(a^{-\vartheta_1}_{c_1},\cdots ,a^{-\vartheta_n}_{c_n }). $$
If $M$ is the monomial corresponding to $W$, we denote by $M^+$ the monomial corresponding to $W^+$.
\end{defi}
Using the symmetry of $\Lambda$ one sees that
\begin{prop}
$$ \langle M\rangle \lambda(\omega_-(W)\setminus\omega_+(W))=
 \langle M^+\rangle \lambda(\omega_+(W)\setminus\omega_-(W)).$$
\end{prop}

\begin{thm}
Let $W_1$ and $W_2$ be two addmissible sequences such that $W_0 = (W_1^*,W_2)$ is admissible too. Denote 
$ \omega_i = \omega(W_i), \omega_{ \pm i}= \omega_\pm (W_i)$  and $M_i$ the corresponding monomials for $i=0,1,2$.
Then 
$$ M_i: \Phi \to \mathcal{G}_{\pi_i,\varrho_i}$$
with $\pi_i = \omega_{+i}, \varrho_i = \omega_{-i}\setminus \omega_{+i}$.
We have $ \pi_1 \cap \pi_2 = \emptyset$ and
$$ \langle M_1\Phi, M_2\Phi\rangle_\lambda \lambda(\varrho_1 \cup \varrho_2) = 
\langle  M_1^+M_2 \rangle\lambda(\varrho_0)$$
\end{thm}
\begin{pf}
 Call  $\omega_{ \pm i}^+= \omega_\pm (W_i^+)=  \omega_{ \mp i}$. 
 Put
 \begin{align*}
 \omega_i & =\alpha_i +\beta_i + \gamma_i\\
 \beta_i & =\omega_{+i}\cap\omega_{-i}\\
 \alpha_i& = \omega_{+i}\setminus \omega_{-i} = \omega_{-i}^+\setminus\omega_{+i}^+\\
 \gamma_i& = \omega_{-i}\setminus \omega_{+i} = \omega_{+i}^+\setminus\omega_{-i}^+ = \varrho_i
 \end{align*}
 Using
 $$\langle a_\zeta M_1\rangle\lambda(\zeta+\gamma_1\rangle = \langle M^+_1 a^+_\zeta\rangle\lambda(\alpha_1)$$
 we have
 \begin{multline*}\langle M_1\Phi, M_2\Phi\rangle_\lambda \lambda(\varrho_1 \cup \varrho_2)=
  \int_\zeta (M_1\Phi)(\zeta)((M_2\Phi)(\zeta)\lambda(\zeta+(\varrho_1\cup\varrho_2))
  \\
 = \int_\zeta \langle M_1^+ a^+_\zeta \rangle\langle a_\zeta M_2\rangle \lambda(\alpha_1+
 (\gamma_1\cup\gamma_2)\setminus \gamma_1)
 = \langle \Phi, M_1^+M_2 \rangle_\lambda\lambda(\varrho_0),\end{multline*}
 by theorem 1.8.1. As $W$ is admissible we have
 $$\alpha_1+(\gamma_1\cup\gamma_2)\setminus \gamma_1 = \alpha_1+\gamma_2\setminus\gamma_1 = \gamma_0 = \varrho_0$$
 \end{pf}
 \subsection{Integration of Normal Ordered Monomials}
\begin{defi} An admissible sequence
 $$ W =(a^{\vartheta_n}_{c_n},\cdots ,a^{\vartheta_1}_{c_1 }), $$ 
 is called \emph {normal ordered}  if all the creators $a^+_c$ are at left of the annihilators $a_c$. ,i.e.
 $$ \vartheta_i = +1,\vartheta_j =-1 \Longrightarrow i>j.$$
 \end{defi}
 Using the commutation relations it is clear, that any normal odered monomial can be brought into the form
 \begin{multline*} a^+(dx_{s_1})\cdots a^+(dx_{s_l})a^+(dx_{t_1})\cdots a^+(dx_{t_m})\\a(x_{t_1})\cdots a(x_{t_m}) 
 a(x_{u_1})\cdots a(x_{u_n}) = a^+_{\sigma + \tau} a_{\tau + \upsilon}.\end{multline*}
 with
 \begin{align*}\sigma& = \{s_1,\cdots ,s_l\},&\tau& = \{t_1,\cdots ,t_m\}, & \upsilon &= \{u_1,\cdots ,u_n\}.
 \end{align*}
 Assume five finite , pairwise disjoint index sets $\pi,\sigma,\tau,\upsilon,\varrho$ and consider the admissible
 monomial $a_\pi a_{\sigma + \tau}^+ a_{\tau + \upsilon} a^+_\varrho$. For the associated graph $S$ we have
 $S_+ = \sigma +\tau+\varrho,\; S_-= \pi + \tau+ \upsilon$.So $S_-\setminus S_+= \pi+\upsilon$
   Following the last subsection
 $\langle a_\pi a_{\sigma + \tau}^+ a_{\tau + \upsilon} a^+_\varrho\rangle\lambda_{\pi+\upsilon}$ is for fixed
 $\#\pi,\#\sigma,\#\tau,\#\upsilon,\#\varrho$ a measure on $X^{\#(\pi+\sigma+\tau+\upsilon+\varrho)}$.
 Letting the numbers $\#\pi,\#\sigma,\#\tau,\#\upsilon,\#\varrho$ run from 0 to $\infty$ we arrive at a measure
 $\mathfrak{m}$ on $\mathfrak{X}^5$
 $$ \mathfrak{m}= \mathfrak{m}(\pi,\sigma,\tau,\upsilon,\varrho)=
 \langle a_\pi a_{\sigma + \tau}^+ a_{\tau + \upsilon} a^+_\varrho\rangle\lambda_{\pi+\upsilon}.$$
 Using thm 1.9.1 we obtain
 \begin{multline*}\mathfrak{m}= \int_\omega(a_{\sigma+\tau} a_\pi^+\Phi)(\omega)(a_{\tau+\upsilon}a_\varrho^+ \Phi)(\omega)
 \lambda_{\omega + \sigma +\tau+\upsilon}
 \\= \int_\omega \varepsilon (\sigma +\tau +\omega,\pi)\varepsilon(\tau+\upsilon+\omega,\varrho )
  \lambda_{\omega + \sigma +\tau+\upsilon}
 \end{multline*}
 If $\varphi \in \mathcal{K}(\mathfrak{X}^5)$ then
 $$\int\mathfrak{m}(\pi,\sigma,\tau,\upsilon,\varrho)\varphi( \pi,\sigma,\tau,\upsilon,\varrho)=\int
 \varphi(\sigma + \tau+\omega,\sigma,\tau,\upsilon,\tau+\upsilon+\omega)\lambda_{\omega + \sigma +\tau+\upsilon}.$$
 Assume a Hilbert space $\mathfrak{k}$ with countable base. We write the scalar product $x,y \mapsto \langle x,
 y \rangle $ often in the form $x^+y$ introducing the dual vector $x^+$. We denote by $B(\mathfrak{k})$ the space
 of bounded linear operarors on $\mathfrak{k}$.We provide $B(\mathfrak{k})$ with the operator norm topology.
 If $ A \in B(\mathfrak{k})$, then $A^+$ denotes the adjoint.
 
 Assume a function $F: \mathfrak{X}^3 \to B(\mathfrak{k})$ locally $\lambda$-integrable, i.e. locally integrable
 with respect to $e(\lambda)^{\otimes 3}$, and $f,g \in \mathcal{K}_s(\mathfrak{X},\mathfrak{k})$ ( continuous in the 
 normtopology of $\mathfrak{k}$).The integral
 \begin{multline*}\int\mathfrak{m}(\pi,\sigma,\tau,\upsilon,\varrho)f^+(\pi) F(\sigma,\tau,\upsilon)g(\varrho)\\
 =\int f^+(\sigma+\tau+\omega)F(\sigma,\tau,\upsilon)g(\tau+\upsilon+\omega)\lambda_{\omega + \sigma +\tau+\upsilon}
 =\langle f |\mathcal{B}(F)| g\rangle \end{multline*}
 exists and defines a sesquilinear form on $\mathcal{K}_s(\mathfrak{X},\mathfrak{k})$.
 
 Using again theorem 1.9.1 we find
 $$ \mathfrak{m}(\pi,\sigma,\tau,\upsilon,\varrho)= \int (a^+_\pi \Phi(\omega)
 (a^+_{\sigma+\tau}a_{\tau+\upsilon}a^+_\varrho\Phi)(\omega)\lambda_{\upsilon+\omega}.$$
 Remark
 \begin{align*}(a^+_\sigma \Phi_\varrho)(\omega) &= 
 \sum_{\alpha \subset \omega} \varepsilon(\alpha,\sigma)\Phi_\varrho(\omega\setminus\alpha)\\
  (a^+_\tau a_\tau \Phi_\varrho)(\omega) &= 
 \sum_{\alpha \subset \omega}\varepsilon(\alpha,\tau) \Phi_\varrho(\omega)\end{align*}
 
 and obtain
 $$ \int\mathfrak{m}(\pi,\sigma,\tau,\upsilon,\varrho)f^+(\omega) F(\sigma,\tau,\upsilon)g(\varrho)=
 \int f^+(\omega) ((\mathcal{O}F)g)(\omega) = \langle f, \mathcal{O}(F)g \rangle $$
 with
 $$((\mathcal{O}(F)g )(\omega)= \sum_{\alpha \subset \omega}\sum_{\beta \subset \omega \setminus\alpha}
 \int_\upsilon \lambda_\upsilon F(\alpha,\beta,\upsilon)g(\omega\setminus\alpha+ \upsilon).$$
 So $\mathcal{O}(F)$ is a mapping from $\mathcal{K}_s(\mathfrak{X})$ into the locally $\lambda$-integrable functions
 on $\mathfrak{X}$.
 
 Using the results of the last subsection ,we have
 $$\mathfrak{m}(\pi,\sigma,\tau,\upsilon,\varrho)= \mathfrak{m}(\varrho,\upsilon,\tau,\sigma,\pi)$$
 and  obtain
 $$\langle f, \mathcal{O}(F)g \rangle = \overline{\langle g,\mathcal{O}(F^+) f\rangle} = 
 \langle \mathcal{O}(F^+)f,g \rangle $$
 with
 $$F^+(\sigma,\tau,\upsilon)=(F(\upsilon,\tau,\sigma))^+.$$
 
 Consider
 $$\mathfrak{m}(\pi,\sigma_1,\tau_1,\upsilon_1,\sigma_2,\tau_2,\upsilon_2,\varrho)
 = \langle a_\pi a^+_{\sigma_1+\tau_1}a_{\tau_1+\upsilon_1}a^+_{\sigma_2+\tau_2}a_{t_2+\upsilon_2} a^+_\varrho \rangle
 \lambda_{\pi + \upsilon_1 +\upsilon_2}.$$
 Assume $F,G:\mathfrak{X}^3 \to B(\mathfrak{k})$ to be $\lambda$-measurable and define
 \begin{multline*} \langle f| \mathcal{B}(F,G)| g \rangle \\=
 \int\mathfrak{m}(\pi,\sigma_1,\tau_1,\upsilon_1,\sigma_2,\tau_2,\upsilon_2,\varrho)
 f^+(\pi)F(\sigma_1,\tau_1,\upsilon_1)G(\sigma_2,\tau_2,\upsilon_3)g(\varrho)\end{multline*}
 provided the integral exists in norm. 
 
 As it might be guessed from Wick's theorem, there exists an $H$ such that $\mathcal{B}(F,G) = \mathcal{B}(H)$. 
 In fact we have the following theorem basically due to P.A. Meyer.
 \begin{thm}[Meyer's formula]

$$H(\sigma,\tau,\upsilon) = \sum
\int_\kappa \lambda_\kappa 
F(\alpha_1,\alpha_2+\beta_1+\beta_2,\gamma_1+\gamma_2+\kappa)
G(\kappa+\alpha_2+\alpha_3,\beta_2+\beta_3+\gamma_2,\gamma_3)$$
where the sum runs through all  indices $\alpha_1,\cdots ,\gamma_3$ with
$$\begin{array}{rcl} \alpha_1+\alpha_2+\alpha_3 &=& \sigma\\
\beta_1+\beta_2+\beta_3 &=&\tau\\\gamma_1+\gamma_2+\gamma_3&=&\upsilon\end{array}.$$
 \end{thm}
 This formula was proven in \cite{vW2007} for $C_c$-functions and can be easily extended to our more general case.
 Our formula differs fom Meyer's formula as the letters $\sigma,\tau,\cdots $ denote index sets and not sets
 of real numbers. It is more general because it holds for for any positive positive measure and not only for the
 Lebesgue measure. We shall noe repeat the proof and we shall not use it in this paper.
 
 We have
 $$ \langle f|\mathcal{B}(F,G)|g \rangle = \langle \mathcal{O}(F^+)f,\mathcal{O}(G)g\rangle. $$
 So a sufficient condition for the existence of $\langle f|\mathcal{B}(F,G)|g \rangle$ is that $\mathcal{O}(F^+)$
 and $\mathcal{O}(G)$ are bounded operators from $\mathcal{K}_s(\mathfrak{X},\mathfrak{k})$ provided with the 
 $L^2(\mathfrak{X},\mathfrak{k},\lambda)$-norm into $L^2(\mathfrak{X},\mathfrak{k},\lambda)$.

\section{Skorokhod Integral}
\subsection{Definition}
Denote $\mathfrak{R}=\mathfrak{W}(\mathbb{R})$ , so
$$ \mathfrak{R}= \{\emptyset\} + \mathbb{R}+\mathbb{R}^2 + \cdots .$$
We provide $\mathfrak{R}$ with the Lebesgue measure $dw$ . So for a symmetric function
\begin{equation*}\begin{split}&\int_{\Delta(\mathfrak{R})}f(w)dw = \int _\Delta f(t_\alpha) dt_\alpha\\&=
f(\emptyset)+\sum_{n=1}^\infty \frac{1}{n!}\idotsint_{\mathbb{R}^n}dt_1\cdots dt_n f(t_1,\cdots ,t_n)\\&=
f(\emptyset)+\sum_{n=1}^\infty \idotsint_{t_1<\cdots< t_n}dt_1\cdots dt_n f(t_1,\cdots ,t_n)
\end{split}
\end{equation*}
\begin{defi}Assume a Banach algebra $\mathfrak{B}$ and a function $x$
\begin{eqnarray*} x: \mathbb{R} \times \mathfrak{R}^k &\to &\mathfrak{B}\\
(t,w_1,\cdots ,w_k) & \mapsto & x_t(w_1,\cdots ,w_k) \end{eqnarray*}
symmetric in any of the variables $w_1,\cdots ,w_k$ and locally integrable in norm with respect to the Lebesgue
measure on $\mathbb{R} \times \mathfrak{R}^k $. Assume a Lebesgue integrable function $f:\mathbb{R}\to \mathbb{C}$
The \emph {Skorokhod} integral $\oint^j(f)x$ is given by
\begin{multline*}
(\oint^j(f)x)(t_{\alpha_1},\cdots ,t_{\alpha_k})=\\ \sum_{c \in \alpha_j} f(t_c)\;x_{t_c}(t_{\alpha_1},\cdots ,
t_{\alpha_{j-1}},t_{\alpha_j \setminus c},t_{\alpha_{j+1}},\cdots ,t_{\alpha_k})
\end{multline*}
\end{defi}
\begin{rem} The function
$$ (w_1,\cdots ,w_k)\in \mathfrak{X}^k \mapsto (\oint^j(f) x)(w_1,\cdots ,w_k) \in \mathfrak{B}$$
is symmetric in any of the variables $w_i$ and locally integrable.

\end{rem} 
\begin{pf}The symmetry is trivial, for local integrability it is sufficient,that all functions have values
$\geq 0$. Let $g(w_1,\cdots ,w_k)$ be a continuous function $\ge 0$with compact support, symmetric in any of the $w_i$,
then by the sum integral lemma
\begin{equation*}\begin{split}
 &\int_{\Delta(\mathfrak{R})}\cdots \int_{\Delta(\mathfrak{R})}(\oint^j(f)x)(t_{\alpha_1},\cdots ,t_{\alpha_k})
g(t_{\alpha_1},\cdots ,t_{\alpha_k})dt_{\alpha_1}\cdots dt_{\alpha_k}
 \\
&=\int_{\mathbb{R}}\int_{\Delta(\mathfrak{R})}\cdots \int_{\Delta(\mathfrak{R})}f(t_c)
x_{t_c}(t_{\alpha_1},\cdots ,t_{\alpha_k})\\&\,\,\;\;\;
g(t_{\alpha_1},\cdots ,\cdots ,t_{\alpha_j+c},\cdots ,t_{\alpha_k})dt_c dt_{\alpha_1}\cdots dt_{\alpha_k}
 < \infty
\end{split}\end{equation*}
\end{pf}
\subsection{A Skorokhod Integral Equation}
\begin{defi}
Consider the subset
$$ \{(t_{\alpha_1},\cdots ,t_{\alpha_k}) \in \mathfrak{R}^k: \mbox{ all } t_i \mbox{ for }
 i \in \alpha_1 +\cdots +\alpha_k
\mbox { are different } \}.$$
This subset differs from the set $\mathfrak{R}^k$ by a nullset. We define on that
set a mapping $\Xi$ onto $\mathfrak{S}(\mathbb{R} \times \{1,\cdots ,k\})$,where
          $\mathfrak{S}$ denotes the set of finite subsets, by mapping
$$(t_{\alpha_1},\cdots ,t_{\alpha_k})\mapsto \xi = \{(s_1,i_1),\cdots ,(s_n,i_n) \}$$
where
\begin{align*} t_{\alpha_1}+\cdots +t_{\alpha_k}&=\{s_1,\cdots ,s_n\}\\
 i_l =j &\Leftrightarrow s_l \in t_{\alpha_j }.\end{align*}
\end{defi}
\begin{defi}Assume $A_1,\cdots ,A_k,B \in \mathfrak{B}$ and define
\begin{multline*} u(A_1,\cdots ,A_k,B): 
\{s,t \in \mathbb{R}^2,s<t\} \times \mathfrak{R}^k \\\mapsto u^t_s(A_1,\cdots ,A_k,B)(t_{\alpha_1},\cdots ,t_{\alpha_k}) \in
\mathfrak{B}\end{multline*}
for the case holding a.e., that all points
$$ \{s,t\} \cup \{t_i: i \in \alpha_1 +\cdots + \alpha_k\} $$
are different by
\begin{multline*}
u^t_s(A_1,\cdots ,A_k,B)((t_{\alpha_1},\cdots ,t_{\alpha_k})
= \1\{s<s_1<\cdots <s_n<t\}\\
\exp((t-s_n)B)A_{i_n}\exp((s_n-s_{n-1})B)A_{i_{n-1}}\\
\times\cdots\times A_{i_2}\exp((s_2-s_1)B)A_{i_1}\exp((s_1-s)B)
\end{multline*}
where
$$\Xi(t_{\alpha_1},\cdots ,t_{\alpha_k})=\{(s_1,i_1),\cdots ,(s_n,i_n) \}$$
with
$$s_1<\cdots <s_n.$$
\end{defi}
Define
\begin{eqnarray} 
\mathbf{e}:\mathfrak{R}^k &\to&\mathfrak{B}\nonumber\\
\mathbf{e}(w_1,\cdots ,w_k)& =& \left \{
 \begin{array}{ll} 
  1 & \mbox { if }w_1 = \cdots = w_k = \emptyset \\
0 &\mbox  { otherwise }\\
\end{array}
\right..
\end{eqnarray}
Write for short
$$\oint^j_{s,t} = \oint^j(\1_{]s,t[})$$
\begin{thm}
Assume$A_1,\cdots ,A_k,B \in \mathfrak{B}$ and
$$x:(t,w_1,\cdots ,w_k) \in \mathbb{R}\times \mathfrak{R}^k \mapsto x_t(w_1,\cdots ,w_k) \in \mathfrak{B}$$
to be a symmetric function in any of the variables $w_i$ and locally integrable. Consider for $t>s$ the equation
$$ x_t = \mathbf{e} + \sum_{j=1}^k A_j \oint^j_{s,t}x + \int _s^t Bx_u du .$$
Then 
$$ x_t = u_s^t(A_1,\cdots ,A_k,B) $$
is the unique solution of that equation.
\end{thm}
\begin{pf}
The proof is very similar to that of \cite{vW2005a}lemma 6.1. We include it for completeness.
 Using the mapping $\Xi$ of definition 2.2.2
we rewrite the equation for
$$ \xi = \{(s_1,i_1),\cdots ,(s_n,i_n)\};\;\; s_1<\cdots <s_n$$
and obtain
$$x_t(\xi) = 
\mathbf{e}(\xi)+ 
+ \sum_{l=1}^n A_{i_l}
x_{s_l}(\xi\setminus (s_l,i_l))\1\{s<s_l<t\} + B\int_s^tx_u(\xi)du
$$
We obtain
\begin{eqnarray*}x_t(\emptyset)&=& 1 + B \int _s^t x_u(\emptyset) du\\
x_t(\emptyset)&=&\exp((t-s)B)\end{eqnarray*}
We want to prove by induction, that $x_t(\xi)=0$ if $\{s_1,\cdots ,s_n\} \not\subset ]s,t[$.
Assume $n=1$ and $s_1 \notin ]s,t[$,then
$$x_t(\{(s_1,i_1)\}) = B \int _s^t x_u(\{(s_1,i_1)\}) du$$
hence $x_t(\{(s_1,i_1)\})=0$ .If $\{s_1,\cdots ,s_n\} \not\subset ]s,t[$, then at least one of the $s_i$,
either $s_1$ or $s_n$ are not in $]s,t[$. Assume $s_1 \notin ]s,t[$,then
$$ x_t(\xi) = \sum _{l=2}^n A_{i_l}\{s<s_l<t\}x_{s_l}(\xi\setminus (s_l,i_l)) + B\int _s^t x_u (\xi) du .$$
The first term vanishes, as $s_1 \in \{s_1,\cdots ,s_n\} \not\subset ]s,t[$ and therefor $ x_t(\xi)=0$.
Now if $s_1 \in \{s_1,\cdots ,s_n\} \subset ]s,t[$, then $x_{s_l}(\xi\setminus (s_l,i_l)) = 0$ for $l<n$ and 
$x_u(\xi)=0$ for $u<s_n$ and
$$ x_t(\xi)= A_{i_n}x_{s_n}(\{(s_1,i_1), \cdots ,(s_{n-1},i_{n-1})\}) + B\int _{s_n}^t x_u(\xi)du.$$
\end{pf}
In a similar way one proves
\begin{prop}  For $s<t$ the function
$$ s \mapsto y_s = u_s^t(A_1,\cdots ,A_k,B)$$
is the unique solution of the equation
$$ y_s = \mathbf{e} + \sum_{j=1}^k (\oint^j_{s,t}y)A_{i_j} -
 \int _s^t y_u du B.$$
\end{prop}
\begin{rem}
Use again the representation $\Xi$ and writing
$$  u_t^s(A_1,\cdots ,A_k,B)(t_{\alpha_1},\cdots ,t_{\alpha_k})= u_s^t(\xi)$$
with
$$ \xi = \{(s_1,i_1),\cdots ,(s_n,i_n)\};\;\; s_1<\cdots <s_n$$
and assume $s<r<t$ and $s_{j-1} <r<s_j$, then
$$ u_s^t(\xi)=u_r^t(\xi_2)u_s^r(\xi_1)$$
with
\begin{eqnarray*}\xi_1 &=& \{(s_1,i_1),\cdots ,(s_{j-1},i_{j-1})\}\\
\xi_2& =& \{(s_j,i_j),\cdots ,(s_n,i_n)\}.\end{eqnarray*}
\end{rem}

\subsection{Functions of Class $\mathcal{C}^1$}
\begin{defi}
Assume a  function
$$x:(t,t_{\alpha_1},\cdots ,t_{\alpha_k}) 
\in \mathbb{R}\times \mathfrak{R}^k \mapsto x_t(t_{\alpha_1},\cdots ,t_{\alpha_k}) \in \mathfrak{B}$$
symmetric in $t_{\alpha_1},\cdots ,t_{\alpha_k}$. Then $x$ is called of class $\mathcal{C}^0$ if
the function is locally intergable and is continuous in the subspace, where all points
 $t,t_i,i\in \alpha_1 + \cdots \alpha_k$
are different. We call $x$ of class $\mathcal{C}^1$ if it is  of class $\mathcal{C}^0$ and if on the same subspace the 
functions
\begin{eqnarray*}(\partial^cx)_t (t_{\alpha_1},\cdots ,t_{\alpha_k})&=& \frac{d}{dt}x_t(t_{\alpha_1},\cdots ,t_{\alpha_k)}\\
(R^j_\pm x)_t(t_{\alpha_1},\cdots ,t_{\alpha_k}) &=& x_{t\pm0}(t_{\alpha_1},\cdots t_{\alpha_{j-1}}
,t_{\alpha_j}+\{t\},t_{\alpha_{j+1}},\cdots ,t_{\alpha_k})\end{eqnarray*}
for $j=1,\cdots ,k$ exist and are  of class  $\mathcal{C}^0$. Here $d/dt= \partial^c$ is the usual derivative at 
the points of usual differentiability.
Put 
$$ D^jx = R^+_jx -R^j_-x$$
\end{defi}
\begin{prop}
For fixed $s$ is the function $t\mapsto u_s^t(A_i,B)$ and for fixed $t$ is the function $s\mapsto u_s^t(A_i,B)$
of class $\mathcal{C}^1$ and one has
\begin{eqnarray*}\partial^c_t u^t_s& =& Bu^t_s\\
(R^j_+u^._s)_t &=& A_j u^t_s \\
(R^j_-u^._s)_t &=& 0\\
\partial^c_s u^t_s& =& -u^t_sB\\
(R^j_+u^t_.)_s &=&  0\\
(R^j_-u^t_.)_s &=&  u^t_s A_j\\
\end{eqnarray*}
for $j=1,\cdots ,k$.
\end{prop}

We recall the definition of the Schwartz test functions on the real line. They form the space $C^\infty_c(\mathbb{R})$
of infinitely differentiable functions of compact support. 
\begin{defi}Assume $f$ to be locally integrable on $\mathbb{R}$, the Schwartz
derivative is the functional given by
$$(\partial f)(\varphi) = -\int f(t)\varphi '(t) dt $$
for Schwartz test functions $ \varphi$. If the functional is given by
$$(\partial f)(\varphi) = \int g(t)\varphi(t) dt,$$
where $g$ is locally integrable, we write
$$ g = \partial f.$$
\end{defi}
If $f$ is continuous differentiable except a finite set of points
$\{t_1,\cdots ,t_n\}$,then its Schwartz diffential is the measure
$$ \partial f = \partial ^c f  + \sum _{i=1}^n (f(t_i+0)-f(t_i-0)) \varepsilon_{t_i}(dt),$$
where $\partial ^c f$ is the usual derivation outside the jumping points and $\varepsilon_t$ is the point
measure in the point $t$.

We extend the notion of the Skorokhod integral to the weakly continuous
 measure valued function
$ \varepsilon :x \mapsto \varepsilon_x$.Hence
\begin{multline*}(\oint ^j(\varepsilon(dt))x) 
(t_{\alpha_1},\cdots ,t_{\alpha_k})=\\ \sum_{c \in \alpha_j} \varepsilon_{t_c}(dt)\;x_{t_c}(t_{\alpha_1},\cdots ,
t_{\alpha_{j-1}},t_{\alpha_j \setminus c},t_{\alpha_{j+1}},\cdots ,t_{\alpha_k})
\end{multline*}
This expression is scalarly defined, i.e. for any function $f$ with compact support in $\mathbb{R}$ we have
$$\int(\oint^j(\varepsilon(dt)x)f(t) = \oint^j(f)x$$
\begin{prop} If $x$ is of class $\mathcal{C}^1$, then its Schwartz derivative is
$$ (\partial x_t)(dt) =( \partial ^c x)_tdt + \sum _{j=1}^k \oint ^j(\varepsilon(dt) (D^jx).$$
\end{prop}
\begin{pf}
 We calculate
$$-\int_\Delta\cdots\int_\Delta dt_{\alpha_1}\cdots dt_{\alpha_k} g(t_{\alpha_1},\cdots ,t_{\alpha_k})
\int dt \varphi'(t) x_t(t_{\alpha_1},\cdots ,t_{\alpha_k})$$ 
where $g$ is a continuous function of compact support.It is sufficient to calculate the integral outside
the null set , where all the $t_i$ with $i \in \alpha_1+\cdots + \alpha_k$ are different. 
Using the representation (def 2.2.2)
$$ \Xi(t_{\alpha_1},\cdots ,t_{\alpha_k})= \xi = \{(s_1,i_1),\cdots ,(s_n,i_n)\}$$
with $s_1<\cdots <s_n$,we may write
\begin{multline*}-\int dt \varphi'(t) x_t(t_{\alpha_1},\cdots ,t_{\alpha_k}) = - \int dt \varphi'(t) x_t(\xi) \\ =
\int dt\varphi(t) \partial^cx_t(\xi) + \sum_{j=1}^n \varphi(t_j)(x_{s_j+0}(\xi)-x_{s_j-0}(\xi)).\end{multline*}
The second term equals
\begin{eqnarray*}&& \sum_{j=1}^k\sum_{c\in \alpha_j}\varphi(t_c)(x_{t_c+0}(t_{\alpha_1},\cdots ,t_{\alpha_k})-x_{t_c-0}
(t_{\alpha_1},\cdots ,t_{\alpha_k}))\\&=&
\sum_{j=1}^k\sum_{c\in \alpha_j}\varphi(t_c)(D^jx)_{t_c}(t_{\alpha_1},\cdots ,t_{a_j\setminus c},\cdots ,t_{\alpha_k})\\&=&
 \sum_{j=1}^k (\oint^j(\varphi)D^jx)(t_{\alpha_1},\cdots ,t_{\alpha_k})\end{eqnarray*}
 From there one obtains immediately the proposition.
 \end{pf}

 \section{Quantum stochastic processes of class $\mathcal{C}^1$}
\subsection{Definition and Fundamental Theorem}
 Recall definition 2.3.1 of functions of class $\mathcal{C}^1$
  and use instead of the index sets $\alpha_1,\cdots ,\alpha_k$ the index sets
 $\sigma,\tau,\upsilon$ and set $\mathfrak{B}=B(\mathfrak{k})$, where $\mathfrak{k}$ is a Hilbert space.
 Assume $x_t(\sigma,\tau.\upsilon)$ of class $\mathcal{C}^1$ and denote
 \begin{align*}
 (R^1_\pm x)_t(t_\sigma,t_\tau,t_\upsilon)&= x_{t\pm 0}(t_\sigma+\{t\},t_\tau,t_\upsilon)\\
 (R^0_\pm x)_t(t_\sigma,t_\tau,t_\upsilon)&= x_{t\pm 0}(t_\sigma,t_\tau+\{t\},t_\upsilon)\\
 (R^{-1}_\pm x)_t(t_\sigma,t_\tau,t_\upsilon)&= x_{t\pm 0}(t_\sigma,t_\tau,t_\upsilon + \{t\})\\
  (D^i x)_t &= (R^i_+x)_t-(R^i_-x)_t \end{align*}
  Recall the definition of the measure
  $$ \mathfrak{m}(\pi,\sigma,\tau,\upsilon,\varrho) = \langle a_\pi a^+_{\sigma+\tau}a_{\tau+\upsilon}a^+_\varrho \rangle
  \lambda_{\pi+\upsilon}$$
  and the  sesquilinear  form subsection 1.10
  $$ \langle f|\mathcal{B}(F)| g \rangle = \int  \mathfrak{m}(\pi,\sigma,\tau,\upsilon,\varrho)
 f^+(\pi)F(\sigma,\tau,\upsilon)g(\varrho).$$
 \begin{defi} If $x_t$ is of class $\mathcal{C}^1$, we call $\mathcal{B}(x_t)$ a quantum stochastic process of
 class $\mathcal{C}^1$.\end{defi}
 \begin{thm}
  If $x_t$ is of class $\mathcal{C}^1$, then the Schwartz derivative of $\langle f|\mathcal{B}(x_t)|g \rangle$
 for $f,g \in \mathcal{K}_s(\mathfrak{R},\mathfrak{k})$ is a locally integrable function
 \begin{multline*}\partial \langle f|\mathcal{B}(x_t)|g \rangle=
 \langle f|\mathcal{B}(\partial^c x_t)|g \rangle +\\
 \langle a(t)f|\mathcal{B}(D^1x_t)|g \rangle +
 \langle a(t)f|\mathcal{B}(D^0x_t)|a(t)g \rangle +
 \langle f|\mathcal{B}(D^{-1}x_t)|a(t)g \rangle\end{multline*}
 \end{thm}
 \begin{pf}
 Using the proof of proposition 2.3.2 we have 
 \begin{align*} \int\partial_t (\langle f|\mathcal{B}(x_t)|g \rangle)\varphi(t)
 &= -\int \langle f|\mathcal{B}(x_t)|g \rangle \varphi '(t) dt\\
 &= \int \mathfrak{m}f^+(\pi)\partial _tx_t(\sigma,\tau,\upsilon)g(\varrho)\varphi(t)dt\\
 &+\int \mathfrak{m}f^+(\pi)\sum_{c \in \sigma}x_{t_c}(\sigma\setminus c,\tau,\upsilon)g(\varrho)\varphi(t_c)\\ 
  &+\int \mathfrak{m}f^+(\pi)\sum_{c \in \tau}x_{t_c}(\sigma,\tau\setminus  c,\upsilon)g(\varrho)\varphi(t_c)\\ 
  &+\int \mathfrak{m}f^+(\pi)\sum_{c \in \upsilon}x_{t_c}(\sigma,\tau,\upsilon\setminus c)g(\varrho)\varphi(t_c).
  \end{align*} 
  Using the sum integra lemma we obtain for the last three terms
  \begin{align*}& \int  \langle a_\pi a^+_{\sigma+c+\tau}a_{\tau+\upsilon}a^+_\varrho \rangle
  \lambda_{\pi+\upsilon }f^+(\pi)(D^1 x)_c g(\varrho)\varphi(c)\\
  +& \int  \langle a_\pi a^+_{\sigma+\tau+c}a_{\tau+c+\upsilon}a^+_\varrho \rangle
  \lambda_{\pi+\upsilon }f^+(\pi)(D^1 x)_c g(\varrho)\varphi(c)\\
   +& \int  \langle a_\pi a^+_{\sigma+c+\tau}a_{\tau+\upsilon+c}a^+_\varrho \rangle
  \lambda_{\pi+\upsilon +c}f^+(\pi)(D^{-1} x)_c g(\varrho)\varphi(c)
  \end{align*}
 and from there the result .\end{pf}
 \subsection{Ito's Theorem} Recall subsection 1.10 and consider the measure 
 $$\mathfrak{m}(\pi,\sigma_1,\tau_1,\upsilon_1,\sigma_2,\tau_2,\upsilon_2,\varrho)
 = \langle a_\pi a^+_{\sigma_1+\tau_1}a_{\tau_1+\upsilon_1}a^+_{\sigma_2+\tau_2}a_{t_2+\upsilon_2} a^+_\varrho \rangle
 \lambda_{\pi + \upsilon_1 +\upsilon_2}.$$
 Assume $F,G:\mathfrak{R}^3 \to B(\mathfrak{k})$ to be $\lambda$-measurable and define
 \begin{multline*} \langle f| \mathcal{B}(F,G)| g \rangle \\=
 \int\mathfrak{m}(\pi,\sigma_1,\tau_1,\upsilon_1,\sigma_2,\tau_2,\upsilon_2,\varrho)
 f^+(\pi)F(\sigma_1,\tau_1,\upsilon_1)G(\sigma_2,\tau_2,\upsilon_3)g(\varrho)\end{multline*}
 provided the integral exists in norm. 
 \begin{thm}
 Assume $x_t,y_t$ to be of class $\mathcal{C}^1$ and that for
  $f,g \in \mathcal{K}_s(\mathfrak{R},\mathfrak{k})$  the sesquilinear forms $\langle f|\mathcal{B}(F_t,G_t)|g \rangle$
  exist in norm and $t \in \mathbb{R}\mapsto \langle f|\mathcal{B}(F_t,G_t)|g \rangle$ is locally integrable,
  where 
  $$F_t \in \{x_t,\partial^cx_t,R^1_\pm x_t,R^0_\pm x_t ,R^{-1}_\pm x_t\}$$
 $$G_t \in \{y_t,\partial^cy_t,R^1_\pm y_t,R^0_\pm y_t ,R^{-1}_\pm y_t\}$$
 is any of the functions. Then the Schwartz derivative of
  $\langle f|\mathcal{B}(x_t,y_t)|g \rangle$ is a locally integrable function and yields
  \begin{align*}\partial \langle f|\mathcal{B}(x_t,y_t)|g \rangle = &
    \langle f|\mathcal{B}(\partial^c x_t,y_t) +\mathcal{B}(f,\partial^c y_t) + I_{-1,+1,t}|g \rangle \\
  &+   \langle a(t)f|\mathcal{B}(D^1 x_t,y_t) +\mathcal{B}(f,D^1y_t) + I_{0,+1,t}|g \rangle \\
   &+  \langle a(t)f|\mathcal{B}(D^0 x_t,y_t) +\mathcal{B}(f,D^0 y_t) + I_{0,0,t}|a(t)g \rangle \\
  &+  \langle f|\mathcal{B}(D^{-1} x_t,y_t) +\mathcal{B}(f,D^{-1} y_t) + I_{-1,0,t}|a(t)g \rangle \end{align*}
  with
 $$I_{i,j,t} = \mathcal{B}(R^i_+x_t,R^j_+y_t)-\mathcal{B}(R^i_-x_t,R^j_-y_t)
 $$
 \end{thm}
 \begin{pf}
 Define
 \begin{equation*}N(\sigma,\tau,\upsilon) = \begin{cases}1 &\text{ if 
 $\{t_{\sigma +\tau+\upsilon}\}^\bullet $
  has no multiple point}\\
0 & \text{ otherwise }\end{cases}\end{equation*}
As $N$ is a Lebesgue nullfunction we find  for a Schwartz 
testfunction $\varphi$
\begin{multline*}\int \varphi(t)\partial\langle f|\mathcal{B}(x_t,y_t)|g \rangle =
-\int dt\varphi '(t)\langle f|\mathcal{B}(x_t,y_t)|g \rangle\\= 
-\int \lambda_c \varphi '(c) f^+(\pi)(1-N(\sigma_1,\tau_1,\upsilon_1))x_t(\sigma_1,\tau_1,\upsilon_1)\\
(1-N(\sigma_2,\tau_2,\upsilon_2))y_t(\sigma_2,\tau_2,\upsilon_2)g(\varrho)
 \mathfrak{m}(\pi,\sigma_1,\tau_1,\upsilon_1,\sigma_2,\tau_2,\upsilon_2,\varrho)\end{multline*}
 We perform at first the integral over $t_c$ and obtain 
 \begin{align*}-&\int \lambda_c \varphi '(c)x_t(\sigma_1,\tau_1,\upsilon_1)
y_t(\sigma_2,\tau_2,\upsilon_2)\\
= &\int dt \varphi(t) (\partial^c x_t(\sigma_1,\tau_1,\upsilon_1)y_t(\sigma_2,t_2,\upsilon_2)+
( x_t(\sigma_1,\tau_1,\upsilon_1)\partial^c y_t(\sigma_2,t_2,\upsilon_2))\\
&+ \sum_{t_c}
\bigl(x_{t_c+0}(\sigma_1,\tau_1,\upsilon_1)y_{t_c+0}(\sigma_2,\tau_2,\upsilon_2)
-x_{t_c-0}(\sigma_1,\tau_1,\upsilon_1)y_{t_c-0}(\sigma_2,\tau_2,\upsilon_2)\bigr)\end{align*}
where $t_c$ runs through all points of
$t_{\sigma_1+\tau_1+\upsilon_1}\cup t_{\sigma_2+\tau_2+\upsilon_2} $.
 Remark that the points of $ t_{\sigma_1+\tau_1+\upsilon_1}$and $t_{\sigma_2+\tau_2+\upsilon_2} $
 resp. are all different, but both sets might have common points. The sum over $t_c$ equals
 \begin{multline*}\sum_{c \in \sigma_1 +\tau_1 +\upsilon_1}
 (x_{t_c+0}(\sigma_1,\tau_1,\upsilon_1)-x_{t_c-0}(\sigma_1,\tau_1,\upsilon_1))
 y_{t_c-0}(\sigma_2,\tau_2,\upsilon_2)\\
 +\sum_{c \in \sigma_2 +\tau_2 +\upsilon_2}(x_{t_c-0}(\sigma_1,\tau_1,\upsilon_1)(y_{t_c+0}(\sigma_2,\tau_2,\upsilon_1)
 -y_{t_c-0}(\sigma_2,\tau_2,\upsilon_2))\end{multline*}
 as e.g.
  $$x_{t_c+0}(\sigma_1,\tau_1,\upsilon_1)-x_{t_c-0}(\sigma_1,\tau_1,\upsilon_1) =0$$
  for $t_c \notin t_{\sigma_1+\tau_1+\upsilon_1}$.
  
  We discuss  the integrals of the terms of the form
  $$\sum_{c \in \sigma_i +\tau_i +\upsilon_i}
 (x_{t_c\pm 0}(\sigma_1,\tau_1,\upsilon_1)(1-N(\sigma_1,\tau_1,\upsilon_1))
 (y_{t_c\pm 0}(\sigma_2,\tau_2,\upsilon_2)(1-N(\sigma_2,\tau_2,\upsilon_2))$$
 and  assume at first, that $f,g,x _t,y_t,\varphi$ are $\geq 0$ and regard
 \begin{align*}\int f(\pi)\sum_{c \in \sigma_1 +\tau_1 +\upsilon_1}&
 (x_{t_c+0}(\sigma_1,\tau_1,\upsilon_1)(1-N(\sigma_1,\tau_1,\upsilon_1))\\
 &(y_{t_c+0}(\sigma_2,\tau_2,\upsilon_2)(1-N(\sigma_2,\tau_2,\upsilon_2))g(\varrho)\varphi(t_c)\\
 &\mathfrak{m}(\pi,\sigma_1,\tau_1,\upsilon_1,\sigma_2,\tau_2,\upsilon_2,\varrho)\\
  =  I + II + III &
  \end{align*}
 splitting up the sum into three parts
 $$\sum_{c \in \sigma_1 +\tau_1 +\upsilon_1} = \sum_{c \in \sigma_1}+\sum_{c \in \tau_1}+\sum_{c \in \upsilon_1}$$
 We have using the sum integral lemma
 \begin{align*}
 I = &\int f(\pi) \sum_{c \in \sigma_1}(R^1_+x)_{t_c}(\sigma_1\setminus c,\tau_1,\upsilon_1)(1-N(\sigma_1,\tau_1,\upsilon_1))\\
 &y_{t_c+0}(\sigma_2,\tau_2,\upsilon_2)(1-N(\sigma_2,\tau_2,\upsilon_2))g(\varrho)\varphi(t_c)\\
 &\mathfrak{m}(\pi,\sigma_1,\tau_1,\upsilon_1,\sigma_2,\tau_2,\upsilon_2,\varrho)\\
 =& 
  \int f(\pi)(R^1_+x)_{t_c}(\sigma_1,\tau_1,\upsilon_1)(1-N(\sigma_1+c,\tau_1,\upsilon_1))\\
 &y_{t_c+0}(\sigma_2,\tau_2,\upsilon_2)(1-N(\sigma_2,\tau_2,\upsilon_2))g(\varrho)\varphi(t_c)\\
 &\langle a_\pi a^+_{\sigma_1+c+\tau_1}a_{\tau_1+\upsilon_1}a^+_{\sigma_2 +\tau_2}a_{\tau_2+\upsilon_2}
 a_\varrho^+ \rangle \lambda_{\pi+\upsilon_1+\upsilon_2}\\
 =& \int dt \varphi(t) \langle a(t)f|\mathcal{B}(R^1_+x_t,y_t)|g \rangle
 \end{align*}
 The integral over $N(1)$ and $N(2)$ vanishes and $y_{t+0} = y_t$ a.e. with respect to the integrating measure.
 
 In the same way
  \begin{align*}
 II=& 
  \int f(\pi)(R^0_+x)_{t_c}(\sigma_1,\tau_1,\upsilon_1)(1-N(\sigma_1,\tau_1+c,\upsilon_1))\\
 &y_{t_c+0}(\sigma_2,\tau_2,\upsilon_2)(1-N(\sigma_2,\tau_2,\upsilon_2))g(\varrho)\varphi(t_c)\\
 &\langle a_\pi a^+_{\sigma_1+\tau_1+c}a_{\tau_1+c+\upsilon_1}a^+_{\sigma_2 +\tau_2}a_{\tau_2+\upsilon_2}
 a^+_\varrho \rangle \lambda_{\pi+\upsilon_1+\upsilon_2}
\end{align*}
 The function
 $$y_{t+0}(\sigma,\tau,\upsilon)(1-N(\sigma,\tau,\upsilon))
  = \begin{cases} y_t(\sigma,\tau,\upsilon) & \text{ if } t_c \notin t_{\sigma+\tau+\upsilon}\\
 (R^1y)_t(\sigma\setminus b,\tau,\upsilon) & \text{ if } t =t_b, b\in \sigma\\
 (R^0y)_t (\sigma,\tau\setminus b,\upsilon)& \text{ if } t=t_b,b\in \tau\\
 (R^{-1}y)_t (\sigma,\tau,\upsilon\setminus b)& \text{ if } t=t_b,b\in \upsilon
 \end{cases}$$
 is defined everywhere. Recall that it vanishes,if  $ t_{\sigma+\tau+\upsilon}$ has multiple points.
 We have
 $$ a_c a^+_{\sigma_2 +\tau_2} = a^+_{\sigma_2 +\tau_2}a_c + \sum_{b\in\sigma_2}\varepsilon(c,b)
 a^+_{\sigma_2\setminus b +\tau_2}+ \sum_{b \in \tau_2}\varepsilon (c,b)a^+_{\sigma_2 +\tau_2\setminus b}$$
 and split the expression $II$ into three parts
 $$ II =  II_1 + II_2 +II_3. $$
 We obtain 
   \begin{align*}
 II_1=& 
  \int f(\pi)(R^0_+x)_{t_c}(\sigma_1,\tau_1,\upsilon_1)(1-N(\sigma_1,\tau_1+c,\upsilon_1))\\
 &y_{t_c+0}(\sigma_2,\tau_2,\upsilon_2)(1-N(\sigma_2,\tau_2,\upsilon_2))g(\varrho)\varphi(t_c)\\
 &\langle a_\pi a^+_{\sigma_1+\tau_1+c}a_{\tau_1+\upsilon_1}a^+_{\sigma_2 +\tau_2}a_{\tau_2+\upsilon_2+c}
 a^+_\varrho \rangle \lambda_{\pi+\upsilon_1+\upsilon_2}\\
 = & \int dt \varphi(t) \langle a(t)f|\mathcal{B}(R^0_+x_t,y_t)|a(t)g \rangle
\end{align*}
 Using the sum integral lemma
    \begin{align*}
 II_2=& 
  \int f(\pi)(R^0_+x)_{t_c}(\sigma_1,\tau_1,\upsilon_1)(1-N(\sigma_1,\tau_1+c,\upsilon_1))\\
 &\sum _{b \in \sigma_2}\varepsilon(c,b)
 y_{t_c+0}(\sigma_2,\tau_2,\upsilon_2)(1-N(\sigma_2,\tau_2,\upsilon_2))g(\varrho)\varphi(t_c)\\
 &\langle a_\pi a^+_{\sigma_1+\tau_1+c}a_{\tau_1+\upsilon_1}a^+_{\sigma_2\setminus b +\tau_2}a_{\tau_2+\upsilon_2}
 a^+_\varrho \rangle \lambda_{\pi+\upsilon_1+\upsilon_2}\\
 = &\int f(\pi)(R^0_+x)_{t_c}(\sigma_1,\tau_1,\upsilon_1)(1-N(\sigma_1,\tau_1+c,\upsilon_1))\\
 &y_{t_c+0}(\sigma_2 +c,\tau_2,\upsilon_2)(1-N(\sigma_2+c,\tau_2,\upsilon_2))g(\varrho)\varphi(t_c)\\
 &\langle a_\pi a^+_{\sigma_1+\tau_1+c}a_{\tau_1+\upsilon_1}a^+_{\sigma_2+\tau_2}a_{\tau_2+\upsilon_2}
 a^+_\varrho \rangle \lambda_{\pi+\upsilon_1+\upsilon_2}\\
 = & \int dt \varphi(t) \langle a(t)f|\mathcal{B}(R^0_+x_t,R^1_+y_t)|g \rangle
\end{align*}
 and again
  \begin{align*}
 II_3=& 
  \int f(\pi)(R^0_+x)_{t_c}(\sigma_1,\tau_1,\upsilon_1)(1-N(\sigma_1,\tau_1+c,\upsilon_1))\\
 &\sum _{b \in \tau_2}\varepsilon(c,b)
 y_{t_c+0}(\sigma_2,\tau_2,\upsilon_2)(1-N(\sigma_2,\tau_2,\upsilon_2))g(\varrho)\varphi(t_c)\\
 &\langle a_\pi a^+_{\sigma_1+\tau_1+c}a_{\tau_1+\upsilon_1}a^+_{\sigma_2 +\tau_2\setminus b}a_{\tau_2+\upsilon_2}
 a^+_\varrho \rangle \lambda_{\pi+\upsilon_1+\upsilon_2}\\
 = &\int f(\pi)(R^0_+x)_{t_c}(\sigma_1,\tau_1,\upsilon_1)(1-N(\sigma_1,\tau_1+c,\upsilon_1))\\
 &y_{t_c+0}(\sigma_2 +c,\tau_2,\upsilon_2)(1-N(\sigma_2,\tau_2+c,\upsilon_2))g(\varrho)\varphi(t_c)\\
 &\langle a_\pi a^+_{\sigma_1+\tau_1+c}a_{\tau_1+\upsilon_1}a^+_{\sigma_2+\tau_2}a_{\tau_2+c+\upsilon_2}
 a^+_\varrho \rangle \lambda_{\pi+\upsilon_1+\upsilon_2}\\
 = & \int dt \varphi(t) \langle a(t)f|\mathcal{B}(R^0_+x_t,R^0_+y_t)|a(t)g \rangle
 \end{align*}
 Similar
 \begin{align*}
  III=& 
  \int f(\pi)(R^{-1}_+x)_{t_c}(\sigma_1,\tau_1,\upsilon_1)(1-N(\sigma_1,\tau_1,\upsilon_1+c))\\
 &y_{t_c+0}(\sigma_2,\tau_2,\upsilon_2)(1-N(\sigma_2,\tau_2,\upsilon_2))g(\varrho)\varphi(t_c)\\
 &\langle a_\pi a^+_{\sigma_1+\tau_1}a_{\tau_1+\upsilon_1+c}a^+_{\sigma_2 +\tau_2}a_{\tau_2+\upsilon_2}
 a^+_\varrho \rangle \lambda_{\pi+\upsilon_1+c+\upsilon_2}\\
 =  & III_1 + III_2 +III_3
\end{align*}
 with
  \begin{align*}
  III_1=& 
  \int f(\pi)(R^{-1}_+x)_{t_c}(\sigma_1,\tau_1,\upsilon_1)(1-N(\sigma_1,\tau_1,\upsilon_1+c))\\
 &y_{t_c+0}(\sigma_2,\tau_2,\upsilon_2)(1-N(\sigma_2,\tau_2,\upsilon_2))g(\varrho)\varphi(t_c)\\
 &\langle a_\pi a^+_{\sigma_1+\tau_1}a_{\tau_1+\upsilon_1}a^+_{\sigma_2 +\tau_2}a_{\tau_2+\upsilon_2+c}
 a^+_\varrho \rangle \lambda_{\pi+\upsilon_1+c+\upsilon_2}\\
 = & \int dt \varphi(t) \langle f|\mathcal{B}(R^{-1}_+x_t,y_t)|a(t)g \rangle
 \end{align*}
 and
   \begin{align*}
  III_2=& 
  \int f(\pi)(R^{-1}_+x)_{t_c}(\sigma_1,\tau_1,\upsilon_1)(1-N(\sigma_1,\tau_1,\upsilon_1+c))\\
 &\sum_{b \in \sigma_2}\varepsilon(c,b)
 y_{t_c+0}(\sigma_2,\tau_2,\upsilon_2)(1-N(\sigma_2,\tau_2,\upsilon_2))g(\varrho)\varphi(t_c)\\
 &\langle a_\pi a^+_{\sigma_1+\tau_1}a_{\tau_1+\upsilon_1}a^+_{\sigma_2\setminus b +\tau_2}a_{\tau_2+\upsilon_2}
 a^+_\varrho \rangle \lambda_{\pi+\upsilon_1+c+\upsilon_2}\\
  =& \int f(\pi)(R^{-1}_+x)_{t_c}(\sigma_1,\tau_1,\upsilon_1)(1-N(\sigma_1,\tau_1,\upsilon_1+c))\\
 &
 (R^1_+y)_{t_c}(\sigma_2,\tau_2,\upsilon_2)(1-N(\sigma_2+c,\tau_2,\upsilon_2))g(\varrho)\varphi(t_c)\\
 &\langle a_\pi a^+_{\sigma_1+\tau_1}a_{\tau_1+\upsilon_1}a^+_{\sigma_2+\tau_2}a_{\tau_2+\upsilon_2}
 a^+_\varrho \rangle \lambda_{\pi+\upsilon_1+c+\upsilon_2}\\
 = & \int dt \varphi(t) \langle f|\mathcal{B}(R^{-1}_+x_t,R^1_+y_t)|g \rangle
 \end{align*}
 and finally
    \begin{align*}
  III_3=& 
  \int f(\pi)(R^{-1}_+x)_{t_c}(\sigma_1,\tau_1,\upsilon_1)(1-N(\sigma_1,\tau_1,\upsilon_1+c))\\
 &\sum_{b \in \tau_2}\varepsilon(c,b)
 y_{t_c+0}(\sigma_2,\tau_2,\upsilon_2)(1-N(\sigma_2,\tau_2,\upsilon_2))g(\varrho)\varphi(t_c)\\
 &\langle a_\pi a^+_{\sigma_1+\tau_1}a_{\tau_1+\upsilon_1}a^+_{\sigma_2 +\tau_2\setminus b}a_{\tau_2+\upsilon_2}
 a^+_\varrho \rangle \lambda_{\pi+\upsilon_1+c+\upsilon_2}\\
  =& \int f(\pi)(R^{-1}_+x)_{t_c}(\sigma_1,\tau_1,\upsilon_1)(1-N(\sigma_1,\tau_1,\upsilon_1+c))\\
 &
 (R^1_+y)_{t_c}(\sigma_2,\tau_2,\upsilon_2)(1-N(\sigma_2,\tau_2+c,\upsilon_2))g(\varrho)\varphi(t_c)\\
 &\langle a_\pi a^+_{\sigma_1+\tau_1}a_{\tau_1+\upsilon_1}a^+_{\sigma_2+\tau_2}a_{\tau_2 +c+\upsilon_2}
 a^+_\varrho \rangle \lambda_{\pi+\upsilon_1+c+\upsilon_2}\\
 = & \int dt \varphi(t) \langle f|\mathcal{B}(R^{-1}_+x_t,R^0_+y_t)|a(t)g \rangle
 \end{align*}
 Collecting all terms
 \begin{multline*}I + II + III =
 $$ \int dt \varphi(t)\biggl(\langle a(t)f|\mathcal{B}(R^1_+x_t,y_t)|g\rangle+
 \langle a(t)f|\mathcal{B}(R^1_+x_t,y_t)|a(t)g\rangle +\\ \langle a(t)f|\mathcal{B}(R^0_+x_t,R^1_+y_t)|a(t)g\rangle
 +\langle a(t)f|\mathcal{B}(R^0_+x_t,R^0_+y_t)|a(t)g\rangle\\
 +\langle f|\mathcal{B}(R^{-1}_+x_t,y_t)|g\rangle
 +\langle f|\mathcal{B}(R^{-1}_+x_t,R^1_+y_t)|g\rangle
 +\langle f|\mathcal{B}(R^{-1}_+x_t,y_t)|g\rangle 
 \biggl)\end{multline*}
 The assumptions of our theorem guarantee that all expressions exist and we may extend the formuals to vector
 and operator valued functions.
 By analogous calculations
  \begin{align*}\int f(\pi)^+\sum_{c \in \sigma_1 +\tau_1 +\upsilon_1}&
 (x_{t_c\pm 0}(\sigma_1,\tau_1,\upsilon_1)(1-N(\sigma_1,\tau_1,\upsilon_1))\\
 &(y_{t_c+0}(\sigma_2,\tau_2,\upsilon_2)(1-N(\sigma_2,\tau_2,\upsilon_2))g(\varrho)\varphi(t_c)\\
 &\mathfrak{m}(\pi,\sigma_1,\tau_1,\upsilon_1,\sigma_2,\tau_2,\upsilon_2,\varrho)\\
  =& \int dt\varphi(t) K_{\pm,+}(t)\\
  \end{align*}
  and
    \begin{align*}\int f(\pi)^+\sum_{c \in \sigma_2 +\tau_2 +\upsilon_2}&
 (x_{t_c- 0}(\sigma_1,\tau_1,\upsilon_1)(1-N(\sigma_1,\tau_1,\upsilon_1)\\
 &(y_{t_c \pm 0}(\sigma_2,\tau_2,\upsilon_2)(1-N(\sigma_2,\tau_2,\upsilon_2)g(\varrho)\varphi(t_c)\\
 &\mathfrak{m}(\pi,\sigma_1,\tau_1,\upsilon_1,\sigma_2,\tau_2,\upsilon_2,\varrho)\\
  =& \int dt\varphi(t) K_{-,\pm}(t)\\
  \end{align*}
 We have 
  $$ K_{\pm,+} = K_{\pm,+}^{(1)} + K_{\pm,+}^{(2)}$$
  $$ K_{-,\pm} = K_{-,\pm}^{(1)} + K_{-,\pm}^{(2)}$$
  with
  $$K_{\pm,+}^{(1)} = \langle a(t)f|\mathcal{B}(R^1_\pm x_t,y_t)|g \rangle
  +\langle a(t)f|\mathcal{B}(R^0_\pm x_t,y_t)|a(t)g \rangle
  +\langle f|\mathcal{B}(R^{-1}_\pm x_t,y_t)|a(t)g \rangle$$
   $$K_{-,\pm}^{(1)} = \langle a(t)f|\mathcal{B}( x_t,R^1_\pm y_t)|g \rangle
  +\langle a(t)f|\mathcal{B}( x_t,R^0_+y_t)|a(t)g \rangle
  +\langle f|\mathcal{B}(x_t,R^{-1}_+y_t)|a(t)g \rangle$$
  and
  \begin{multline*}K^{(2)}_{\pm,\pm} =  \langle a(t)f|\mathcal{B}(R^0_\pm x_t,R^1_\pm y_t|g \rangle
 +\langle a(t)f|\mathcal{B}(R^0_\pm x_t,R^0_\pm y_t|a(t)g \rangle\\
  + \langle f|\mathcal{B}(R^{-1}_\pm x_t,R^0_\pm y_t|a(t)g\rangle
 + \langle f|\mathcal{B}(R^{-1}_\pm x_t,R^1_\pm y_t|g \rangle.\end{multline*}
 From there one achieves the final result without big difficulty.
 \end{pf}
 \section{ A Quantum Stochastic Differential Equation}
 \subsection{Formulation of the Equation}
 We shall investigate the quantum stochastic differential equation which reads in the Hudson-Parthasarathy
 calculus
 $$d_t U^t_s = A_1 dB^+_tU^t_s +A_0d\Lambda_tU^t_s +A_{-1}dB_tU^t_s + BU^t_sdt,\;U_s^s = 1$$
 where $A_1,A_0,A_{-1},B$ are Operators in $B(\mathfrak{k})$.
 
 In his white noise calculus Accardi formulates it as \emph{normal ordered} equation
  $$ dU^t_s/dt = A_1 a^+_t U^t_s + A_0 a^+_t U^t_s a_t + A_{-1}U^t_s a_t + BU^t_s$$
  
  Our approach is very similar to Accardi's one. We understand $U^t_s$ as a sesquilinear form over
  $\mathcal{K}_s(\mathfrak{R})$ given by
  $$\langle f|U^t_s|g\rangle = \int f^+(\pi) u^t_s(\sigma,\tau,\upsilon)g(\varrho)
  a_\pi a^+_{\sigma+\tau}a_{\tau+\upsilon}a^+_\varrho \lambda_{\pi+\varrho}$$
  where $u^t_s$ is locally integrable in all four variables.
  We formulate the differential equation in the weak sense
 \begin{multline*}(d/dt)\langle f|U^t_s|g\rangle \\= \langle a(t)f|A_1U^t_s|g\rangle +\langle a(t)f|A_0U^t_s|a(t)g\rangle +
 \langle f|A_{-1}U^t_s|a(t)g\rangle +\langle f|BU^t_s|g\rangle\end{multline*}
 or better as integral equation
 \begin{multline*}\langle f|U^t_s|g\rangle = \langle f|g\rangle + \int_s^tdr \langle a(r) f|A_1U^r_s|g\rangle +
  \int_s^tdr \langle a(r) f|A_0U^r_s|a(r)g\rangle \\+ \int_s^tdr \langle  f|A_{-1}U^r_s|a(r)g\rangle +
   \int_s^tdr \langle  f|BU^r_s|g\rangle \tag{$*$}\end{multline*}
   for $ t \ge s$
   We shall show that this equation has a uniques solution, which can be given explicitely.
   \subsection{Existence and Unicity of the Solution}
   \begin{lem}
   The equation $(*)$ is equivalent to the Skorokhod integral equation $(**)$
   \begin{equation*}u^t_s = \mathbf{e} + A_1 \oint ^1_{s,t}u_s^.+ A_0 \oint^0_{s,t} u_s^. + A_{-1}\oint^{-1}_{s,t}u_s^.
   + B\int^t_s dr u^r_s \tag{$**$}\end{equation*}
   \end{lem}
   \begin{pf}
   Cosider e.g.the term
 \begin{align*}&\int_s^tdr \langle a(r) f|A_0U^r_s|a(r)g\rangle \\&=
 \int\1_{[s,t]}(t_c)f^+(\omega+\sigma+\tau+c)A_0 u_s^{t_c}(\sigma,\tau,\upsilon)g(\omega+\tau+\upsilon+c)
 \lambda_{\omega+\sigma+\tau+\upsilon+c}
  \\&=\int\sum_{c\in \tau}\1_{[s,t]}(t_c)f^+(\omega+\sigma+\tau)A_0 u_s^{t_c}(\sigma,\tau\setminus c,\upsilon)
  g(\omega+\tau+\upsilon)
 \lambda_{\omega+\sigma+\tau+\upsilon}
  \\&=\int f^+(\omega+\sigma+\tau)A_0 (\oint ^0_{s,t}u_s^.)(\sigma,\tau,\upsilon)
  g(\omega+\tau+\upsilon)
 \lambda_{\omega+\sigma+\tau+\upsilon}\end{align*}
 Remark that the function $u^t_s(\sigma,\tau,\upsilon)$ is determined by the sesquilinear 
 form $\langle f|U^t_s|g \rangle$ Lebesgue almost 
 everywhere. 
 \end{pf}
 Applying theorem 2.2.1 we obtain immediately
 \begin{thm} Equation $(*)$ has a unique solution namely
 \[U^t_s= \mathcal{B}(u^t_s(A_1,A_0,A_{-1};B) \tag{$**$}.\]
 \end{thm}
 \subsection{ Invarianz of $\Gamma_k$.Continuity}
 \begin{defi} We define the Fock space
$$ \Gamma = L^2_s(\mathfrak{R},\mathfrak{k},\lambda) $$
of all symmetric square integrable functions with respect to Lebesgue measure from 
$\mathfrak{R}$ to $\mathfrak{k}.$ 
If $f$ is a mesurable function on $\mathfrak{R}$ define the operator $N$ by
$(Nf)(w) = (\#w)f(w)$ and define $\Gamma_k$ as the space of those measurable
 symmetric functions
from  $\mathfrak{R}$ to $\mathfrak{k}$ , for which 
$$\int_{\Delta(\mathfrak{R})} \langle f(w)|(N+1)^k\,f(w)\rangle\,dw < \infty $$

We denote by $\|.\|_{\Gamma_k}$ the corresponding norm.
We write for short
$$ \mathcal{K} = \mathcal{K}_s(\mathfrak{R},\mathfrak{k})$$
for the space of all symmetric continuous functions from  $\mathfrak{R}$ to $\mathfrak{k}$
with compact support.
Call $\mathcal{K}^{(n)}$,resp. $\Gamma^{(n)}$ the subspace, where $f(w)= 0$ for $\#w > n$.
\end{defi}
We extend the notions of $a$  and $a^+$. We define $a(\varphi\lambda)f = a(\varphi)f$ and
$a^+(\varphi)f$ for $\varphi\in L^2(\mathbb{R})$ and $f \in \Gamma$.We have the well known relations
\begin{align*} 
a(\varphi):\Gamma^{(n)}& \to \Gamma^{(n-1)}& \|a(\varphi)f\|_\Gamma& \leq \sqrt{n}\|\varphi\|_{L^2} \|f\|_\Gamma 
\\
 a(\varphi)^+:\Gamma^{(n)}&\to \Gamma^{(n+1)}&\|a(\varphi)^+f\|_\Gamma &\leq \sqrt{n+1}\|\varphi\|_{L^2} \|f\|_\Gamma 
 \end{align*}
 \begin{lem}
 Assume a Lebesgue measurable function $\mathfrak{R}\to\mathfrak{k}$, then
 \begin{equation*}\int\lambda_{\xi+\omega}\1\{\#\xi = k\}\,\, \|f(\xi+\omega)\|^2= \left\langle f\left| \binom {N}{k}
 \right|f\right\rangle .\end{equation*}
 \end{lem}
 \begin{pf} The LHS of the last equation equals
 $$ \int \lambda_\omega \sum_{\xi \in \omega} \|f(\omega\|^2 \1\{\#\xi = k \} =  \int \lambda_\omega \binom {\#\omega}{k}
 \|f(\omega)\|^2 :$$
 \end{pf}
 \begin{prop}
 Assume 
 $$ U^t_s= \mathcal{O}(u^t_s(A_1,A_0,A_{-1};B)).$$
 Then there exist constants $C_{n,k}(t-s)$ such that for $ f \in \mathcal{K}^{(n)}$
 $$ \|U^t_s f\|_{\Gamma_k} \le C_{n,k}(t-s) \|f\|_\Gamma $$
 Furthermore for $ t \downarrow s$ and $f \in \mathcal{K}$
 $$\|U^t_s f -f\|_{\Gamma_k} \to 0.$$
 \end{prop}
 \begin{pf}
 Define
 $$ C = \max ( \|A_i\|,i= 1,0,-1,\|B\|),$$
 then
 \begin{multline*} \|u^t_s(\sigma,\tau,\upsilon)\|  \le e^{C(t-s)} C^{\#\sigma+\#\tau+\#\upsilon}\1 \{t_{\sigma+\tau+\upsilon} 
 \subset [s,t]\}\\ = e^{C(t-s)}e(\chi)(\sigma+\tau+\upsilon)\end{multline*}
 with 
 $$ \chi(r)= C \1_{[s,t]}(r).$$
 We have using 1.10
 \begin{multline*}\|\mathcal{O}(u^t_s)f(\omega) \| \\\le \sum _{\sigma \subset \omega}\sum_{\tau \subset \omega\setminus s}
 \int \lambda_\upsilon \exp(C(t-s))(e(\chi)(\sigma+\tau+\upsilon)\Vert f(\omega\setminus\sigma+\upsilon)\Vert
\\ = \exp(C(t-s))(R^t_s\;S^t_s\;T^t_s\;\|f\|)(\omega).\end{multline*}
For $g \in \mathcal{K}^{(n)}(\mathfrak{R},\mathbb{R}), g \ge 0$ we have
$$ (T^t_s g)(\omega)= \int  \lambda_\upsilon e(\chi)(\upsilon)g(\omega+\upsilon)= (\exp(a(\chi))g)(\omega).$$
As $T^t_s: \mathcal{K}^{(n)}\to \mathcal{K}^{(n)}$ , we may estimate the $\Gamma_k$-norm by the 
$\Gamma$-norm.We have
$$\Vert T^t_s g \Vert \le \sum_{l=0}^n (1/l!)\sqrt{n(n-1)\cdots (n-l+1)} C^l(t-s)^{l/2} \Vert g \Vert_\Gamma$$
as
$$ \Vert \chi \Vert_{L^2}= C \sqrt {t-s}.$$
Further more we have $$S^t_s: \mathcal{K}^{(n)} \to \Gamma^{(n)},\;\;\;\;
 (S^t_s g)(\omega) = \sum_{\tau \subset \omega} e(\chi)(\tau)g(\omega) = e(1+\chi)(\omega)g(\omega)$$
and
$$\|S^s_t g\|_\Gamma \le (1+C)^n \|g\|_\Gamma.$$ 
Again
\begin{align*}R^t_s:\Gamma^{(n)}& \to \Gamma&(R^s_tg)(\omega)&=\sum_{\sigma\subset\omega}
e(\chi)(\sigma)g(\omega\setminus\sigma)=
(\exp(a^+(\chi)g)(\omega) \end{align*}
and
$$ \|R^t_sg\| \le \sum_{l=0}^\infty (1/l!)\sqrt{(n+1)(n+2)\cdots (n+l)}C^l(t-s)^{l/2}\|g\|_\Gamma=
c_n(t-s)\|g\|_\Gamma $$
We calculate
$$ (R^s_tg)(\omega+\xi)= \sum_{\sigma\subset\omega+\xi}e(\chi)(\sigma)((\omega+\xi)\setminus\sigma)
= \sum_{\substack{\xi_1+\xi_2=\xi\\\omega_1+\omega_2=\omega}}e(\chi)(\xi_1+\omega_1)g(\omega_2+\xi_2)$$
$$= \sum e(\xi)(\xi_1)e(\chi)(\omega_1)(a_{\xi_2}g)(\omega_2)
=\sum_{\xi_1 + \xi_2 = \xi}e(\chi)(\xi_1)(\exp(a^+(\chi)a_{\xi_2}g)(\omega)$$
and
$$\left\langle R^s_tg\left|\binom{N}{k}\right| R^s_t g \right\rangle=
\int\lambda_{\omega+\xi}((R^s_tg)(\omega+\xi))^2 \1\{\#\xi=k\}$$
$$\le 2^k\int \lambda_{\omega+\xi}\sum_{\xi_1+\xi_2=\xi}e(\chi^2)(\xi_1)(\exp(a^+(\chi)a_{\xi_2}g)(\omega))^2
\1\{\#\xi=k\}$$
$$ \le 2^k\sum_{k_1+k_2=k}\int_{\#\xi_=k_1} \lambda_{\xi_1}e(\chi^2)(\xi_1)( c_{n-k_2}(t-s))^2
\int_{\#\xi_2=k_2}\lambda_{\xi_2+\omega}(a_{\xi_2}g)(\omega))^2$$
$$ \le 2^k\sum_{k_1+k_2=k}(1/k_1!)C^{2k_1}(t-s)^{k_1}(c_{n-k_2}(t-s))^2\left\langle g \left | \binom{N}{k_2}\right|g\right)
$$
and
$$\left\langle g \left | \binom{N}{k_2}\right|g\right)\le \binom {n}{k_2}\|g\|_\Gamma^2.$$
From there follows the first assertion of the proposition.

We investigate the second assertion.As
$$(\mathcal{O}(u^t_sf-f)(\omega) = \sum_{\sigma \subset\omega}
\sum_{\tau\subset\omega\setminus\sigma}\int \lambda_\upsilon u^t_s(\sigma,\tau,\upsilon)
\1\{ \sigma+\tau+\upsilon\ne \emptyset\}f(\omega\setminus\sigma+\upsilon)$$
we may estimate the norm by
\begin{multline*} \sum_{\sigma \subset\omega}
\sum_{\tau\subset\omega\setminus\sigma}\int \lambda_\upsilon \exp(C(t-s))e(\chi)(\sigma+\tau+\upsilon)
\|f(\omega\setminus\sigma+\upsilon)\| \1\{\sigma+\tau+\upsilon \ne \emptyset\}\\
=\exp (C(t-s))(R^t_sS^t_sT^t_s-1)\|f\|)(\omega).\end{multline*}
We have
\begin{multline*}\|T^t_sg-g\|_\Gamma\\ \le  \sum_{l=1}^n (1/l!)\sqrt{n(n-1)\cdots (n-l+1)} C^l(t-s)^{l/2} \Vert g \Vert_\Gamma
=O(\sqrt{t-s})\end{multline*}
We have
\begin{multline*}((S^t_s-1)g)(\omega)=((e(1+\chi)(\omega)-1)g(\omega)\\= \sum_{c \in \omega}\chi(c)e(1+\chi)(\omega\setminus c)g(\omega)
\le \sum_{c \in \omega}\chi(c)(1+C)^{n-1}g(\omega).\end{multline*}
As $f\in \mathcal{K}^n$, there exists a compact interval $K \subset \mathbb{R},\;[s,t]\subset K$, such that
$g(\omega) \le e(\1_K)(\omega)$ for $\#\omega \le n$ , if $g(\omega) \le 1$ for all $\omega$. We have 
 $$\sum_{c \in \omega}\chi(c)(1+C)^{n-1}g(\omega)\le \sum_{c \in \omega} \chi(c)(1+C)^n e(\1_K)(\omega\setminus c)$$
 The norm may be estimated by
 $$\sqrt{n+1}\sqrt{t-s}(1+C)^n \exp(|K|/2).$$
 We have 
 $$((R^s_t-1)g)(\xi+\omega)= \sum_{\xi_1+\xi_2=\xi}e(\chi)(\xi_1)(\exp(a^+(\chi))a_{\xi_2}g)(\omega)-a_\xi g(\omega)$$
 $$= \sum_{\substack{\xi_1+\xi_2=\xi\\\xi_1\ne \emptyset}}e(\chi)(\xi_1)(\exp(a^+(\chi))a_{\xi_2}g)(\omega)
 +(\exp(a^+(\chi))-1)a_\xi)g(\omega)= I + II.$$
 The $\Gamma$-norm square of $I$ can be estimated by
 $$2^k\sum_{\substack{k_1+k_2 =k \\ k_1 \ne 0}} 
(1/k_1!)C^{2k_1}(t-s)^{k_1}(c_{n-k_2}(t-s))^2\binom{n}{k_2}\|g\|_\Gamma^2 =O(t-s).$$
The norm of $II$ can be estimated by 
\begin{multline*}
\sum_{l=1}^\infty (1/l!)\sqrt{(n+1)(n+2)\cdots (n+l)}C^l(t-s)^{l/2}(\int \lambda_{\xi+\omega}g(\xi+\omega)^2)^{1/2}\\
=O(\sqrt{t-s})\end{multline*}
From these results one obtains the second assertion of the proposition easily.
\end{pf}
\subsection{Time Consecutive Intervals}
We consider again $u^t_s(A_1,A_0,A_{-1};B)$. We showed in 4.4, that $\mathcal{O}(u^t_s):\mathcal{K}^{(n)} \to
\Gamma$ is bounded. So $\mathcal{B}(u^t_r,u^r_s)$ exists (see 1.10).
\begin{prop}
For $ s<r<t$ we have
$$\mathcal{B}(u^t_r,u^r_s)=\mathcal{B}(u^t_s).$$
\end{prop}
\begin{pf}
If $\{r,s,t,T_{\sigma+\tau+\tau}\}^\bullet$ is without multiple points, then we showed in remark 2.2.1
$$u^t_s(\sigma,\tau,\upsilon) = u^t_r(\sigma_2',\tau_2',\upsilon_2')u^r_s(\sigma_1',\tau_1',\upsilon_1')$$
with
$$ \sigma_2'=\{ c \in \sigma: r < t_c <t\},\;\; \sigma_1'=\{ c \in \sigma: r < t_c <t\}$$
etc..Hence
\begin{equation*}u^t_s(\sigma,\tau,\upsilon) = \sum_{\substack{\sigma_1+\sigma_2=\sigma\\\tau_1+\tau_2=\tau\\\upsilon_1+\upsilon_2
=\upsilon}} u^t_r(\sigma_2,\tau_2,\upsilon_2)u^r_s(\sigma_1,\tau_1,\upsilon_1)\tag{$*$}\end{equation*}
as from this terms only those terms survive, where $\sigma_1=\sigma_1',\sigma_2=\sigma_2',\cdots .$ 
We have 
\begin{multline*}\langle f|\mathcal{B}(u^t_r,u^r_s)|g \rangle = \int f^+(\pi)u^t_r(\sigma_2,\tau_2,\upsilon_2)
u^r_t(\sigma_1,\tau_1,\upsilon_1)g(\varrho)\\
\langle a_\pi a^+_{\sigma_2+\tau_2}a_{t_2+\upsilon_2}a^+_{\sigma_1+\tau_1}a_{\tau_1+\upsilon_1}a^+_\varrho\rangle
\lambda_{\pi+\upsilon_1+\upsilon_2}\end{multline*}
We claim that the last term equals
 \begin{multline*}\int f^+(\pi)u^t_r(\sigma_2,\tau_2,\upsilon_2)
u^r_t(\sigma_1,\tau_1,\upsilon_1)g(\varrho)\\
\langle a_\pi a^+_{\sigma_2+\tau_2}a^+_{\sigma_1+\tau_1}a_{\tau_2+\upsilon_2}a_{\tau_1+\upsilon_1}a^+_\varrho\rangle
\lambda_{\pi+\upsilon_1+\upsilon_2}\tag{$**$}\end{multline*}
Take $c\in \tau_2+\upsilon_2$, e.g.$c \in \upsilon_2$, then
$$ a_{\tau_2+\upsilon_2}a^+_{\sigma_1+\tau_1} =  a_{\tau_2+\upsilon_2 \setminus c}
a^+_{\sigma_1+\tau_1}a_c + \sum_{b\in\sigma_2+\tau_2} \varepsilon(c,b)a_{\tau_2+\upsilon_2\setminus c}
a^+_{(\sigma_1+\tau_1)\setminus c}.$$
But
\begin{multline*}\int f^+(\pi)u^t_r(\sigma_2,\tau_2,\upsilon_2)
u^r_t(\sigma_1,\tau_1,\upsilon_1)g(\varrho)\\\varepsilon(c,b)
\langle a_\pi a^+_{\sigma_2+\tau_2}a_{t_2+\upsilon_2\setminus c}
a^+_{(\sigma_1+\tau_1)\setminus b}a_{\tau_1+\upsilon_1}a^+_\varrho\rangle  
\lambda_{\pi+\upsilon_1+\upsilon_2}=0\end{multline*}
as
$$\int\lambda_c \varepsilon(c,b)\1\{r<t_c<t\}\1\{s<t_b<r\}=\int \lambda_c\1\{r<t_c<t\}\1\{s<t_c<r\}=0.$$
With the help of these considerations we prove successively ($**$) and with the help of ($*$) and the 
sum integral lemma we obtain the result.
\end{pf}

\subsection{Unitarity}
\begin{thm}
Assume $u^t_s = u^t_s(A_1,A_0,A_{-1};B)$.
The mapping 
$$\mathcal{O}(u^t_s):\mathcal{K}\to \Gamma$$
can be extended to a unitary mapping
$$U^t_s: \Gamma \to \Gamma ,$$
iff the operators $A_i,i=1,0,-1;B$ fulfill the following conditions:
There exist a unitary operator $\Upsilon$ such that
\begin{align*}
A_0 &= \Upsilon-1\\
A_1&=-\Upsilon A_{-1}^+\\
B+B^+& =-A_1^+A_1 = - A_ {-1}A_{-1}^+.
\end{align*}
\end{thm}
\begin{pf}
We recall proposition 2.3..
For fixed $s$ is the function $t\mapsto u_s^t$ and for fixes $t$ is the function $s\mapsto u_s^t$
of class $\mathcal{C}^1$ and one has
\begin{eqnarray*}\partial^c_t u^t_s& =& Bu^t_s\\
(R^j_+u^._s)_t &=& A_j u^t_s \\
(R^j_-u^._s)_t &=& 0\\
\partial^c_s u^t_s& =& -u^t_sB\\
(R^j_+u^t_.)_s &=&  0\\
(R^j_-u^t_.)_s &=&  u^t_s A_j
\end{eqnarray*}
for $j=1,0,-1$.
We recall Ito's formula theorem 3.2.1
 Assume $x_t,y_t$ to be of class $\mathcal{C}^1$ and that for
  $f,g \in \mathcal{K}_s(\mathfrak{R},\mathfrak{k})$  the sesquilinear forms $\langle f|\mathcal{B}(F_t,G_t)|g \rangle$
  exist in norm and $t \in \mathbb{R}\mapsto \langle f|\mathcal{B}(F_t,G_t)|g \rangle$ is locally integrable,
  where 
  $$F_t \in \{x_t,\partial^cx_t,R^1_\pm x_t,R^0_\pm x_t ,R^{-1}_\pm x_t\}$$
 $$G_t \in \{y_t,\partial^cy_t,R^1_\pm y_t,R^0_\pm y_t ,R^{-1}_\pm y_t\}$$
 is any of the functions. 

Then the Schwartz derivative of
  $\langle f|\mathcal{B}(x_t,y_t)|g \rangle$ is a locally integrable function and yields
  \begin{align*}\partial \langle f|\mathcal{B}(x_t,y_t)|g \rangle = &
    \langle f|\mathcal{B}(\partial^c x_t,y_t) +\mathcal{B}(f,\partial^c y_t) + I_{-1,+1,t}|g \rangle \\
  &+   \langle a(t)f|\mathcal{B}(D^1 x_t,y_t) +\mathcal{B}(f,D^1y_t) + I_{0,+1,t}|g \rangle \\
   &+  \langle a(t)f|\mathcal{B}(D^0 x_t,y_t) +\mathcal{B}(f,D^0 y_t) + I_{0,0,t}|a(t)g \rangle \\
  &+  \langle f|\mathcal{B}(D^{-1} x_t,y_t) +\mathcal{B}(f,D^{-1} y_t) + I_{-1,0,t}|a(t)g \rangle \end{align*}
  with
 $$I_{i,j,t} = \mathcal{B}(R^i_+x_t,R^j_+y_t)-\mathcal{B}(R^i_-x_t,R^j_-y_t)$$

 We want to calculate the Schwartz derivativs of 
 $$ t \mapsto \langle f|\mathcal{B}( (u^t_s)^+,u^t_s)|g \rangle$$
 $$ s \mapsto \langle f|\mathcal{B}( u^t_s,(u^t_s)^+)|g \rangle.$$
 The derivativs exist, because $u:\mathcal{K }\to \Gamma$. We obtain
 \begin{multline*}\partial_t \langle f|\mathcal{B}( (u^t_s)^+,u^t_s)|g \rangle=
 \langle f |\mathcal{B}(u^+,C_1u)|g \rangle
 +\langle a(t)f |\mathcal{B}(u^+,C_2u)|g \rangle\\
+\langle a(t)f |\mathcal{B}(u^+,C_3u)|a(t)g \rangle+
\langle f |\mathcal{B}(u^+,C_4u)|a(t)g \rangle\end{multline*}
\begin{multline*}\partial_s \langle f|\mathcal{B}( u^t_s,(u^t_s)^+)|g \rangle=
 \langle f |\mathcal{B}(u,C_5u^+)|g \rangle
 +\langle a(t)f |\mathcal{B}(u,C_6u^+)|g \rangle\\
+\langle a(t)f |\mathcal{B}(u,C_7u^+)|a(t)g \rangle+
\langle f |\mathcal{B}(u,C_8u^+)|a(t)g \rangle\end{multline*}
with
\begin{align*}C_1 &= B+B^+ +A_1^+A_1&C_5 & =B + B^+ + A_{-1}A_{-1}^+\\
C_2 &= A_{-1}^+ + A_1+ A_0^+A_1& C_6 &= A_1 + A_{-1}^+ +A_0A_{-1}^+\\
C_3 & = A_0^+ +A_0 + A_0^+A_0& C_7 &= A_0+A_0^+ +A_0A_0^+\\
C_4 & = A_1^+ + A_{-1} + A_1^+ A_0 & C_8 &= A_{-1} + A_1^+ +A_{-1}A_0^+
\end{align*}
The operator $\mathcal{O}(u^t_s)$ is unitary iff both derivativs vanish and they vanish,iff $C_i=0,i=1,\cdots ,8.$
The equations $C_3 =0$ and $C_7=0$ imply
$$ (1+A_0^+)(1+A_0)= (1+A_0)(1+A_0^+)=1.$$
So 
$$ \Upsilon = 1+A_0$$
is unitary.
The equations are not independent.We have $C_2^+ = C_4$ and $C_6^+ =C_8$.Furthermore
$$C_2 = A_{-1}^+ +(1+A_0^+)A_1 =  A_{-1}^+ +\Upsilon^+ A_1 = \Upsilon^+ C_6.$$
So $C_2 =0$ implies $A_1 = - \Upsilon A_{-1}^+$, and from there $C_1 =C_5.$
\end{pf}
\begin{defi}
 For $ t<s $ we define
$$ U^t_s = (U^s_t)^+.$$
\end{defi}
\begin{prop}
For $r,s,t, \in \mathbb{R}$ we have
$$ U^t_rU^r_s = U^t_s.$$
\end{prop}
\begin{pf}
For $ s<r<t$ and $t<r<s$ the assertion follows from proposition 4.4.1. For $s<t<r$ we calculate
$$ \langle f |U^t_r U^r_s g \rangle = \langle U^r_tf|( U^r_t)^+ (U^t_s g)\rangle = \langle f | U^t_s g \rangle.$$
The other variants can be calculated in a similar way.\end{pf}
\subsection{Estimation of the $\Gamma_k$-norm}
Recall subsection 1.10
 $$ \int\mathfrak{m}(\pi,\sigma,\tau,\upsilon,\varrho)f^+(\omega) F(\sigma,\tau,\upsilon)g(\varrho)=
 \int f^+(\omega) ((\mathcal{O}F)g)(\omega) = \langle f, \mathcal{O}(F)g \rangle $$
 with
 $$((\mathcal{O}(F)g )(\omega)= \sum_{\alpha \subset \omega}\sum_{\beta \subset \omega \setminus\alpha}
 \int_\upsilon \lambda_\upsilon F(\alpha,\beta,\upsilon)g(\omega\setminus\alpha+ \upsilon).$$
 So $\mathcal{O}(F)$ is a mapping from $\mathcal{K}_s(\mathfrak{X})$ into the locally $\lambda$-integrable functions
 on $\mathfrak{X}$. Extend it to those functions $g$, such that the integral exists in norm for almost all $\omega$.
\begin{lem}
Assume a locally integrable function $F:\mathfrak{R}^3 \to B(\mathfrak{K})$
 such that for $g \in \mathcal{K}$ we have $\|\mathcal{O}(F)g\|_\Gamma
 \le \text{const}\|g\|_\Gamma$ and $\|\mathcal{O}(F^+)g\|_\Gamma\le \text{const}\|g\|_\Gamma$. Let $T:\Gamma \to
 \Gamma$ the operator , such that $\mathcal{O}(F)$ is the restriction of $T$ to $\mathcal{K}$. Then 
  $\mathcal{O}(F^+)$ is the restriction of $T^+$ to $\mathcal{K}$. Assume $f \in \Gamma$ such that
   $\mathcal{O}(\|F\|)\|f\|  \in L^2(\mathbb{R})$, then 
   $$Tf = \mathcal{O}(F)f$$
   \end{lem}
   \begin{pf}
   Assume $g,h \in \mathcal{K}$,then 
   \begin{multline*}\langle h |Tg\rangle=
    \langle h|\mathcal{O}(F)g\rangle = \int h^+(\omega) F(\sigma,\tau,\upsilon)g(\varrho)\mathfrak{m}\\= 
   \overline{\int g^+(\omega) F^+(\sigma,\tau,\upsilon)h(\varrho)\mathfrak{m}} = \langle \mathcal{O}(F^+)h|g \rangle 
   =\langle T^+h|g\rangle
   \end{multline*}
   with 
   $$ \mathfrak{m}= \langle a_\omega a^+_{\sigma+\tau}a_{\tau+\upsilon}a^+_\varrho\rangle\lambda_{\omega+\upsilon}.$$
So $\mathcal{O}(F^+)$ is the restriction of $T^+$ to $\mathcal{K}$.

We have
$$ \int \|h(\omega)\| \|F(\sigma,\tau,\upsilon)\| \|g(\varrho)\| \mathfrak{m}< \infty.$$
Hence
$$ \int h^+(\omega) (\mathcal{O}(F)f)(\omega)d\omega = \overline{\int f^+(\omega) (\mathcal{O}(F^+)h)(\omega)d\omega }
= \overline{\langle f| T^+h\rangle}= \langle h|Tf \rangle.$$
As this holds for any $h \in \mathcal{K}$ the assertion follows.
\end{pf}
 \begin{lem}
 Asssume $u^t_s = u^t_s(A_i,B)$ satisfying the unitarity conditions and $U^t_s$ the  corresponding unitary operator.
 Assume $G_1,\cdots ,G_k \in B(\mathfrak{k})$ and $ s=t_0 < t_1 <\cdots <t_k<t_{k+1}=t$ and
 \begin{multline*} F(\sigma,\tau,\upsilon) = \\\sum_{\substack{\sigma_0+\sigma_1+\cdots +\sigma_k = \sigma\\
 \tau_0+\tau_1+\cdots +\tau_k=\tau\\\upsilon_0+\upsilon_1+\cdots +\upsilon_k=\upsilon}}
 u^t_{t_k}(\sigma_k,\tau_k,\upsilon_k)G_k u^{t_k}_{t_{k-1}}(\sigma_{k-1},\tau_{k-1},\upsilon_{k-1})G_{k-1}\\\cdots 
 G_2 u^{t_2}_{t_1}(\sigma_1,\tau_1,\upsilon_1)G_1 u^{t_1}_s(\sigma_0,\tau_0,\upsilon_0)
 .\end{multline*}
 Then for $f \in \mathcal{K}$
 $$\mathcal{O}(F)g = U^t_{t_k}G_kU^{t_k}_{t_{k-1}}G_{k-1}\cdots G_2U^{t_2}_{t_1}G_1U^{t_1}_sg$$
 \end{lem}
 \begin{pf}
 The case $k=0$ is clear. We prove by induction from $k-1$ to $k$. Put for short
 $$ u(i) = u^{t_{i+1}}_{t_i}(\sigma_i,\tau_i,\upsilon_i).$$
 Then 
 $$ F = \sum u(k)G_k\cdots u(1)G_1u(0)= \sum _
 {\substack{\sigma_k+\sigma ' =\sigma\\\tau_k+\tau '=\tau\\\upsilon_k+\upsilon '=\upsilon}}
 u(k)G_k F'(\sigma ',\tau ',\upsilon ')$$
 with
 $$F'(\sigma ',\tau ',\upsilon ') = \sum_{\substack{\sigma_0+\cdots +\sigma_{k-1}=\sigma '\\
 \tau_0+\cdots+ \tau_{k-1} = \tau '\\\upsilon_0+\cdots +\upsilon_{k-1} = \upsilon '}}
 u(k-1)G_{k-1}\cdots u(1)G_1u(0)$$
 Go back to the  proof of 4.4.1. Put 
 $$ C = \max( \|A_i\|,i=1,0,-1;\|B\|;\|G_i\|,i=1,\cdots ,k).$$
 For $g\in \mathcal{K}$ we have
 $$\|(\mathcal{O}(\|F\|)\|g\|)\|_\Gamma \le c\|g\|_\Gamma$$
 $$\|(\mathcal{O}(\|F'\|)\|g\|)\|_\Gamma \le c'\|g\|_\Gamma.$$
 For $h\in \mathcal{K}$ 
 $$ \int h^+(\omega) (\mathcal{O}(F)g)(\omega)\lambda_\omega = \int h^+(\pi)F(\sigma,\tau,\upsilon)g(\varrho)
 \langle a_\pi a^+_{\sigma+\tau}a_{\tau+\upsilon}a^+_\varrho\rangle\lambda_{\pi+\upsilon}.$$
 Using the same argument as in the proof of 4.4.1 the last term equals
 $$ \int h^+(\pi)u(k) G_k F'(\sigma ',\tau ',\upsilon ')g(\varrho)
 \langle a_\pi a^+_{\sigma_k+\tau_k}a_{\tau_k+\upsilon_k}a^+_{\sigma ' + \tau '} a_{\tau ' +\upsilon '}a^+_\varrho 
 \rangle \lambda_{\pi+\upsilon_k+\upsilon '}.$$
 Now following theorem 1.9.1
 $$\langle a_\pi a^+_{\sigma_k+\tau_k}a_{\tau_k+\upsilon_k}a^+_{\sigma ' + \tau '} a_{\tau ' +\upsilon '} a^+_\varrho
 \rangle=
 \int_\omega \langle a_\pi a^+_{\sigma_k+\tau_k}a_{\tau_k+\upsilon_k}a^+_\omega\rangle
 \langle a_\omega a^+_{\sigma ' + \tau '} a_{\tau ' +\upsilon '} a^+_\varrho
 \rangle $$
 As 
 $$\int F'(\sigma ',\tau ',\upsilon ')g(\varrho)
 \langle a_\omega a^+_{\sigma ' + \tau '} a_{\tau ' +\upsilon '} a^+_\varrho\rangle \lambda_{\upsilon '}
 = (\mathcal{O}(F')g)(\omega) = f(\omega)$$
 the integrability conditions are fulfilled and one obtains
 $$ \int h^+(\omega) (\mathcal{O}(F)g)(\omega)\lambda_\omega
 =\int h^+(\pi) u(k) G_k  f(\omega)\langle a_\pi a^+_{\sigma_k+\tau_k}a_{\tau_k+\upsilon_k}a^+_\omega\rangle
 \lambda_{\pi+\upsilon_k}.$$
 The conditions of the preceding lemma are fulfilled, the last expression equals
 $$ \langle h| U^t_{t_k}f\rangle
  =\langle h| U^t_{t_k}G_kU^{t_k}_{t_{k-1}}G_{k-1}\cdots G_2U^{t_2}_{t_1}G_1U^{t_1}_sg\rangle$$
 using the hypothesis of induction.
 \end{pf}
 \begin{thm}For any $k$ there exists a polynomial $P$ of degree $ \le k$, such that for $g\in \Gamma_k$
 $$ \|U^t_s g\|_{\Gamma_k}^2 \le P(|t-s|)\|g\|_{\Gamma_k}^2.$$
 \end{thm}
 \begin{pf}
 Following lemma 4.3.1 and proposition 4.3.1 we have for $f \in \mathcal{K}$
 $$ \langle U^t_s f|\binom{N}{k}|U^t_s f\rangle = \int \|(U^t_sf)(\omega+\xi)\|^2\lambda_{\omega+\xi}<\infty$$
 Hence $\int \|(U^t_sf)(\omega+\xi)\|^2\lambda_{\omega}< \infty$ for allmost all $\xi$.We have
 $$ (\mathcal{O}(u^t_s))(\omega) = \sum_{\omega_1+\omega_2+\omega_3 = \omega}\int \lambda_\upsilon
 u^t_s(\omega_1,\omega_2,\upsilon)f(\omega_2+\omega_3 +\upsilon)$$
 and
 \begin{multline*} (\mathcal{O}(u^t_s))(\omega+\xi)\\ = \sum_{\substack{\omega_1+\omega_2+\omega_3 = \omega\\ 
 \xi_1+\xi_2+\xi_3=\xi}}
 \int \lambda_\upsilon
 u^t_s(\omega_1+\xi_1,\omega_2+\xi_2,\upsilon)f(\omega_2+\xi_2+\omega_3+\xi_3 +\upsilon)\\
 = (\mathcal{O}((u^t_s)_{\xi_1,\xi_2}))f_{\xi_2+\xi_3})(\omega)\end{multline*}
 with
 $$(u^t_s)_{\xi_1,\xi_2}(\sigma,\tau,\upsilon)= u^t_s(\sigma+\xi_1,\tau+\xi_2,\upsilon)$$
 $$  f_{\xi_1+\xi_2}(\varrho)=f(\xi_1+\xi_2+\varrho).$$
 Assume that the multiset $\{ s,t,t_\xi,t_{\sigma+\tau+\upsilon}\}^\bullet$ is without multiple points
 and order the set
 $$\{(t_i,1):i\in\xi_1\}+\{t_i,0):i\in \xi_0\} = \{(t_1,i_1),\cdots ,(t_l,i_l)\}$$
 with $t_1<\cdots < t_l$ and $i_j\in\{1,0\}$.Then
  \begin{multline*}
   (u^t_s)_{\xi_1,\xi_2}(\sigma,\tau,\upsilon) = 
  \\\sum_{\substack{\sigma_0+\sigma_1+\cdots +\sigma_l = \sigma\\
 \tau_0+\tau_1+\cdots +\tau_l=\tau\\\upsilon_0+\upsilon_1+\cdots +\upsilon_l=\upsilon}}
 u^t_{t_l}(\sigma_l,\tau_l,\upsilon_l)A_{i_l} u^{t_l}_{t_{l-1}}(\sigma_{l-1},\tau_{l-1},\upsilon_{l-1})A_{i_{l-1}}
 \\\cdots 
 A_{i_2} u^{t_2}_{t_1}(\sigma_1,\tau_1,\upsilon_1)A_{i_1} u^{t_1}_s(\sigma_0,\tau_0,\upsilon_0)
 .\end{multline*}
 Using the last lemma we obtain for $h\in\mathcal{K}$
 $$\mathcal{O}((u^t_s)_{\xi_1,\xi_2})h = U^t_{t_l}A_{i_l}\cdots A_{i_2}U^{t_2}_{t_1}A_{i_1}U^{t_1}_sh$$
 If $C= \max(\|A_i\|,\|B\|)$, then
 $$\|\mathcal{O}((u^t_s)_{\xi_1,\xi_2})h\|_\Gamma \le C^{\#(\xi_1+\xi_2)}\1\{t_{\xi_1+\xi_2}\subset [s,t]\}\|h\|_\Gamma.$$
 Finally
 \begin{multline*}  \langle U^t_s f|\binom{N}{k}|U^t_s f\rangle = \int \lambda_{\xi_1+\xi_2+\xi_3}\int \lambda_\omega
 \|\mathcal{O}(u^t_s)_{\xi_1,\xi_2}f_{\xi_2+\xi_3}h\|(\omega)^2\\
 \le \int\lambda_{\xi_1+\xi_2+\xi_3}C^{2\#(\xi_1+\xi_2)} 
 \1\{ t_{\xi_1+\xi_2}\subset s,t\}
 \|f_{\xi_2+\xi_3}\|^2_\Gamma\\
 \le C^{2k}\sum_{k_1+k_2+k_3=k}(t-s)^{k_1}/(k_1!)\langle f|\binom{N}{k_2+k_3}|f \rangle
 \end{multline*}
 The proof for $t<s$ goes in the same way.\end{pf}
 
 \section{The Hamiltonian}
 \subsection{Definition of the One Parameter Group $W(t)$}
 
Denote by $\vartheta(t)$ the rightshift on $\mathbb{R}$ and extend it to $\mathfrak{R}$,
$$\vartheta(t)(t_1,\cdots,t_n)=(t_1+t,\cdots,t_n+t)$$
If $\{t_1,\cdots,t_n\}^\bullet$ is a multiset we define 
$$\vartheta(t)\{t_1,\cdots,t_n\}^\bullet=\{t_1+t,\cdots,t_n+t\}^\bullet$$
In the notation $\{t_1,\cdots,t_n\}^\bullet = t_{\alpha}$ we write $\vartheta(t)t_{\alpha} = t_{\alpha}
+t\, e_{\alpha}$ with $e_{\alpha}=\{1,\cdots,1\}^\bullet$.\\

If $f$ is a function on $\mathfrak{R}$, then $(\vartheta(t)f)(w)= f(\vartheta(t)w)$. If
$\mu$ is a measure on $\mathfrak{R}$,then $\vartheta(t)\mu$ is defined by the property
$$\int (\vartheta(t)\mu(dw))\,(\vartheta(t)f(w)) = \int \mu(dw)f(w).$$
If $\mu(dw)=g(w)dw$ then $\vartheta(t)\mu(dw)= (\vartheta(t)g)(w)dw$.
Simlar notations hold for $\mathfrak{R}^k$.
We have
$$ (\vartheta(t)\varepsilon_x)(dy) = \varepsilon_{t-x}(dy).$$
\begin{lem}
If $W$ is an admissible sequence, then $\langle W\rangle \lambda_{\omega_-\setminus\omega_+}$ is a shift 
invariant measure.
\end{lem}
\begin{pf}
Following 1.9 this measure is a sum of measures where each one is a tensor product of measures of the form
$$ \Lambda(dt_1,\cdots ,dt_n): \int \Lambda(dt_1,\cdots ,dt_n) f(t_1,\cdots ,t_n) = \int dt f(t,\cdots ,t).$$
So it is clearly invariant due to the invariance of Lebesgue measure $dt$.
\end{pf}

Using the defining formulas one obtains immediately
\begin{lem} One has
$$ \Theta(r)u^t_s= u^{t-r}_{s-r}$$
for \  $r,s,t \in \mathbb{R}$.\end{lem}

Define a unitary operator $\Theta(t)$ on $\Gamma$ by $\Theta(t)f = \Theta(t)f$.
\begin{prop} The operators $U^t_s, \,s,t \in \mathbb{R}$ form a cocycle with respect
to $\Theta(t)$, i.e. 
$$\Theta(r)U^t_s\Theta(-r) = U^{t-r}_{s-r}.$$
\end{prop}
\begin{pf}
Use the invariance of 
$$ \mathfrak{m} = \langle a_\pi a^+_{\sigma+\tau}a_{\tau+\upsilon}a^+_\varrho\rangle\lambda_{\pi+\upsilon}.$$
and obtain
\begin{multline*} f^+(\pi)u^{t-r}_{s-r}(\sigma,\tau,\upsilon)g(\varrho)\mathfrak{m}= \int (\Theta(-r)f^+)(\pi)
u^t_s(\sigma,\tau,\upsilon)( \Theta(-r)g)(\varrho)\mathfrak{m}\\= \langle \Theta(-r)f|U^t_s \Theta(-r)g \rangle
= \langle f| \Theta(r)U^t_s\Theta(-r)g\rangle.\end{multline*}
\end{pf}
\begin{prop}Define for $t \in \mathbb{R}$

$$ W(t) = \Theta(t)\,U^t_0,$$
then $W(t)$ is a  unitary strongly continuous one parameter semigroup on $\Gamma$.
\end{prop}

\begin{pf}We have $W(0)=1$ and
$$W(s+t)= \Theta(t+s)U^{t+s}_0 = \Theta(s)\Theta(t)U^{t+s}_s\Theta(-t)\Theta(t)U^t_0 =
W(s)W(t)$$
and
$$W(t)^+ = U^0_t\,\Theta(-t)=\Theta(-t)U^{-t}_0=W(-t)$$
  \end{pf}
An immediate consequence of theorem 4.6.1 is
\begin{prop}The operators $W(t)$ map the space $\Gamma_k$ into itself and 
$$\|W(t)f\|^2_{\Gamma_k }\leq P(|t|) \|f\|^2_{\Gamma_K}$$
where $P$ is a polynomial of degree $ \leq k$.\end{prop}
\subsection{Definition of $\hat{\mathfrak{a}},\hat{\mathfrak{a}}^+,\hat{\partial}$}

If $\varphi$ is an integrable function on the real line, we define
$$\Theta(\varphi) = \int \varphi(t)\Theta(t)dt$$
an operator mapping $\Gamma_k$ into $\Gamma_k$ for all $k$ .

If $\nu$ is a measure on $\mathbb{R}$ and $f$ a locally integrable function , symmetric function on $\mathfrak{R}$,
then 
$$ (a^+(\nu)f\lambda)(\omega) = \sum _{c\in \omega} \nu(c)f(\omega\setminus c)\lambda(\omega\setminus c)$$
We shall use again the convention of L.Schwartz  and denote $f\lambda$ by$f$.
So $$ (a^+(\nu)f)(\omega) = \sum _{c\in \omega} \nu(c)f(\omega\setminus c).$$
We set
\begin{align*} \mathfrak{a}&=a(\varepsilon_0)=a(0)& \mathfrak{a}^+&= a^+(\varepsilon_0) \end{align*}
We use gothic $\mathfrak{a}^+$ in order to distinguish from $a^+(dx)=a^+(\varepsilon(dx))$ used in the preceding text.
 $a(\varepsilon_0)=a(0)g$ is the same as before.We have 
 $$ (\mathfrak{a}^+f)(\omega) = \sum_{c \in \omega}\varepsilon_0(dt_c)f(t_{\omega\setminus c})
 =\sum_{c \in \omega}\varepsilon_0(dt_c)f(t_{\omega\setminus c})\lambda(dt_{\omega\setminus c}).$$
 The duality relation (subsection 1.3)
 becomes
 $$ \int g(\omega)(\mathfrak{a}^+f)(\omega)= \int (\mathfrak{a}g)(\omega)f(\omega)\lambda(\omega).$$
\begin{lem}
Assume $f \in L^2(\mathbb{R}^n)$ and $\varphi \in (L^1\cap L^2)(\mathbb{R})$. Then $\Theta(\varphi)$
maps the singular measure $\mathfrak{a}^+ f $ into a absolute continuous measure identified with its density and
we have
 $$ (\Theta(\varphi)\mathfrak{a}^+f)(t_\omega)= 
  \sum_{c\in \omega} \varphi(-t_c) (\Theta(-t_c)f)(\omega\setminus c)$$
The application $\Theta(\varphi)\mathfrak{a}^+$ can be extended to a mapping $\Gamma_k \to \Gamma_{k-1}$ and we have
 $$ \|\Theta(\varphi)\mathfrak{a}^+f\|_{\Gamma_{k- 1}} \le \|\varphi\|_{L^2}\|f\|_{\Gamma_k}.$$
 We have
 $$ \int g(\omega)^+(\Theta(\varphi)\mathfrak{a}^+f)(t_\omega)d\omega =
 \int (\mathfrak{a}\Theta(\varphi^+) g)^+(\omega)f(\omega)d\omega.$$
 with $\varphi^+(t)=\overline {\varphi(-t})$. One obtains
 $$(\mathfrak{a}\Theta(\varphi)f )(t_1,\cdots ,t_n)= \int \varphi(s) ds f(s,t_1+s,\cdots ,t_n+s).$$
The application can be extended to a mapping $\Gamma_k \to \Gamma_{k-1}$ and we have
 $$ \|\mathfrak{a}\Theta(\varphi)f\|_{\Gamma_{k-1}} \le \|\varphi\|_{L^2}\|f\|_{\Gamma_k}.$$
\end{lem}
\begin{pf}One has
 \begin{multline*}\int ds \,\varphi(s)\Theta(s)
\sum_{c \in \omega}\varepsilon_0(dt_c)f(t_{\omega\setminus c})\lambda(dt_{\omega\setminus c})\\
=
\int ds \,\varphi(s)\Theta(s)
\sum_{c \in \omega}\varepsilon_{-s}(dt_c)(\Theta(s)f)(t_{\omega\setminus c})\lambda(dt_{\omega\setminus c})\\=
\int ds\, \varphi(-s)\Theta(-s)
\sum_{c \in \omega}\varepsilon_{s}(dt_c)(\Theta (-s)f)(t_{\omega\setminus c})\lambda(dt_{\omega\setminus c})\\=
\sum_{c\in \omega} \varphi(-t_c) (\Theta(-t_c)f)(\omega\setminus c)\lambda(\omega) \end{multline*}
So $\Theta(\varphi)$ works as a mollyfier,as it is called in
the theory of Schwartz's distributions, making a function out of a singular measure .
The other results folllow by simple calculations.
\end{pf}

\begin{lem} If $\varphi \in (L^1 \cap L^2)(\mathbb{R})$ and $f \in L^2(\mathbb{R}^{n+1})$,then
$$ \begin{array}{l}x\in\mathbb{R}^{n+1} \mapsto g(x)\\
g(x)(t_1,\cdots ,t_n)= (\Theta(\varphi)f)((0,t_1,\cdots ,t_n)+x)\end{array}$$
maps $\mathbb{R}^{n+1}$ into $L^2(\mathbb{R}^n)$ and $x\mapsto g(x)$ is a continuous function bounded by
$\|\varphi\|_{L^2(\mathbb{R})}\|f\|_{L^2(\mathbb{R}^{n+1})}$.\end{lem}
\begin{pf}We have with $e=(1,1,\cdots ,1)$
$$\begin{array}{l}\|g(x)-g(y)\|^2_{L^2(\mathbb{R}^n)}\\ = \int dt_1\cdots dt_n \|\int ds\,\varphi(s)(f((0,t_1, \cdots ,t_n)+se +x)
-f((0,t_1, \cdots ,t_n)+se +y)\|^2\\
\leq\|\varphi\|_{L^2(\mathbb{R})}^2\int dt_1\cdots dt_n \int ds \\
\quad\|f(s+x_0,t_1+s+x_1,\cdots ,t_n+s+x_n)
-f(s+y_0,t_1+s+y_1,\cdots ,t_n+s+y_n)\|^2\\
=\|\varphi\|_{L^2(\mathbb{R})}^2 \int dt_0\cdots dt_n\|f(t_0+x_0,\cdots ,t_n+x_n)
-f(t_0+y_0,\cdots ,t_n+y_n)\|^2\\
=\|\varphi\|_{L^2(\mathbb{R})}^2 \|(T(x)-T(y))f\|^2_{L^2(\mathbb{R}^{n+1})}\end{array}$$
where $T(x)$ denotes the translation by $x$. The bound for $\|g(x)\|$ can be shown in the same way.
\end{pf}

\begin{lem} 
Assume $f \in L^2(\mathbb{R}^n)$ and $\eta \in L^2(\mathbb{R})$ a continuous bounded function on 
$\mathbb{R}\setminus\{0\}$ with right and left limits at $0$ or with other words $\eta$ is a continuous bounded
function on $\mathbb{R}_0$. We show 
\begin{align*}x \in \mathbb{R}^{n+1}&\mapsto g(x) \in L^2(\mathbb{R}^n)\\
 g(x)(t_1,\cdots ,t_n) &= (\Theta(\eta)\mathfrak{a}^+f)((0,t_1,\cdots ,t_n)+x)\end{align*}
 and 
 $$\|g(x)\|_{L^2(\mathbb{R}^n)}\le (n+1)\|\eta\|_\infty \|f\|_{L^2(\mathbb{R}^n)}$$
 for all $x \in L^2(\mathbb{R}^{n+1})$ and $x \mapsto g(x)$ is continuous in $\{(x_0,x_1,\cdots ,x_n): x_0 \ne 0 \}$.
 We have
 $g(0\pm,x_1,\cdots ,x_n)$ exist and
 $$ g(0+,x_1,\cdots ,x_n)-g(0-,x_1,\cdots ,x_n) = -(\eta(0+)-\eta(0-))T(x_1,\cdots ,x_n)f$$
 where $T(x)$ is the translation by $x$.
\end{lem}
\begin{pf}
We have
$$(\Theta(\eta)\mathfrak{a}^+f)(t_0,t_1,\cdots,t_n) = k_0(t_0,t_1,\cdots,t_n)+\cdots+k_n(t_0,t_1,\cdots,t_n)$$
with 
\begin{align*}k_0(t_0,t_1,\cdots,t_n)& =\eta(-t_0)
f(t_1-t_0,t_2-t_0,\cdots,t_n-t_0)\\
k_1(t_0,t_1,\cdots,t_n)&=\eta(-t_1)f(t_0-t_1,t_2-t_1,\cdots,t_n-t_1)\\
\vdots\\ 
k_n(t_0,t_1,\cdots,t_n)&=\eta(-t_n)f(t_0-t_n,t_1-t_n,\cdots,t_{n-1}-t_n)\end{align*}
Define
$$g_i(x)(t_1,\cdots ,t_n)=k_i((0,t_1,\cdots ,t_n)+x)$$
and discuss at first $g_i$ with $i \ne 0$, e.g. $g_n$.We have
\begin{align*}
&g_n(x)(t_1,\cdots ,t_n)
=k_n((0,t_1,\cdots ,t_n)+x)\\
&=\eta(-x_n-t_n)f(x_0-x_n-t_n,x_1-x_n+t_1-t_n,\cdots ,x_{n-1}-x_n+t_{n-1}-t_n)
\\&=\eta(-x_n-t_n)(T(x')f)(-t_n,t_1-t_n,\cdots ,t_{n-1}-t_{n-1})\end{align*}
with
$$x'=(x_0-x_n,x_1-x_n,\cdots ,x_{n-1}-x_n).$$
From there one obtains
$$\|g_n(x)\|\leq\|\eta\|_\infty\|f\|.$$
We have
\begin{multline*}\int dt_1\cdots dt_n\|k_n((0,t_1,\cdots ,t_n)+x)-k_n((0,t_1,\cdots ,t_n)+y)\|^2\\
\leq2\int dt_1\cdots dt_n
|\eta(-x_n-t_n)-\eta(-y_n-t_n)|^2 \\\|(T(x')f)(-t_n,t_1-t_n,\cdots ,t_{n-1}-t_{n-1})\|^2\\
+2\|\eta\|_\infty^2\int dt_1\cdots dt_n \|(T(x')-T(y'))f(-t_n,t_1-t_n,\cdots ,t_{n-1}-t_{n-1})\|^2\end{multline*}
For $y\rightarrow x$ the first term goes to zero by the theorem of Lebesgue, in the second term
observe, that $z\mapsto T(z)f$ is normcontinuous.So $g_n(x)$  and $g_i(x),i\ne 0$ is continuous for all $x$.
We have
$$g_0(x)= \eta(-t_0)f(t_1+x_1-x_0,t_2+x_2-x_0,\cdots ,t_n+x_n-x_0).$$
From there one obtains the result.
\end{pf}

We double the point $0$ and introduce
$$\mathbb{R}_0 = ]-\infty,-0]+[+0,\infty[$$
with the usual topology. A function $f$ on $\mathbb{R}_0$ is continuous, if its restriction to $\mathbb{R}\setminus
\{0\}$ is continuous and if both limits $f(\pm 0)$ exist.We define
$$\mathfrak{R}_0 = \{\emptyset\} + \mathbb{R}_0 + \mathbb{R}_0^2 +\cdots .$$
We introduce on $\mathbb{R}_0$ and on $\mathfrak{R}_0$ the Lebesgue measure $\lambda$.
A continuous function on $\mathbb{R}\setminus\{0\}$ which has left and right limits at $0$, can be considered 
as a function on $\mathbb{R}_0$.
 We define the  measures
$ \varepsilon_{\pm 0}$ and for symmetric functions $f$ on $\mathfrak{R}_0$ the operators
$\mathfrak{a}_{\pm}= a(\varepsilon_{\pm 0})$ and $\mathfrak{a}^+_{\pm}= a^+(\varepsilon_{\pm 0})$ and 
shall use similar conventions as above.We put
\begin{align*} \hat \varepsilon_0 =& \tfrac {1}{2}(\varepsilon_{+0}+\varepsilon_{-0})&
\hat{\mathfrak{a}} = &\tfrac {1}{2}(\mathfrak{a}_+ +\mathfrak{a}_-)&
\hat{\mathfrak{a}}^+ = &\tfrac {1}{2}(\mathfrak{a}^+_+ +\mathfrak{a}^+_-).
\end{align*}
 A 
$\delta$-sequence is a sequence of functions $\varphi_n \in C_c^{\infty}$ such that 
\begin{align*} \int \varphi_n(t)dt &=1,&
 \int |\varphi_n(t)|dt \leq C &< \infty, &
 \mathrm{supp}(\varphi_n) \subset ]-\varepsilon_n,\varepsilon_n[\end{align*}
and $\varepsilon_n \downarrow 0$.
\begin{defi}We call a $\delta$-sequence $\varphi_n$ a 
\emph{symmetric}  $\delta$-sequence if the $\varphi_n$ are real and if 
$\varphi_n(t) = \varphi_n(-t)$ for all $n$ and $t$. \end{defi}

\begin{prop} Assume two functions $f$ and $\eta$ fulfilling the conditions of lemma 5.2.3, then
$\hat{\mathfrak{a}}\Theta(\eta)\mathfrak{a}^+f$ exists ,and
$$\|\hat{\mathfrak{a}}\Theta(\eta)\mathfrak{a}^+f\|_\Gamma \le \|\eta\|_\infty \|f\|_{\Gamma_2}$$
and
if $\varphi_n$ is a symmetric $\delta$-sequence, then $\mathfrak{a}\Theta(\varphi_n)\Theta(\eta)\mathfrak{a}f$
exists and
 $$\mathfrak{a}\Theta(\varphi_n)\Theta(\eta)\mathfrak{a}^+f \to
  \hat{\mathfrak{a}}\Theta(\eta)\mathfrak{a}^+f$$
\end{prop}
\begin{pf}
We apply  lemma 5.2.3. We have
$$(\Theta(\eta)\mathfrak{a}^+f)(s,t_1,\cdots ,t_n) = g(s,0,\cdots ,0)(t_1,\cdots ,t_n)$$
and $g(s,0,\cdots ,0) \to g(0+,0,\cdots,0)$ 
for $s \downarrow 0$. From there one concludes the existence of 
 $\hat{\mathfrak{a}}\Theta(\eta)\mathfrak{a}^+f$.We have
 $$(\mathfrak{a}\Theta(s)\Theta(\eta)\mathfrak{a}^+f )(t_1,\cdots ,t_n)= g(s,\cdots ,s)(t_1,\cdots ,t_n)$$
 and $g(s,\cdots ,s)\to g(0+,0,\cdots ,0)$ for $s\downarrow 0$. 
 From there one obtains the rest of the proposition.
\end{pf}

 Denote 
$$\zeta(t) = \mathbf{1}\{t>0\}e^{-t}$$
and 
$$Z=\Theta(\zeta)=\int_0^\infty e^{-t}\Theta(t)\,dt$$
\begin{defi} The vector space $D \subset \Gamma $ is defined by
$$D=\{f=Z(f_0+\mathfrak{a}^+
f_1):f_0\in\Gamma_1,f_1\in\Gamma_2\}$$
 \end{defi}
As a consequence of lemmata 5.2.2 and 5.2.3
for $f \in D $ the mapping $f\to\hat{\mathfrak{a}}f $ exists and defines a application
$\hat{\mathfrak{a}}:D\to \Gamma$, which stay bounded, if $\|f_0\|_{\Gamma_1}$ and $\|f_1\|_{\Gamma_2}$ stay
bounded.
Assume $\varphi_n$ to be a symmetric $\delta$-sequence and $f \in D$,then
$ \mathfrak{a}\Theta(\varphi_n)f \to \hat {\mathfrak{a}}f$
in the norm of $\Gamma$.

If $f \in \mathcal{C}^1(\mathbb{R}^n)$ and $\varphi_n$ a $\delta$-sequence. Then 
\begin{multline*}\Theta(\varphi_n')f(t)=\int \varphi_n'(s)f(t+se)ds\\=-\int \varphi_n(s)f'(t+se)ds
\rightarrow - \sum \frac{\partial f}{\partial t_i}(t) = - (\partial f)(t)\end{multline*}
This motivates
\begin{defi} We define
$$ \hat{\partial}= -\lim \Theta(\varphi_n'),$$
where $ \varphi_n$is a \emph{symmetric} $\delta$-equence.\end{defi}
 The sense of the convergence and the domain of 
$\hat{\partial}$ will be explained in the following proposition.
\begin{prop}Assume $\varphi_n$ to be a symmetric $\delta$-sequence, then for $n \rightarrow
\infty$ and  $f=Z(f_0+\hat{\mathfrak{a}}^+f_1)\in D$
$$-\Theta(\varphi_n')f \rightarrow \hat{\partial}f = -(f_0+\hat{\mathfrak{a}}
^+f_1)+f$$
weakly over $D$,i.e.for $g \in D$
$$ -\int g^+(\omega)(\Theta(\varphi_n')f)(\omega) \rightarrow \int g^+(\omega)(\hat{\partial}f)(\omega).$$
\end{prop}
\begin{pf}One has 
$$-\Theta(\varphi_n')Z=-\Theta(\varphi_n'\star\zeta)$$
where $\star$ denotes the
 usual convolution and
$$\varphi_n '\star \zeta = \varphi_n \star \zeta '= \varphi_n \star (\varepsilon_0 - \zeta) =
\varphi_n  -\varphi_n  \star \zeta .$$
So
$$ -\Theta(\varphi_n')Z = -\Theta(\varphi_n) +\Theta(\varphi_n)Z$$
and for $f=Z(f_0+\hat{\mathfrak{a}}^+f_1)$ we have 
$$-\Theta(\varphi_n')f = -\Theta(\varphi_n)(f_0+{\mathfrak{a}}^+f_1)+ \Theta(\varphi_n)f$$
Now $\Theta(\varphi_n)f$ and $\Theta(\varphi_n)f_0$ converge to $f$,resp.$f_0$ in $\Gamma$-
norm. We investigate for $g \in D$
\begin{multline*}\int g(t_0,t_1,\cdots,t_m)^+\Theta(\varphi_n)(\mathfrak{a}^+f_1)(t_0,t_1,\cdots,t_m)dt_0\cdots dt_m\\ =
\int (\mathfrak{a}\Theta(\varphi_n)g)(t_1,\cdots ,t_n)f(t_1,\cdots ,t_n)dt_1\cdots dt_n \\\to
\int (\hat{\mathfrak{a}}g)(t_1,\cdots ,t_n)f(t_1,\cdots ,t_n)dt_1\cdots dt_n\\
= \int g(t_0,\cdots ,t_n)(\hat{\mathfrak{a}}^+f)(t_1,\cdots t_n)\end{multline*}
by proposition 5.2.3.
\end{pf}
\begin{prop} The sesquilinear form
$$ f,g \in D\mapsto \langle f|i\hat{\partial}g \rangle = \int f^+(\omega) (i\hat{\partial}g)(\omega)$$
is symmetric,i.e.$\langle f|i\hat{\partial}g \rangle = \langle i\hat{\partial}f|g \rangle$.\end{prop}
\begin{pf} One has
\begin{multline*}\int g^+(t_0,t_1,\cdots,t_m)(\Theta(\varphi_n')f)(t_0,t_1,\cdots,t_m)dt_1\cdots dt_m\\=
-\int(\Theta(\varphi_n')g^+)(t_0,t_1,\cdots,t_m)
f(t_0,t_1,\cdots,t_m)dt_1\cdots dt_m\end{multline*}
as $\varphi_n(-s)= \varphi_n(s)$ and hence $\varphi_n(-s)'=-\varphi_n(s)$. This establishes the 
symmetry of $i\hat{\partial}$.\end{pf}

\subsection{Discussion of the Resolvent}
 We go  back to the group $W(t)$.
 \begin{defi}For $ z \in \mathbb{C},\Im z \neq 0$ we define the resolvent $R(z)$
 $$R(z) = \begin{cases}-i\int_0^\infty e^{izt}W(t)dt&\mbox{ for }\Im z > 0\\ 
 i\int_{-\infty}^0 e^{izt}W(t)dt&\mbox{ for }\Im z < 0 \end{cases}$$
 Furthermore we set 
 $$ S(z) = \begin{cases} -i\int_0^\infty e^{izt}W(t)a(t) & \text{ for } \Im z>0\\
 i\int^0_{-\infty }e^{izt}W(t)a(t) & \text{ for } \Im z<0 \end{cases}$$
 and
 \begin{equation*}\kappa(z)= \begin{cases}-i\mathbf{1}\{t>0\}e^{izt+Bt}&\text{ for }\Im z >0
 \\ i\mathbf{1}\{t>0\}e^{izt-B^+t} &\text { for } \Im z < 0 \end{cases}\end{equation*}
 and
  $$\tilde{R}(z) = \Theta(\kappa(z))= \left\{\begin{array}{rcl}-i\int_0^\infty e^{izt}e^{Bt}\Theta(t)
  dt&\mbox{ for }&\Im z > 0\\ \\
 +i\int_{-\infty}^0 e^{izt}e^{-B^+t}\Theta(t)dt&\mbox{ for }&\Im z < 0 \end{array}\right.$$
 \end{defi}
\begin{prop}
We have for $f \in \mathcal{K}$ 
 \begin{multline*}R(z)f=\tilde{R}(z)f\\ + \begin{cases}i\tilde{R}(z)\mathfrak{a}^+A_1R(z)f+
  i\tilde{R}(z)\mathfrak{a^+}A_0S(z)f+
 i\tilde{R}(z)A_{-1}S(z)f\\
 - i\tilde{R}(z)\mathfrak{a}^+A_{-1}^+R(z)f- i\tilde{R}(z)\mathfrak{a}^+A_0^+S(z)f-
 i\tilde{R}(z)A_1^+S(z)f\end{cases}\end{multline*}
 The upper line holds for $\Im z >0$, the lower one for $\Im z <0$.
\end{prop}
\begin{pf}
Directly from the definition for $t>s$
\begin{multline*}u^t_s(\sigma,\tau,\upsilon) = e^{B(t-s)}\mathbf{e}(\sigma,\tau,\upsilon)\\
+ \sum_{c\in \sigma}e^{B(t-t_c)}\1\{t_c \in [s,t]\}A_1u_s^{t_c}(\sigma\setminus c,\tau,\upsilon)\\+
\sum_{c\in \tau}e^{B(t-t_c)}\1\{t_c \in [s,t]\}A_0u_s^{t_c}(\sigma,\tau\setminus c,\upsilon)\\+
\sum_{c\in \upsilon}e^{B(t-t_c)}\1\{t_c \in [s,t]\}A_{-1}u_s^{t_c}(\sigma,\tau,\upsilon\setminus c)\end{multline*}
\begin{multline*} = e^{B(t-s)}\mathbf{e}(\sigma,\tau,\upsilon)\\
+\sum_{c\in \sigma}u^t_{t_c}(\sigma\setminus c,\tau,\upsilon)A_1e^{B(t_c-s)}\1\{t_c \in [s,t]\}\\
+\sum_{c\in \tau}u^t_{t_c}(\sigma,\tau\setminus c,\upsilon)A_0e^{B(t_c-s)}\1\{t_c \in [s,t]\}\\
+\sum_{c\in \upsilon}u^t_{t_c}(\sigma,\tau,\upsilon\setminus c)A_{-1}e^{B(t_c-s)}\1\{t_c \in [s,t]\}.
\end{multline*}
Hence
\begin{multline*}(u^t_s)^+(\sigma,\tau,\upsilon) =( u^t_s(\upsilon,\tau,\sigma)^+=
e^{B^+(t-s)}\mathbf{e}(\sigma,\tau,\upsilon)\\
+\sum_{c\in \sigma}e^{B^+(t_c-s)}\1\{t_c \in [s,t]\}A_{-1}^+(u_s^{t_c})^+(\sigma\setminus c,\tau,\upsilon)\\
+\sum_{c\in \tau}e^{B^+(t_c-s)}\1\{t_c \in [s,t]\}A_0^+(u_s^{t_c})^+(\sigma,\tau\setminus c,\upsilon)\\
+\sum_{c\in \upsilon}e^{B^+(t_c-s)}\1\{t_c \in [s,t]\}A_1^+(u_s^{t_c})^+(\sigma,\tau,\upsilon\setminus c)
\end{multline*}
Assume $f,g\in\mathcal{K}$ and using the same arguments as in the proof of theorem 2.2.1we obtain
\begin{multline*} \langle f| U^t_s g \rangle= \langle f |e^{B(t-s)}g \rangle+ \int _s^t dr 
(\langle a(r)f|e^{B(t-r)}A_1 U^r_sg \rangle \\+ \langle a(t)f|e^{B(t-r)}A_0U^r_sa(r)g\rangle +
 \langle f|e^{B(t-r)}A_{-1}U^r_s a(r)g \rangle)
\end{multline*}
and
\begin{multline*} \langle f| (U^t_s)^+ g \rangle= \langle f |e^{B^+(t-s)}g \rangle+ \int _s^t dr 
(\langle a(r)f|e^{B^+(r-s)}A_{-1}^+ (U^t_r)^+g \rangle \\+ \langle a(r)f|e^{B^+(r-s)}A_0^+(U^t_r)^+a(r)g\rangle +
 \langle f|e^{B^+(r-s)}A_1^+(U^t_r)^+a(r) g \rangle)
\end{multline*}
Finally for $t>0$
\begin{multline*} \langle f| U^t_0 g \rangle= \langle f |e^{Bt}g \rangle+ \int _0^t dr 
(\langle a(r)f|e^{B(t-r)}A_1 U^r_0g \rangle \\+ \langle a(r)f|e^{B(t-r)}A_0U^r_0a(r)g\rangle +
 \langle f|e^{B(t-r)}A_{-1}U^r_0 a(r)g \rangle)
\end{multline*}
and for $t<0$
\begin{multline*} \langle f|U^t_0 g\rangle =
\langle f| (U^0_t)^+ g \rangle= \langle f |e^{-B^+t}g \rangle+ \int _t^0 dr 
(\langle a(r)f|e^{B^+(r-t)}A_{-1}^+ U^r_0g \rangle \\+ \langle a(r)f|e^{B^+(r-t)}A_0^+U^r_0a(r)g\rangle +
 \langle f|e^{B^+(r-t)}A_1^+U^r_0a(r) g \rangle)
\end{multline*}
We want now to calculate the resolvent for $\Im z > 0$
$$ \langle f|R(z)g \rangle = -i \int _0^\infty dt e^{izt}\langle f|\Theta(t)U^t_0 g \rangle = 
-i \int _0^\infty dt e^{izt}\langle \Theta(-t)f|U^t_0 g \rangle$$
and consider e.g. the term
\begin{multline*}
-i\int_0^\infty dt e^{izt}\int _0^t dr\langle a(r) \Theta(-t)f |e^{B(t-r)}A_0 U^r_0 a(r)g\rangle\\ 
= -i \iint_0^\infty dt\; ds\;e^{iz(t+s)}\langle a(s)\Theta(-t-s)f|e^{Bt}A_0U^s_0 a(s)g \rangle \\
 = -i \iint_0^\infty dt\;ds\; e^{iz(t+s)}\langle a(0)e^{B^+t}\Theta(-t)f|A_0 \Theta(s)U_0^s a(s)g \rangle\\
= \langle \mathfrak{a}\tilde{R}(z)^+g|iA_0 S(z)g\rangle = i\langle f| \tilde{R}(z)\mathfrak{a}^+A_0 S(z)g \rangle. 
\end{multline*}
By similar calculations one finishes the proof.
\end{pf}
\begin{cor}
If $ f \in \mathcal{K}$ we may write
$$R(z) = \tilde{R}(z)(g_0(z)+\mathfrak{a}^+g_1(z))$$
with
$$g_0(z)= f + \begin{cases}iA_{-1}S(z)f& \text{ for } \Im z > 0\\
-iA_1^+ S(z)f& \text{ for } \Im z <0 \end{cases}$$
and 
$$g_1(z)= \begin{cases}iA_1R(z)f+iA_0S(z)f & \text{ for }\Im z > 0\\
-iA_{-1}^+R(z)f -iA_0^+ S(z)f&\text{ for }\Im z < 0\end{cases}.$$
\end{cor}
\begin{lem}
We have
$$ \tilde{R}(z)+ i Z = -i(1+iz + C(z))\tilde{R}(z)Z = -i(1+iz + C(z))Z\tilde{R}(z) $$
with
$$ C(z) = \begin{cases} B &\text{ for }\Im z > 0\\
-B^+& \text{ for } \Im z < 0 \end{cases}$$
\end{lem}
\begin{pf}
Use the equation
$$ \kappa(z)+i\zeta = -i(1+iz+C(z))(\zeta\star\kappa )=-i(1+iz+C(z))(\kappa\star\zeta)$$
\end{pf}
\begin{lem}
\label{VII5} We have that $f=Z(f_0+\mathfrak{a}^+f_1)$ is in $D$ iff $$f=\tilde{R}(z)(g_0+\mathfrak{a}^+g_1),$$
 with $g_0 \in \Gamma_1$ and $g_1 \in \Gamma_2$. 
\end{lem}
\begin{pf}
Assume at first that $f=\tilde{R}(z)(g_0+a^+*g_1)$,
using lemma 5.2.1 we see that $f \in \Gamma_1$ and by lemma 5.3.1
$$ f= -iZ(g_0+\mathfrak{a}^+g_1)-i(1+iz+C(z))Zf$$
and $f\in D$.By interchanging the roles of $Z$ and $\tilde{R}(z)$ one shows that every  $f \in D$ is of that form.
\end{pf}
\begin{prop}
The resolvent 
$$R(z):\mathcal{K}\to D$$.
\end{prop}
\begin{pf}
 By proposition 5.1.3
 $$\|W(t)f\|_{\Gamma_k}^2 \leq P(|t|)\|f\|^2_{\Gamma_k},$$
 where $P(|t|)$ is a polynomial in of degree  $\leq k$.So e.g. for $\Im z > 0$
 $$\|R(z)f\|_{\Gamma_k}\leq \int_0^\infty dt\,\exp{(-t\Im z)}\sqrt{P(|t|)}\|f\|_{\Gamma_k}.$$
 If $f\in\mathcal{K}$, the function $f\in\Gamma_k$ for all $k$. The functions $a(t)f,t\in\mathbb{R}$ are
 uniformely bounded in any $\Gamma_k$-norm. Hence
 $R(z)f$ and $S(z)f$
 are in $\Gamma_k$ for all $k$. The proposition follows from cor.5.3.1.
\end{pf}
\begin{prop}
Using the functions $g_0(z),g_1(z)$ of corollary 5.3.1 and $C(z)$ of lemma 5.3.1 we have for $f\in \mathcal{K}$
$$ i\hat{\partial}R(z)f = -g_0 -\hat{\mathfrak{a}}^+g_1 + (z -iC(z))R(z)f.$$
\end{prop}
\begin{pf}
The Schwartz derivative
$$ \kappa(z)' = -i\varepsilon_0 + (iz + C(z))\kappa(z)$$
and if $\varphi_n$ is a symmetric $\delta$-sequence, then

\begin{multline*}\Theta(\varphi_n')\tilde{R}(z)= \Theta(\varphi_n'\star\kappa(z))
=\Theta(\varphi_n\star\kappa(z)')\\=
-i\Theta(\varphi_n)+(iz+C(z))\Theta(\varphi_n)\tilde{R}(z)\end{multline*}
 So 
 $$-i\Theta(\varphi_n')\tilde{R}(z)(g_0+\mathfrak{a}^+g_1)= -\Theta(\varphi_n)(g_0+\mathfrak{a}^+g_1)+
 (z-iC(z))\Theta(\varphi_n)f$$
 converging weakly over $D$ to
 $$i\hat{\partial}f=-g_0-\hat{\mathfrak{a}}g_1+(z-iC(z))f$$
 From there one obtains immediately the result.
\end{pf}
\begin{prop}

 For $f \in \mathcal{K}$we have
$$ \hat{\mathfrak{a}}R(z)f = S(z)f + \begin{cases}\tfrac{1}{2}A_1R(z)f+\tfrac{1}{2} A_0S(z)f$$&\text{ for }\Im z > 0\\
\tfrac{1}{2}A_{-1}^+R(z)f+\tfrac{1}{2} A_0^+S(z)f &\text{ for }\Im z < 0.\end{cases}$$
\end{prop}
\begin{pf}
We prove at first the case $\Im z >0$. Following  the arguments of the proof of  
theotem 4.6.1 we have for $t>0$
\begin{multline*}
(U_{-t}^0f)(\omega + c) = ( U_{-t}^0 a(t_c)f)(\omega)\\+\1\{ t_c \in [-t,0]\}(
(U_{t_c}^0A_1U_{-t}^{t_c}f)(\omega)+(U_{t_c}^0A_0U_{-t}^{t_c}a(t_c)f)(\omega))\end{multline*}
and
\begin{multline*}(R(z)f)(\omega + c)\\ = -i\int_0^\infty dt e^{izt}(\Theta(t)U_0^tf)(\omega+ c)
 = -i\int_0^\infty dt e^{izt}(U_{-t}^0\Theta(t)f)(\omega + c)\\
 = -i\int_0^\infty dt(U_{-t}^0\Theta(t)a(t_c+t)f)(\omega)
-i \1\{t_c< 0\}U_{t_c}^0\\\left(A_1 \int_{-t_c}^\infty dt e^{izt}(U_{-t}^{t_c}\Theta(t)f)(\omega)+
A_0 \int_{-t_c}^\infty dt e^{izt}(U_{-t}^{t_c}\Theta(t)a(t_c+t)f)(\omega)\right).\end{multline*}
One concludes
\begin{align*} (R(z)f)(0+,t_1,\cdots ,t_n) =& (S(z)f)(t_1,\cdots ,t_n)\\
(R(z)f)(0-,t_1,\cdots ,t_n) = &(S(z)f)(t_1,\cdots ,t_n)+ A_1(R(z)f)(t_1,\cdots ,t_n)\\&+
A_0(S(z)f)(t_1,\cdots ,t_n)\end{align*}
Similar for $\Im z <0$
\begin{align*} (R(z)f)(0+,t_1,\cdots ,t_n) = &(S(z)f)(t_1,\cdots ,t_n)+ A_{-1}^+(R(z)f)(t_1,\cdots ,t_n)\\&+
A_0^+(S(z)f)(t_1,\cdots ,t_n)\\
 (R(z)f)(0-,t_1,\cdots ,t_n) =& (S(z)f)(t_1,\cdots ,t_n)\end{align*}
\end{pf}

\subsection{Characterization of the Hamiltonian}

By Stone's theorem there exists for the unitary strongly continuous one parameter group
a unique spectral family $(E_\lambda)_{\lambda \in \mathbb{R}}$ such that
$$W(t) = \int e^{-it\lambda} dE_\lambda$$
The infinitesimal generator of $W(t)$ is $-iH$ with
$$ H=\int \lambda dE_\lambda $$
and $H$ is the Hamiltonian of $W(t)$. Define the norm
$$\|f\|_H^2 = \int \lambda^2 \langle f|dE_\lambda f \rangle .$$ 
The domain $D_H$ of $H$ is the set all those $f \in \Gamma$, for which $\|f\|_H<\infty$
and $H$ is bounded with repect to that norm. The resolvent 
$$R(z) = \int 1/(z-\lambda) dE_\lambda$$
 defines a bounded
operator from $\Gamma$ onto $D_H$. The equation
\[(z-H)R(z)=1\]
shows that the mapping $R(z)$ is bijective. As $\mathcal{K}$ is dense in $\Gamma$ we have
that $R(z)\mathcal{K}$ is dense in $D_H$ with respect to the $D_H$-norm.We make at first an \emph{Ansatz}$\hat{H}$
for $H$
\begin{defi}
 Assume four operators $ M_0,M_{\pm1},G \in B(\mathfrak{k})$ such that
\begin{align*}M_0^+&=M_0 &M_1^+=& M_{-1}& G^+ &=G ,\end{align*} then define by
$$\hat{H} = i\hat{\partial} + M_1 \hat{\mathfrak{a}}^++M_0 \hat{\mathfrak{a}}^+\hat{\mathfrak{a}}+
M_{-1}\hat{\mathfrak{a}} +G.$$
an application from $D$ into the singular measures on $\mathfrak{R}_0$.
We identify again measures with densities with respect to Lebesgue measure with their densities. \end{defi}
By proposition 5.3.2 we have that $R(z)f \in D$ for $f\in \mathcal{K}$ . The foolowing lemma is a direct
consequence of prop.5.2.3 and the assumptions about the coefficients.

\begin{lem}The  
sesquilinear form
\begin{multline*}f,g \in D \mapsto \langle f|\hat{H}g\rangle \\= \langle f|( i\hat{\partial}+G)g\rangle +
\langle \hat{\mathfrak{a}}f|M_1g\rangle + \langle \hat{\mathfrak{a}}f|\hat{\mathfrak{a}}M_0g\rangle
+ \langle f|\hat{\mathfrak{a}}M_{-1}g \rangle\end{multline*}
exists and is symmetric.
\end{lem}
\begin{prop} Assume $f=Z(f_0+\mathfrak{a}^+f_1)\in D$. Then $\hat{H}f$ is a absolute continuous mesure with
density in $\Gamma$ iff
 $$ -if_1 +M_1 f+M_0\hat{\mathfrak{a}}f=0$$\end{prop}
 \begin{pf} An immediate consequence of the formula 
 $$\hat{H}f=-i(f_0+\hat{a}^+f_1)+ M_1 \hat{\mathfrak{a}}^+f+M_0 \hat{\mathfrak{a}}^+\hat{\mathfrak{a}}f+
M_{-1}\hat{\mathfrak{a}}f +Gf.$$
 \end{pf}
 \begin{defi} Define by $D_0$ the subspace of those functions $ \in D$,
 which obey the condition of the last lemma and denote by $H_0$ the restriction
 of $\hat{H}$ to $D_0$f\end{defi}
\begin{cor} As $\hat{H}$ is symmetric on $D$ ,it is symmetric in $D_0$ too.\end{cor} 
The conditions for the unitarity of the 
operators $\mathcal{O}(u^t_s)(A_i,B))$were (thm.4.5.1) 
that the operators $A_i,i=1,0,-1;B$ fulfill the following conditions:
There exist a unitary operator $\Upsilon$ such that
\begin{align*}
A_0 &= \Upsilon-1\\
A_1&=-\Upsilon A_{-1}^+\\
B+B^+& =-A_1^+A_1 = - A_ {-1}A_{-1}^+.
\end{align*}
\begin{thm}The operator $\hat{H}$ fulfills the equation
\[\hat{H}R(z)f = -f + zR(z)f\tag{$*$}\]
for all $f \in \mathcal{K}$ iff
\begin{align*}A_1 &= \frac{1}{i-M_0/2}M_1\tag{$**$}\\
A_0&=\frac{M_0}{i-M_0/2}\\
A_{-1}&= M_{-1}\frac{1}{i-M_0/2}\\
B&= -iG - \tfrac{i}{2}M_{-1}\frac{1}{i-M_0/2}M_1
\end{align*}
As a consequence
$$\Upsilon = \frac {i+M_0/2}{i-M_0/2}$$
If equaton $(**)$ is fulfilled, then $R(z)$ maps $\mathcal{K}$ into $D_0$.
The domain $D_H$ of the Hamiltonian $H$ of $W(t)$ contains $D_0$ and 
the restriction of $H$ to $D_0$ coincides with the restriction $H_0$  of
$ \hat{H} $
to $D_0$ and $D_0$ is dense in $\Gamma$ and $H$ is the closure of $H_0$.\end{thm}
\begin{pf}
Assume at first $\Im z > 0$. By prop.5.3.3 and 5.3.4
\begin{align*} i\hat{\partial}Rf &= -f -iA_1SF- \hat{\mathfrak{a}}^+(iA_1Rf+iA_0Sf)+(z-iB)Rf\\
 \hat{\mathfrak{a}}Rf&=Sf +\tfrac{1}{2}A_1Rf+\tfrac{1}{2}A_0Sf\end{align*}
 Then
 $$ \hat{H}Rf= -f +zRf+C_1\hat{\mathfrak{a}}^+Rf + C_2\hat{\mathfrak{a}}^+Sf+C_3Sf+C_4Rf$$
 with
 \begin{align*}
 C_1&= -iA_1+M_1+\tfrac{1}{2}M_0A_1\\
 C_2&= -iA_0 + M_0 + \tfrac{1}{2}M_0A_0\\
C_3&= -iA_{-1} + M_{-1} + \tfrac{1}{2}M_{-1}A_0\\
C_4&= -iB+G+ \tfrac{1}{2}M_{-1}A_1
\end{align*}
The equations $C_1=C_2=C_3=C_4=0$ are equivalent to $(**)$.

For $\Im z < 0$ we obtain
\begin{align*}i\hat{\partial}R&= -(f-iA_1^+Sf)+\hat{\mathfrak{a}}^+(iA_{-1}Rf+A_0^+Sf)+(z+iB^+)Rf\\
 \hat{\mathfrak{a}}Rf& = Sf +\tfrac{1}{2}A_{-1}^+Rf+\tfrac{1}{2}A_0^+Sf\end{align*}
 Again
 $$ \hat{H}Rf= -f +zRf+C_1'\hat{\mathfrak{a}}^+Rf + C_2'\hat{\mathfrak{a}}^+Sf+C_3'Sf+C_4'Rf$$
with
\begin{align*}
C_1'&= i A_{-1}^++M_1 + \tfrac{1}{2}M_0A_{-1}^+\\
C_2'&= i A_0^++M_0 + \tfrac{1}{2}M_0A_0^+\\
C_3'&= i A_1^++M_{-1} + \tfrac{1}{2}M_{-1}A_0^+\\
C_4'&= iB^++G + \tfrac{1}{2} M_{-1}A_{-1}^+
\end{align*}
Equations $C_1'=C_2'=C_3'=C_4'=0$ are equvalent to $(**)$ as well, as it should be.

Formula $(*)$ shows, that  $R(z)$ maps $\mathcal{K}$ into $D_0$.

 As $H_0$ coincides
on $R(z)\mathcal{K}   $ with $H$  and $R(z)\mathcal{K}$ is dense in $D_H$ w.r.t. the $D_H$-norm
 , it has a unique extension to $D_H$ coinciding with $H$,
and $H$ is the closure of $H_0$. As $R(z)\mathcal{K}$ is contained in $D_0$, the space
$D_0$ is dense in $\Gamma$.

Consider the matrix elements for 
$\xi \in \mathcal{K}$ and $f \in D_0$
$$ \langle \xi|R(z)H_0f\rangle = \langle R(\overline{z})\xi|H_0f\rangle.$$
Now $ R(\overline{z})\xi$ is in $D_0$ and using the symmetry of $H_0$ the last expression
equals
$$\langle H_0  R(\overline{z})\xi|f \rangle = \langle -\xi + \overline{z} R(\overline{z})\xi
|f \rangle.$$
As $\mathcal{K}$ is dense in $\Gamma$ we obtain
$$R(z)H_0f = -f + zR(z)$$
and with $H_0f=g$
$$\int \frac{1}{z-\lambda}dE_\lambda g = \int \frac{\lambda}{z-\lambda} dE_\lambda f$$
and
$$ \|R(z)g\|^2 = \int \frac{1}{|z-\lambda|^2}\langle g| dE_\lambda g\rangle=
\int \frac{\lambda^2}{|z-\lambda|^2}\langle f| dE_\lambda f\rangle.$$
Put $z=i/\varepsilon$ and multiply by $1/\varepsilon^2$ and obtain
$$ \int \frac{1/\varepsilon^2}{1/\varepsilon^2 + \lambda^2}\langle g| dE_\lambda g\rangle=
\int \frac{\lambda^2/\varepsilon^2}{1/\varepsilon^2 + \lambda^2}\langle f| dE_\lambda f\rangle.$$
Go $\varepsilon \downarrow 0 $ and get
$$\|g\|^2 =\int \langle g|dE_\lambda g \rangle = \int \lambda^2\langle f|dE_\lambda f \rangle = \|f\|_H^2
< \infty .$$
So $f \in D_H$ and 
$$R(z)H_0f = R(z)Hf$$
and 
$$H_0f = Hf.$$
This finishes the proof.
\end{pf}


\begin{thebibliography}{99}




\bibitem{acc89}L.Accardi:Noise and dissipation in quantum theory.Reviews in Math.Physics
2(1990)127-176.
\bibitem{acc96}Accardi L., Lu Y.G., Obata N.:
Towards a non--linear extension of stochastic calculus,
in: Publications of the Research Institute for Mathematical
Sciences, Kyoto, RIMS Kokyuroku 957, N. Obata (ed.) (1996) 1--15
\bibitem{acc99}L.Accardi,Y.-G. Lu, I.V.Volovich: White noise approach to classical 
and quantum stochastic calculus.Preprint 375, Centro Vito Volterra, Universita Roma 2, (1999).
\bibitem{acc2002}L.Accardi,Y.G.Lu,I.V.Volovich: Quantum theory and its
stochastic limit.  Springer 2002.
\bibitem{Ayed2005}W.Ayed:White noise approach to quantum stochastic calulus.Tesis.Tunis 2005
\bibitem{attal92} S. Attal : Problemes d'unicite dans les representations d'operateurs 
sur l'espace de Fock. Seminaire de Probabilites XXVI: LNM 1526, p.617-632 (1992)
\bibitem{attal95} S.Attal: Non-commutative chaotic expansions of Hilbert-Schmidt oprerators on Fockspace. C.M.P. ca 1995
\bibitem {bourbaki1965}N.Bourbaki:Integration chap.1-4,Paris 1965
\bibitem{bourbaki1967}N.Bourbaki:Integration chap.5,Paris 1965
\bibitem{chebotarev}
A.M.Chebotarev,Quantum stochastic differential equation is unitary equivalent 
to a symmetric boundary problem in Fock space. Inf.Dim.Anal.Quantum Prob.
1(1998),175-199.
\bibitem{fagnola}  F.~Fagnola, "On quantum stochastic differential equations with
unbounded coefficients", Prob. Theory and Related Fields, 1990, {\bf
86}, 501--516.
\bibitem{gelfand}
I.Gelfand and N.Vilenkin,\emph{Generalized functions vol.1}.Academic
Press, New York 1964
\bibitem{gough}J.Gough.Causal structures of quantum stochastic integrators.Theoret. and
Math. Physics 111(1997),563-575.
\bibitem{gough97}J.~Gough,  "Noncommutative Ito and Stratonovich noise and stochastic
evolution", Theor. and Math. Phys., 1997, {\bf 113}, N2, 276--284.
\bibitem{gregoratti}
M.Gregoratti,\emph{The hamiltonian associated to some quantum stochastic 
differential equations}.Thesis. Milano 2000.
\bibitem{gregoratti1}M.Gregoratti:The Hamiltonian operator associated with some quantum stochastic evolutions.
Comm.Math.Physics 222(2001),no.1,181-200.
\bibitem{lindsay93} J.M. Lindsay:Quantum and non-causal stochastic calculus. Prob.Theory Relat.Fields 97,p.65-80 (1993).
\bibitem{lindsaymaassen88}J.M.Lindsay, H.Maassen: An integral kernel approach to noise. LNM 1303, p.192-208 (1988)
\bibitem{maasen85} H. Maassen: Quantum Markov processes on Fock space
 described by integral kernels. LNM 1136,p.361-374 (1985)
\bibitem{meyer}
P.A.Meyer \emph{Quantum Probability for Probabilists}
Lecture Notes in Mathematics \textbf{1538}(Springer, Berlin,Heidelberg
1993).
\bibitem{obata94}N.Obata:White noise calculus and Fock space.Springer LNM 1577(1994)
\bibitem{obata08}N.Obata and U.C.Ji:Quantum Stochastic Gradients.2008. To appear.
\bibitem{partha}
K.R.Parthasarathy:
\emph{An Introduction to Quantum Stochastic Calculus}(Birkhaeuser, Basel,Boston,Berlin 1992).
\bibitem{schwartz}L.Schwartz:\emph{Theorie des distributions I}(Herrmann,Paris 1951)
\bibitem{vw96} W.von Waldenfels: Continous Maassen Kernels and the Inverse
Oscillator. Seminaire des Probabilites XXX, Springer LNM
1626(1996)
\bibitem{vw2000}W.von Waldenfels:Continuous kernel processes in quantum probability.Quantum
Probability Communications.Vol.XII,2003,World Scientific,p.237-260.
\bibitem{vw2003}W. von Waldenfels:Description of the damped oscillator by a singular Friedrichs
kernel.In Quantum prob.and rel.fields 2003.
\bibitem{vW2005}W.von Waldenfels:The Hamiltonian of a simple pure number process.Quantum Probability and Infinite
iimensional Analysis.XVIII. World Scientific 2005.
\bibitem{vW2005a}W.von Waldenfels.
Symmetric differentiation and Hamiltonian of a quantum stochastic process. Infinite dimensional analysis
quantun probability ... .2005,p.73-116.
\bibitem{vW2007}W.von Waldenfels.Creation and annihilation operators on locally compact spaces.
 Mathematical analysis of random phenomena.World Scientific. Proceedings Hammamet 2005, p.191-212.
\end{thebibliography}
\end{document}